# Coordinated Lane-Level Variable Speed Limits and Ramp Metering for Successive Weaving Segments Considering Merging/Diverging Risks: A Hybrid Model Predictive Control and Multi-Agent Reinforcement Learning Approach

Guodong Ma[a], Baofeng Sun[a,*], Wenyu Yang[a], Zhihong Yao[b]

a. School of Transportation, Jilin University, Changchun 130022, China;

b. School of Transportation and Logistics, Southwest Jiaotong University, Chengdu 610031, China.

**Abstract**

Successive weaving segments (SWSs) on urban expressways present critical bottlenecks prone to frequent congestion and collisions, demanding fine-grained active traffic management (ATM). However, existing strategies struggle to balance high performance from data-driven adaptive optimization with high resilience and transferability from model-driven rigid protection. To address this, we propose a hybrid model- and data-driven ATM framework that coordinates SWSs. First, we reconstruct a lane-level macroscopic traffic flow model, L-METANET, to capture free and forced lane-changing behaviors. Second, we combine XGBoost-SHAP with random parameters binary logit (RPBL) to fit analytical equations for merging/diverging collision risks, thereby formulating the system cost and reward functions. Finally, we design a hierarchical framework that coordinates lane-level variable speed limits and ramp metering, MPC-STMAPPO, integrating model predictive control (MPC) and multi-agent reinforcement learning (MARL). The upper MPC layer utilizes L-METANET for long-horizon rolling optimization to output baseline commands, while the lower layer employs an ST-MAPPO algorithm incorporating Mamba cells and a graph attention to generate residual actions for short-horizon agile fine-tuning. Real-world experiments on the 18-km Eastern Expressway in Changchun, China, demonstrate that: (i) L-METANET accurately reproduces lane-changing-induced flow redistribution and capacity drops, aligning state evolutions with the ground truth; (ii) XGBoost-SHAP-RPBL model achieves AUC values above 0.80 for most tasks, outperforming traditional logit models; and (iii) Compared to multiple MPC- and MARL-based baselines, MPC-STMAPPO delivers superior training convergence speed and excels across various metrics. Furthermore, in transfer tests under various randomly fluctuating demands, MPC-STMAPPO significantly outperforms pure MARL in generalization capability, revealing high industrial deployment value.

* Corresponding author at: No. 5988 Renmin Street, Changchun 130022, China.
E-mail address: magd22@mails.jlu.edu.cn (Guodong Ma), sunbf@jlu.edu.cn (Baofeng Sun), yangwy24@mails.jlu.edu.cn (Wenyu Yang), zhyao@swjtu.edu.cn (Zhihong Yao).

## 1. Introduction

Urban expressways serve as the backbone of modern urban transportation networks, accommodating ever-increasing traffic demands. However, their operational efficiency and safety often suffer at critical bottleneck nodes, with weaving segments representing one of the most typical problematic zones (Gao et al., 2025; Ouyang et al., 2023; Rim et al., 2023; Zhao et al., 2021). As hubs connecting the mainline with on- and off-ramps, weaving segments facilitate a complex "three-stream weaving" phenomenon, where mainline through traffic, on-ramp merging flows, and off-ramp diverging flows interact at high frequencies (Chen and Ahn, 2018; Golob et al., 2004; Yuan et al., 2024). The resulting forced lane-changing maneuvers easily induce turbulence, causing a sharp drop in traffic capacity accompanied by extremely high collision risks. Fortunately, the advancement of intelligent transportation systems provides new opportunities to alleviate this persistent issue through Active Traffic Management (ATM) based on vehicle-road-cloud integration. Especially with the increasing penetration rate of connected and automated vehicles (CAVs), utilizing precise perception and control capabilities to actively intervene in weaving segments has become a critical means to enhance road network resilience.

Among various ATM strategies, researchers widely recognize ramp metering (RM) and variable speed limits (VSL) as effective tools for mitigating congestion (Ma et al., 2021; Zhang et al., 2025). In this study, ATM specifically refers to either an individual strategy or their combination. Early literature predominantly focused on isolated RM (Airaldi et al., 2025; Peng and Xu, 2023; Zhang et al., 2024a) or standalone VSL control (Jin et al., 2024; Yang et al., 2026; Zhang et al., 2024c). However, because extensive evidence validates the efficacy of their joint operation in alleviating freeway and expressway bottleneck congestion, recent studies increasingly evaluate their coordinated control effects (Han et al., 2025; He et al., 2024; Qiu et al., 2025; Zhang et al., 2025). Most existing coordinated control studies propose segment-level ATM strategies that assume lane homogeneity and apply uniform, coarse-grained control commands across all lanes. This approach works well for basic mainline segments because rare lateral lane-changing maneuvers allow researchers to neglect lane heterogeneity. However, weaving segments severely challenge this assumption. Concurrent forced and free lane-changing behaviors driven by merging and diverging demands impose significant lane-level heterogeneity on weaving traffic. Specifically, early lane-changing friction from exiting vehicles and incoming on-ramp traffic frequently plunges the outer lanes into low-speed congestion, whereas the inner lanes often remain completely smooth. Such a one-size-fits-all segment-level approach lacks differentiated management, thereby wasting inner-lane capacity and failing to guarantee outer-lane safety. Consequently, within a connected and automated environment, we must shift the control granularity from the segment level down to the lane level. Implementing fine-grained lane-level coordinated ATM (Lu et al., 2023; Lu et al., 2024) provides an inevitable solution to counter the lane heterogeneity inherent in weaving segments.

Beyond control granularity, the design of optimization objectives also demands refinement. Existing coordinated control strategies predominantly prioritize traffic efficiency as their primary or sole goal, often incorporating safety considerations insufficiently or unscientifically. Traffic safety in weaving segments warrants dedicated attention due to its unique spatiotemporal evolution characteristics. Specifically, while merging collision risks persist at the entrance, aggressive lane-crossing maneuvers by vehicles exiting the mainline introduce a distinct diverging risk at the exit. Furthermore, the heterogeneous underlying triggers of merging and diverging risks necessitate targeted attention. Neglecting these potential weaving-induced risks or treating them simplistically can easily result in overly aggressive control commands, thereby inducing severe safety hazards.

Translating these complex, lane-level, and multi-objective requirements into real-time control algorithms poses a dual challenge for solution methodologies. Model-driven algorithms, particularly Model

Predictive Control (MPC), find widespread application due to their capacity to manage physical constraints and multi-objective optimization (Chen et al., 2025; Zhang et al., 2024b); however, their performance relies heavily on prediction model accuracy. Within weaving segments, highly stochastic and nonlinear lateral lane-changing behaviors challenge conventional macroscopic traffic flow models like METANET or the cell transmission model, as failing to capture these dynamics accurately induces severe model mismatch (Chen et al., 2021). Conversely, while data-driven methods, notably deep reinforcement learning (DRL), excel at capturing environmental nonlinearities and uncertainties (Afifah and Guo, 2025; Jin et al., 2025; Kang et al., 2024), they suffer from slow training convergence, action oscillations, and a lack of safety boundary guarantees. For safety-critical weaving scenarios, a single approach can rarely satisfy both control stability and high performance. Consequently, developing a hybrid solution framework that merges the complementary advantages of model-and-data-driven approaches represents a critical pathway to breaking current technical bottlenecks (Airaldi et al., 2025; Sun et al., 2024).

To address these research challenges, we propose a lane-level coordinated VSL and RM strategy for weaving segments that explicitly accounts for both merging and diverging collision risks. First, we develop an improved lane-level macroscopic traffic flow model, designated as L-METANET, which explicitly captures the lateral friction and weaving resistance induced by forced lane-changing maneuvers, thereby overcoming the inherent limitations of conventional models in depicting lane heterogeneity. Second, we construct an analytical risk prediction model, XGBoost-SHAP-RPBL, driven dually by controllability and interpretive accuracy, and integrate it into a multi-objective optimization function to achieve proactive safety control. Building upon these components, we design a hybrid model-and-data-driven hierarchical coordinated control framework. Within this framework, the upper layer leverages MPC to calculate baseline LVSL value and RM rates under strict physical constraints, securing system control resilience and guiding the lower-layer learning process. Concurrently, the lower layer utilizes multi-agent reinforcement learning (MARL) to execute residual compensation on these reference policies, ensuring high system performance. While balancing control safety and traffic efficiency, this framework significantly enhances system adaptability to complex traffic turbulence, providing rigorous theoretical guidance and technical support for ATM in expressway weaving segments.

## 2. Literature review and main contributions

This section reviews ATM literature across modeling methodologies and control algorithms. Section 2.1 evaluates the transition to lane-level coordinated VSL and RM strategies for lane heterogeneity, alongside control formulations and closed-loop safety. Section 2.2 contrasts model-driven and data-driven control architectures, delineating their respective limitations to motivate a hybrid alternative. Finally, while Section 2.3 synthesizes research gaps to establish the study's motivations and contributions.

### 2.1. The models of coordinated variable speed limit and ramp metering

As core instruments within ATM frameworks, RM and VSL offer highly effective coordinated approaches to mitigate traffic congestion and collisions. Specifically, RM prevents mainline breakdown by regulating on-ramp inflow rates, whereas VSL smooths mainline speeds to reduce speed differentials between mainline and ramp vehicles, simultaneously throttling upstream traffic volumes entering bottleneck zones to alleviate localized congestion and collision risks. Early engineering practices and theoretical explorations typically decoupled these two strategies into independent subsystems. However, accelerating urbanization and mounting traffic loads have exposed the limitations of this decoupled approach, which often traps the system in local optima as traffic densities approach critical thresholds. To overcome this, the academic community has progressively established an integrated VSL and RM coordinated control paradigm. The core mechanism of this paradigm relies on VSL to preemptively generate low-density spatiotemporal windows upstream, creating sufficient merging space for the ramp traffic released by downstream RM; this coordination significantly suppresses the backward propagation of traffic shockwaves and yields system-wide global optimality. For instance, van de Weg et al. (van de Weg et al., 2019) applied coordinated VSL and RM on a two-lane freeway with two on-ramps and off-ramps, demonstrating superior throughput performance. Similarly, Ma et al. (Ma et al., 2021) confirmed that integrated VSL-RM control outperforms standalone VSL or RM strategies in reducing freeway collision risks and traffic conflicts.

Early RM and VSL strategies relied heavily on model-driven approaches, which necessitate a critical prerequisite: constructing a predictive macroscopic traffic flow model that accurately captures traffic evolution dynamics, enabling controllers to anticipate future traffic states and optimize control actions accordingly. Within the existing research framework, the cell transmission model (CTM) (Daganzo, 1994) and the METANET model (Messmer and Papageorgiou, 1990) represent the two primary model choices. CTM utilizes microscopic or mesoscopic rules to simulate vehicle lane-changing and car-following behaviors across discrete grids, effectively capturing the nonlinear characteristics of traffic flow. However, researchers have noted that such models perform poorly in characterizing capacity drops and stop-and-go traffic waves at freeway bottlenecks (Spiliopoulou et al., 2014). In contrast, the second-order macroscopic METANET model, rooted in fluid dynamics, provides an analytical differential equation form that explicitly describes the spatiotemporal evolution of flow, density, and speed; its low computational complexity establishes it as a preferred predictive tool (Chen et al., 2021; Spiliopoulou et al., 2014).

Notably, despite its advantages, the original METANET model operates on a lane homogeneity assumption and lacks fine-grained lane differentiation (Messmer and Papageorgiou, 1990), meaning that researchers can only apply it to execute segment-level VSL. This uniform management may suffice for conventional mainline segments devoid of frequent ramp operations. However, weaving segments strictly demand fine-grained, lane-level VSL to simultaneously balance traffic efficiency and driving safety, which severely limits the applicability of segment-level models to these areas. This limitation arises because forced lane-changing maneuvers within the weaving section occur concurrently with free lane-changing behaviors upstream (Arman and Tampere, 2022; Wang et al., 2025; Zhou et al., 2024), imposing profound lane-level heterogeneity. Specifically, outer lanes accommodate accelerating merging flows from on-ramps and decelerating diverging flows toward off-ramps, frequently triggering low-speed turbulence, whereas inner lanes predominantly carry high-speed through traffic. Applying uniform control while ignoring such pronounced cross-lane variations inevitably wastes inner-lane capacity or compromises outer-lane safety regulation, significantly degrading coordinated control efficacy. Consequently, recent pioneering studies have begun exploring lane-level coordinated VSL and RM strategies. Researchers generally pursue this objective via two distinct methodological pathways, namely lane-level macroscopic traffic flow prediction models (Chen et al., 2021; Chen et al., 2025; Sai, 2025) and model-free DRL (Lu et al., 2023).

Although lane-level ATM research has recently emerged, existing studies predominantly focus on basic mainline segments or merging zones, leaving weaving segments largely unexplored. Merging zones only involve interactions between two traffic streams, specifically mainline through traffic and incoming ramp traffic, which concentrates conflict points within a limited area. Conversely, weaving segments exhibit a complex "three-stream weaving" structure that integrates mainline through, on-ramp merging, and off-ramp diverging flows, triggering unique “X-shaped” lane-changing conflicts and traffic dynamics fundamentally distinct from merging zones. Consequently, researchers cannot directly transfer control models and algorithms developed for merging areas to weaving segments, necessitating further investigation into fine-grained lane-level management for these complex sections.

The design of control objectives and constraints is equally crucial. Existing literature often prioritizes traffic efficiency over driving safety, a bias that severely constrains the performance of ATM systems under extreme operating conditions. Conventional control strategies typically optimize for minimizing total travel time or maximizing total throughput efficiency; however, this efficiency-centric orientation inadequately accounts for safety, necessitating a more comprehensive control strategy to address pressing traffic safety and efficiency challenges (Zhang et al., 2025). Crucially, the core bottleneck preventing risk integration is not a disregard for safety, but rather the difficulty of constructing a mathematical model that describes the safety characteristics of macroscopic traffic flow in real time, which demands three essential attributes: (1) Accuracy: it must precisely map the relationship between macroscopic aggregated traffic flow indicators and microscopic crash risks; (2) Controllability: the risk metrics must rely on parameters directly or indirectly manageable via VSL and RM, enabling the controller to purposely mitigate the targeted risks; (3) Observability: the input parameters must be rapidly collectable via existing detectors, such as inductive loops, to satisfy real-time computational efficiency requirements; (4) Heterogeneity: Due to variations in traffic demand, geometric length, and channelization design, the underlying macroscopic causal mechanisms of collision risks differ inherently across distinct weaving segments. Consequently, it is essential to develop customized, site-specific risk assessment models tailored to individual weaving locations.

The integration of risk perception into macroscopic ATM is fundamentally impeded by two distinct dimensions, specifically regarding the observation level and the modeling methodology. At the observation level, a significant disconnect persists between microscopic risk dynamics and macroscopic control capabilities. While classic microscopic safety indicators, such as time-to-collision (TTC) and deceleration rate to avoid a crash (DRAC), offer high fidelity in capturing risk, collecting these metrics across an entire network in real time remains unfeasible. Furthermore, macroscopic control methods, including RM and VSL, cannot directly regulate vehicle-level headways or decelerations, causing a severe decoupling between control variables and microscopic optimization objectives. Conversely, traditional macroscopic risk assessment models, despite exhibiting favorable controllability and observability, fail to resolve this bottleneck because they lack the precision required to characterize microscopic collision dynamics and fail to account for the causal heterogeneity across distinct weaving segments. At the modeling methodology level, recent artificial intelligence research has introduced various real-time crash risk evaluation models leveraging machine learning and deep learning techniques (Cheng et al., 2022; Wang et al., 2024; Yang et al., 2021; Zhang et al., 2025). Although these data-driven approaches excel at quantifying real-time collision risks within intelligent transportation systems (Zhang et al., 2025), they introduce a new algorithmic hurdle. Specifically, end-to-end deep learning methods lack analytical, closed-form expressions, making them exceptionally difficult to embed within optimization-based control loops for proactive management. Consequently, resolving these concurrent limitations across both observation levels and modeling methodologies remains an unresolved challenge that warrants further investigation.

### 2.2. The algorithms of coordinated variable speed limit and ramp metering

Translating lane-level, multi-objective ATM into deployable real-time control algorithms poses challenges, particularly in managing the stochastic uncertainties of traffic environments and balancing offline training costs against online computational efficiency. To address these twin challenges, existing control methodologies primarily fall into either model-driven approaches exemplified by MPC (Chen et al., 2025; Han et al., 2021; Mao et al., 2022; Othman et al., 2022; Sirmatel and Yildirimoglu, 2023; Zhang et al., 2024b) or data-driven methods spearheaded by DRL (Afifah and Guo, 2025; Greguric et al., 2022; Jin et al., 2025; Kang et al., 2024; Li and Lasenby, 2024; Wu et al., 2020). Crucially, these two paradigms exhibit highly complementary characteristics.

Regarding tackling environmental uncertainties, MPC and DRL offer contrasting paradigms. MPC has long dominated ATM control due to its capacity to explicitly handle multivariable physical constraints. Its core mechanism leverages a macroscopic traffic flow model to construct a prediction horizon, calculating the optimal control law for the current time step via rolling horizon optimization. Consequently, MPC provides strong physical interpretability and strictly prevents system states from violating safety and physical hard constraints. For instance, Mao et al. (Mao et al., 2022) developed a model-based extended VSL controller integrating an extended CTM with variable-length cells and an MPC scheme optimized via an improved genetic algorithm; this approach reduced total travel time by 14.57% compared to conventional VSL without variable control zones. Because MPC performance heavily depends on prediction model accuracy, researchers continuously strive to develop and calibrate more precise high-order models. These advancements include: (1) incorporating VSL and/or RM rates (Chavoshi et al., 2023; Frejo et al., 2019; Hegyi et al., 2005; Yu and Abdel-Aty, 2014); (2) incorporating mixed traffic flow dynamics (Rahmanidehkordi and Ghasemi, 2024); (3) accounting for lane-changing disruptions (Cheng et al., 2025); and (4) addressing lane-level heterogeneity as detailed in Section 2.1 (see the comprehensive review of Wang et al. (Wang et al., 2022)). These refinements remain vital to elevating predictive accuracy and, consequently, boosting MPC efficacy. Nevertheless, such macroscopic models remain idealized abstractions of reality that struggle to capture all stochastic disturbances inherent in real-world traffic. Because real-world traffic evolution rarely adheres strictly to deterministic differential equations, inevitable model mismatches cause deviations between open-loop predictions and closed-loop executions, hindering closed-loop systems from achieving global optimality. Thus, while model refinement is important, perfectly replicating real-world traffic laws within imperfect environments remains inherently difficult (Li and Lasenby, 2024). To counteract control degradation driven by model bias, researchers have increasingly turned toward model-free DRL methods. DRL approximates state-action mappings through iterative trial-and-error interactions with simulation environments, leveraging the potent nonlinear approximation capabilities of deep neural networks to manage highly nonlinear and stochastic dynamics. Yet, despite demonstrating superior adaptability to environmental uncertainties over traditional control methods, these pure data-driven approaches face some skepticism regarding practical engineering deployment. First, they often lack safety guarantees; unlike MPC, DRL struggles to explicitly embed physical constraints into its optimization process, meaning its exploration-driven learning mechanism can output hazardous commands, which remains unacceptable in safety-critical traffic systems. Second, interpretability remains a major challenge because neural-network-based black-box policies lack physical transparency, hindering trust from traffic authorities and complicating liability attribution after accidents. Finally, poor generalization limits transferability, as RL policies easily overfit their training environments. Minor variations in network topology or traffic flow characteristics often necessitate retraining, stripping the algorithm of the flexible transfer capabilities inherent in MPC.

Beyond their distinct approaches to environmental uncertainty, model-driven and data-driven algorithms exhibit highly complementary computational paradigms, particularly regarding training costs and online computational real-time performance (Sun et al., 2024). The primary advantage of MPC lies in its elimination of offline pre-training. However, this advantage incurs a heavy online computational burden, as the controller treats each optimization step as an entirely new problem, even when encountering identical traffic scenarios previously. Because macroscopic traffic flow models exhibit severe nonlinearities, the dimension of decision variables scales exponentially within multi-lane, long-horizon prediction scenarios, making single-step optimization problems highly time-consuming to solve. Conversely, RL utilizes an offline-training-and-online-inference paradigm. Once trained, the policy requires only a single forward propagation during online deployment, enabling millisecond-level execution that perfectly accommodates the high-frequency demands of real-time control. Nevertheless, RL shifts the computational burden entirely to the offline training stage. Under complex multi-objective and multi-agent configurations, the models frequently suffer from severe convergence bottlenecks. This problem becomes particularly acute in lane-level weaving control scenarios, where high-dimensional state-action spaces combine with environmental non-stationarity. These challenges render MARL highly susceptible to local optima or non-convergent oscillations, making the training process excessively time-consuming and structurally unstable.

In summary, model-driven MPC and data-driven DRL complement each other across multiple operational dimensions. (1) For uncertainty, MPC provides structural robustness but depends heavily on model accuracy, whereas RL offers flexibility but lacks reliable safety boundaries. (2) For computational efficiency, MPC suffers from heavy online computational burdens despite eliminating pre-training requirements, while RL enables rapid online inference at the expense of highly challenging offline training. (3) For transferability, MPC generalizes to any scenario while maintaining high performance given an accurate prediction model, though such perfect models remain elusive; conversely, RL struggles with transfer applications but frequently outperforms MPC within familiar environments. Consequently, isolated methodologies can no longer satisfy the stringent demands for high precision, rigorous safety, and high real-time execution in weaving segment management. Future research paradigms thus favor a deep integration of both approaches (Airaldi et al., 2025; Sun et al., 2024), which yields a hybrid control architecture. This framework leverages MPC to establish physics-based baseline control that satisfies safety constraints, thereby guaranteeing baseline system stability and interpretability. Concurrently, it employs RL to execute real-time data-driven residual compensation on the baseline policy, effectively mitigating model mismatches while elevating online computational efficiency. This integrated philosophy facilitates a profound paradigm shift from passive adaptation to proactive learning.

### 2.3. Research gaps and main contributions

Although extensive literature establishes a solid theoretical foundation, several critical technical gaps persist when addressing the high-risk, high-disturbance scenarios inherent to weaving segments:

(1) **Mismatches in modeling granularity fail to capture pronounced lane-level heterogeneity.** Existing coordinated control studies rely heavily on segment-level METANET models. These frameworks cannot accurately characterize the intense lane heterogeneity arising from forced lane-changing within the weaving section and free lane-changing conflicts upstream, which inevitably induces severe model mismatch.

(2) **Decoupled risk perception and closed-loop control hinder effective collision risk mitigation.** Although classic microscopic risk metrics, including TTC, THW, and DRAC, can precisely quantify risk dynamics, they lack real-time availability, direct controllability, and observability. Conversely, macroscopic risk assessment models fail to establish clear relationships with collisions and struggle to capture the causal heterogeneity of risk across different weaving segments. Consequently, these limitations impede the construction of a complete closed-loop controller that bridges risk perception and control execution.

(3) **Isolated model-driven or data-driven architectures show inadequate adaptability to uncertainties, transfer deployment, and real-time computation.** Model-driven MPC suffers from model mismatch and heavy online computational burdens, failing to flexibly respond to stochastic disturbances and high-frequency execution demands. Conversely, data-driven DRL struggles with slow training convergence and poor generalization in transfer scenarios, which compromises training efficiency and fails to safeguard the strict physical performance lower bound of the system.

To bridge these technical gaps, our study delivers the following main contributions:

(1) We develop a multi-lane macroscopic traffic flow model, L-METANET, tailored for weaving segments by incorporating free and forced lane-changing behavior. L-METANET explicitly reproduces lane-changing-driven cross-lane flow redistribution and capacity drops, serving as a high-fidelity predictive model for fine-grained lane-level ATM.
(2) We developed a multi-objective coordinated LVSL and RM control strategy that explicitly accounts for merging and diverging collision risks within weaving segments. By integrating XGBoost-SHAP feature selection with a random parameters binary logit model, we designed an analytical risk assessment model to ensure direct controllability and real-time data availability while effectively capturing the heterogeneity of macroscopic risk precursors across different weaving segments. This closed-form risk formulation is directly embedded into the MPC cost function and MARL reward function.
(3) We design a hybrid model-and-data-driven hierarchical coordinated control ATM framework. The upper layer implements MPC alongside L-METANET to compute baseline control commands under strict physical boundaries, guaranteeing system stability and guiding lower-layer training. Concurrently, the lower layer deploys a spatiotemporal attention-reinforced MAPPO algorithm, embedding Mamba blocks and graph attention mechanisms into the actor and value networks to capture long-term temporal dependencies and spatial agent interactions. By incorporating the upper-layer MPC control outputs as environmental states and reference control commands, the lower layer executes residual compensation over a short control horizon, ensuring high learning efficiency and robust disturbance rejection capabilities.

The remainder of this paper is organized as follows. Section 3 outlines the problem formulation and the overall research framework. Section 4 presents the lane-level macroscopic traffic flow model and the analytical risk assessment model tailored for expressway weaving segments. Section 5 constructs the MPC and MARL-based hierarchical coordinated control framework integrating lane-level VSL and RM. Section 6 conducts extensive simulation experiments to validate the performance of the proposed macroscopic traffic flow model, macroscopic risk assessment model for weaving segments, and hierarchical coordinated ATM framework. Finally, Section 0 concludes the paper.

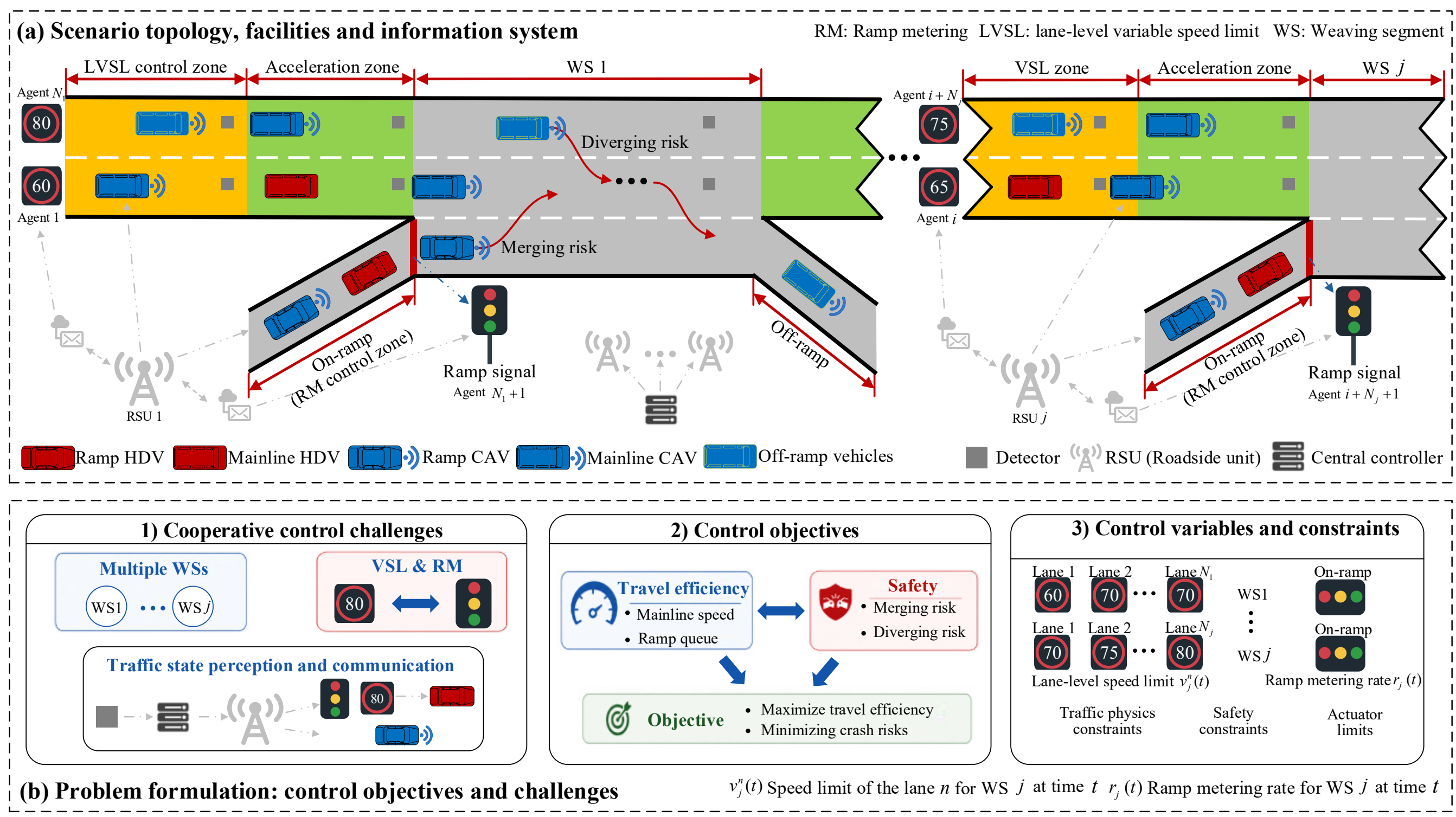

Fig. 1 Scenario definition and problem formulation

## 3. Problem formulation and research framework

### 3.1. Scenario definition and problem description

This study addresses the macroscopic coordinated control problem across continuous multi-weaving segments on urban expressways (Fig. 1). As illustrated in Fig. 1(a), each core control unit is partitioned longitudinally into four consecutive zones: an upstream LVSL zone, an acceleration zone, a weaving bottleneck segment, and a downstream off-ramp diverging segment. The infrastructure architecture comprises two primary layers: (1) Perception layer: Densely deployed detectors gather real-time, lane-specific macroscopic traffic states, boundary demands at origin cross-sections and on-ramps, and queue lengths at on-ramps; (2) Control layer: A central controller manages integrated mainline LVSL and on-ramp RM strategies across fine-grained lane-segment units. Control commands are dispatched via vehicle-to-everything (V2X) communications to connected and automated vehicles (CAVs), and via overhead dynamic speed displays and conventional ramp signals to human-driven vehicles (HDVs).

Under a mixed traffic environment (Fig. 1(b)), executing this coordinated control framework entails three multidimensional challenges: (1) Competing control objectives: maximizing throughput efficiency (maintaining high speeds and minimal queues) inherently conflicts with mitigating traffic safety risks. Intensive lateral cutting-in maneuvers within weaving sections severely compress safety gaps, accelerating collision risks. Balancing traffic efficiency against merging/diverging risks remains a primary trade-off; (2) Cascading congestion propagation: due to the dense spacing of urban expressway ramps, downstream bottleneck breakdowns generate backward-propagating stop-and-go traffic waves. These waves trigger cascading failures across successive upstream weaving segments, rendering conventional, isolated bottleneck management sub-optimal; (3) Spatiotemporal variable coupling: RM physically restricts inflows to protect the mainline but risks urban surface-street spillback, whereas VSL regulates mainline speeds to proactively create upstream low-density merging windows and smooth longitudinal speed differentials. Because these control variables are highly coupled, a lack of joint, fine-grained design triggers excessive mainline delays or ramp deadlocks.

To address these coordinated challenges, the control variables, optimization objectives, and constraints of this study are synthesized into the following core elements (as illustrated in Fig. 1(b)):

(1) Control variables: For any given control time step $t$ ($t$ is MARL time step) and weaving segment $j\left(j \in \{1,2,\ldots,J\}\right)$, the joint control action sequence of the system comprises two primary components: LVSL value $v_j^n(t)$: The control commands applied to lane $n$ within the upstream LVSL zone of weaving segment $j$. RM rate $r_j(t)$: The regulation rate command applied to the on-ramp signaling lights of weaving segment $j$.

(2) Optimization objectives: The central controller pursues multi-objective Pareto optimality to simultaneously optimize system-level efficiency and safety performance: Maximizing system efficiency: This objective seeks to maximize the average mainline travel speed while minimizing the total vehicle queueing time at the on-ramps; Minimizing collision risks: This objective aims to suppress macroscopic merging and diverging risks within the weaving sections, thereby actively reshaping a safer traffic fluid state.

(3) Constraints: The solution of the control sequence must strictly comply with rigorous physical and safety boundaries: 1)Traffic physics constraints: The spatiotemporal evolution of traffic states must strictly adhere to the L-METANET equations governing the conservation of mass and momentum relaxation; 2)Safety boundary constraints: These incorporate the maximum allowable queue length capacity at individual on-ramps alongside spatial variation bounds for both VSL and ramp signaling across adjacent lanes or segments; 3)Actuator physical limits: Restricted by traffic regulations and mechanical characteristics, these impose maximum rate-of-change thresholds as well as absolute upper and lower bounds for both VSL and RM rates.

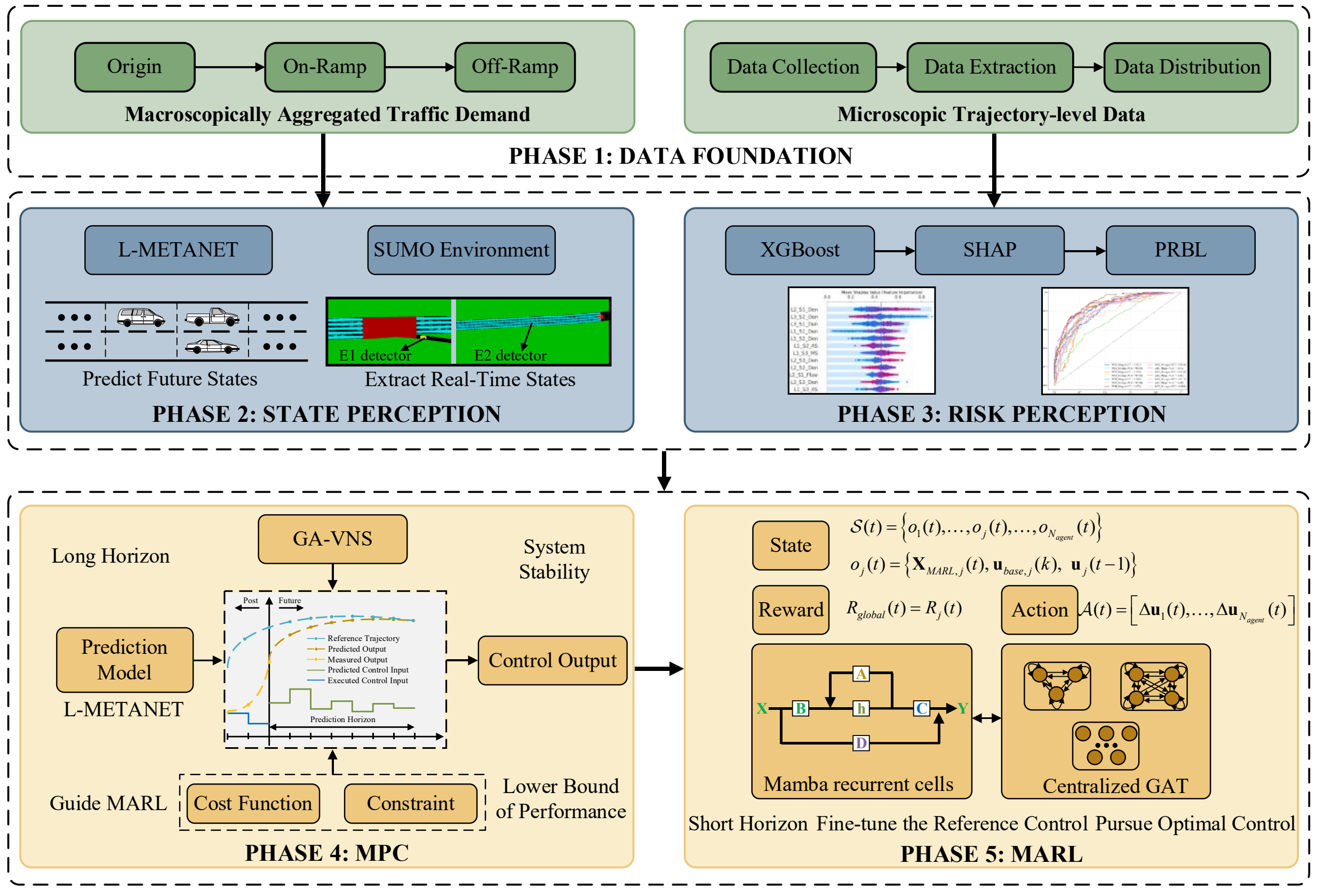

Fig. 2 Overall framework of this study

3.2. Overall research framework

To address the aforementioned challenges in successive weaving segments coordinated control, the ATM framework decouples the control logic into distinct, interconnected operational layers spanning data foundation, state perception, risk assessment, upper-layer baseline control, and lower-layer residual compensation, as illustrated in Fig. 2.

(1) Data foundation layer: This layer is categorized into macroscopic traffic demand data and microscopic vehicle trajectory data. The former calibrates the L-METANET model and defines boundary traffic demands for environmental evolution in the SUMO simulation. In contrast, the latter facilitates the fitting and training of the risk assessment model.

(2) State perception layer: Serving as the data and model bedrock of the entire control framework, this layer leverages the L-METANET model alongside the underlying SUMO simulation environment to deliver high-precision state perception and traffic predictions for both MPC and MARL operations.

(3) Risk assessment layer: By integrating the XGBoost algorithm, SHAP interpretable analysis, and the random parameters binary logit (RPBL) model, this layer constructs an analytical risk evaluation and prediction framework. It simultaneously guarantees high risk-quantification accuracy and operational availability for control execution, quantifying dynamic collision risks across the network in real time to serve as the cost function for upper-layer MPC and the reward function for lower-layer MARL.

(4) Upper-layer baseline control layer (MPC-based): Operating at a longer macro-control period to capture long-term traffic evolution trends, this layer utilizes L-METANET as its predictive model. Under strict physical constraints, it computes a sequence of safe baseline control outputs via rolling horizon optimization, establishing a robust safety floor for system operations.

(5) Lower-layer residual compensation control layer (MARL-based): Executed at a shorter micro-control step, this layer compensates for the sluggish response of MPC to microscopic disturbances and mitigates the impacts of predictive model mismatches. Distributed agents deployed at individual control nodes execute rapid inference based on real-time local states using an enhanced ST-MAPPO algorithm, outputting high-frequency residual control actions.

Ultimately, grounded in real-time state perception and risk assessment, the system dynamically superimposes the upper-layer MPC baseline commands with the lower-layer MARL residual adjustments. The combined signals undergo a safety boundary clipping process to comply with physical and regulatory constraints before being dispatched to the traffic actuators.

## 4. Macroscopic traffic flow and risk assessment models for expressway weaving segments

### 4.1. L-METANET: a lane-level macroscopic traffic flow model

This section extends the classic METANET model across multiple dimensions, most notably by expanding its granularity to the lane level, thereby establishing a generalized methodology for traffic flow modeling on urban expressway networks. Building upon this formulation, we also present a rigorous parameter calibration approach.

#### *4.1.1 Formulation of L-METANET*

The conventional METANET model (Messmer and Papageorgiou, 1990) treats multi-lane cross-sections as homogeneous entities, neglecting critical speed and density differentials between inner fast and outer slow lanes. This homogeneity assumption fails on urban expressways, lacking the granular state variables required for lane-level ATM. Furthermore, original METANET cannot characterize localized ramp disturbances on the outermost lane or the lateral lane-changing friction inherent to weaving segments. To bridge these gaps, this section downscales traffic flow variables to the lane level and introduces five critical modifications to the classic METANET framework, establishing a generalized macroscopic modeling methodology tailored for successive weaving segments:

(1) Expansion to the lane level: Discretizing cross-sectional traffic states into distinct, lane-specific variables to track fine-grained flow dynamics;
(2) Incorporation of asymmetric speed impacts: Modeling the directional, unbalanced speed degradation caused by lateral merging and diverging behaviors;
(3) Consideration of VSL control and compliance of HDVs: Incorporating the explicit impacts of VSL combined with HDV driver compliance rates;
(4) Accounting for upstream and downstream state interdependencies triggered by changes in lane drops or additions;
(5) Dynamic flow estimation of lateral behaviors: Quantifying the localized, cross-lane flow exchanges driven concurrently by free and forced lane-changing maneuvers.

**Improvement 1:** Expansion to the lane level

To accurately characterize lane-segment-level traffic state evolution and facilitate lane-specific VSL, the corridor is discretized spatially into $M$ segments of length $\Delta x_m$, where $\Delta x_m$ denotes the length of segment $m$. This partition establishes the mathematical foundation for subsequent variable-length segment modeling. Time is discretized into steps of duration $\Delta t$. For each traffic state variable, a lane index $n\left(n=1,2,\cdots,N_m\right)$ is introduced, where $N_m$ denotes the outermost lane and represents the total number of lanes in segment $m$. The schematic layout of the segment-lane discretization and aggregated traffic flow variables is illustrated in Fig. 3. The lane-level macroscopic traffic flow model extended state evolution is governed by Eqs. (1) to (5).

$$k_{m,n}^{i+1} = k_{m,n}^{i} + \frac{\Delta t}{\Delta x_m}\left(\hat{q}_{m,n}^{i} + r_{m,n}^{i} - q_{m,n}^{i}\right) \tag{1}$$

$$\hat{q}_{m,n}^{i} = \sum_{n' \in N^{m,n}} \left(q_{m-1,\tilde{n}}^{i} + \phi_{m-1,\tilde{n}}^{i,fr} + \phi_{m-1,\tilde{n}}^{i,fo}\right) \cdot R_{m-1,\tilde{n}\to n}^{i} \tag{2}$$

$$v_{m,n}^{i+1} = v_{m,n}^{i} + \frac{\Delta t}{\tau}[V(k_{m,n}^{i}) - v_{m,n}^{i}] + \frac{\Delta t}{\Delta x_m} \cdot v_{m,n}^{i} \cdot (v_{m-1,\tilde{n}}^{i} - v_{m,n}^{i}) - \frac{\eta \cdot \Delta t \cdot \left(k_{m+1,\tilde{n}}^{i} - k_{m,n}^{i}\right)}{\tau \cdot \Delta x_m \cdot \left(k_{m,n}^{i} + \kappa\right)} \tag{3}$$

$$V(k_{m,n}^{i}) = v_{f,m,n} \cdot \exp\left[-\frac{1}{a_{m,n}}\left(\frac{k_{m,n}^{i}}{k_{cr,m,n}}\right)^{a_{m,n}}\right] \tag{4}$$

$$q_{m,n}^{i} = k_{m,n}^{i} \cdot v_{m,n}^{i} \tag{5}$$

where Eq. (1) is the traffic flow conservation equation: $k_{m,n}^{i}$, $q_{m,n}^{i}$ respectively represent the density and flow rate of lane $n$ in segment $m$ at time step $i$; $r_{m,n}^{i}$ and $s_{m,n}^{i}$ respectively represent the inflow merging from the on-ramp and the outflow diverging to the off-ramp connected to lane $n$ of segment $m$; $\hat{q}_{m,n}^{i}$ is the longitudinal flow entering lane $n$ of segment $m$ from all lanes in segment $m-1$ connected to lane $n$ of segment $m$ at time step $i$, with its calculation formula given in Eq. (2); $\tilde{n}$ is the lane number of the upstream segment connected to lane $n$. Generally, when there is no increase or decrease of lanes on the mainline, $\tilde{n} = n$, but this mapping relationship will change when on/off-ramps or lane changes exist, which should strictly depend on the actual road connection layout; $\phi_{m-1,\tilde{n}}^{i,fr}$ represents the flow rate variation in lane $\tilde{n}$ of segment $m-1$ caused by free lane-changing behaviors, $\phi_{m-1,\tilde{n}}^{i,fr} = \phi_{m-1,\tilde{n}-1\to\tilde{n}}^{i,fr} + \phi_{m-1,\tilde{n}+1\to\tilde{n}}^{i,fr} - \phi_{m-1,\tilde{n}\to\tilde{n}-1}^{i,fr} - \phi_{m-1,\tilde{n}\to\tilde{n}+1}^{i,fr}$; $\phi_{m-1,n'}^{i,fo}$ represents the flow rate variation in lane $\tilde{n}$ of segment $m-1$ caused by forced lane-changing behaviors, $\phi_{m-1,\tilde{n}}^{i,fo} = \phi_{m,\tilde{n}+1\to\tilde{n}}^{i,fo} + \phi_{m,\tilde{n}-1\to\tilde{n}}^{i,fo} - \phi_{m,\tilde{n}\to\tilde{n}+1}^{i,fo} - \phi_{m,\tilde{n}\to\tilde{n}-1}^{i,fo} - s_{m-1,\tilde{n}}^{i}$; $R_{m-1,\tilde{n}\to n}^{i}$ is the proportion of traffic flowing from lane $\tilde{n}$ to lane $n$ within segment $m-1$; $N^{m,n}$ is the set of all lane numbers in segment $m-1$ connected to lane $n$ of segment $m$. Eq. (3) is the dynamic speed equation: the first term represents the average speed of lane $n$ in segment $m$ at time step $i$; the second term is the relaxation term, representing the difference between the desired speed and the current average speed, which reflects the tendency of drivers to travel at the desired speed during driving; the third term is the convection term, representing the speed variation caused by vehicles entering the current cell $m$ from the upstream cell $m-1$, reflecting the impact of upstream speed fluctuations; the fourth term is the anticipation term, which adjusts the speed of the current cell $m$ through feedback on the density variation of the downstream cell $m+1$, reflecting the impact of downstream density changes, where $\tau$ is the driver reaction time, $\eta$ is the speed-density relationship coefficient, and $\kappa$ is the elasticity coefficient. Eq. (4) is the steady-state speed equation, and Eq. (5) represents the three-parameter relationship equation: $v_{f,m,n}$, $k_{cr,m,n}$, and $a_{m,n}$ are the free-flow speed, critical density, and the model parameter controlling the shape of the fundamental diagram for lane $n$ in segment $m$, respectively.

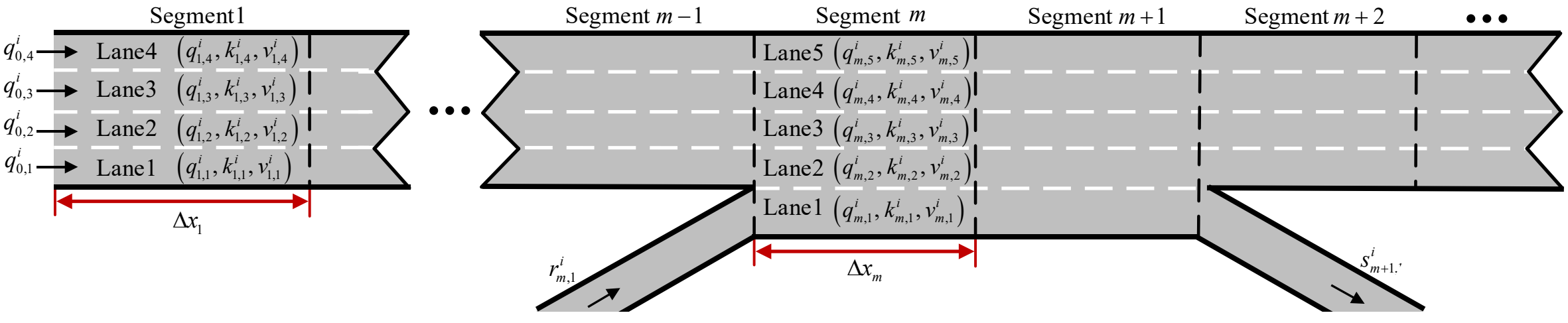

Fig. 3 Schematic diagram of segment-lane division and aggregated traffic flow variables

Subsequently, it is also necessary to model the flow equations for the origin segments. These origin segments (including on-ramps and mainline entries) serve to receive and discharge network traffic demands, where the actual entering flow rate depends on the minimum of the desired inflow (demand) and the maximum downstream acceptable capacity (supply). The relationship between traffic demand and the entering flow rate within the L-METANET macroscopic traffic flow model is captured via a queue model, with the specific calculation procedures formulated in Eqs. (6) to (11).

$$q_0^i = z_m^i \hat{q}_0^i \tag{6}$$

$$z_m^i = \begin{cases} 0, 0 \; is \; mainline \\ R_m^i, 0 \; is \; ramp \\ 1, 0 \; is \; main \; entrance \end{cases} \tag{7}$$

$$\hat{q}_0^i = \min\left(q_0^{i,1}, q_0^{i,2}\right) \tag{8}$$

$$q_0^{i,1} = d_0^i + \frac{w_0^i}{\Delta t} \tag{9}$$

$$q_0^{i,2} = Q_0^{cap,i} \times \min\left(1, \frac{k_{jam} - k_1^i}{k_{jam} - k_{cr}}\right) \tag{10}$$

$$w_0^{i+1} = w_0^i + \Delta t \cdot \left(d_0^i - q_0^i\right) \tag{11}$$

where Eqs. (6) to (11) are the flow rate equations for the origin segment ($q_0^i$), which depend on the minimum value between the traffic demand of the origin segment $q_0^{i,1}$ and the maximum acceptable flow rate of the mainline segment $q_0^{i,2}$. It should be noted that when the origin segment is the mainline origin segment, $q_0^i = \hat{q}_0^i$; when the origin segment is an on-ramp entrance, it is also influenced by the RM rate $R_m^i \in (0,1)$. Eq. (9) is the calculation formula for $q_0^{i,1}$, which is determined by the average arrival flow of the origin segment $d_0^i$ and the queue length of the origin segment $w_0^i$ per unit time. Eq. (10) is the calculation formula for $q_0^{i,2}$, which is determined by the mainline capacity of the origin segment $Q_0^{cap,i}$ and the current state of the road ($k_{cr}$, $k_1^i$, and the jam density $k_{jam}$). Eq. (11) is the calculation formula for $w_0^{i+1}$, which is determined by $w_0^i$, $d_0^i$, and $q_0^i$.

**Improvement 2:** Asymmetric speed impacts from lateral merging and diverging behaviors.

In the METANET model, the impacts of on-ramp merging, off-ramp diverging, and weaving behaviors on speed are neglected. When the modeling granularity downscales to the lane level, the dynamic boundaries of traffic flow undergo a fundamental shift: for any specific lane $n$, the physical essence of vehicles merging from an on-ramp into the mainline is approximately equivalent to vehicles changing from an adjacent lane into this lane, both of which manifest as the insertion of external vehicles into the subject lane. Such forced insertions into the gaps ahead of the target lane disrupt the original steady-state car-following behavior. To re-establish safe spacing, following vehicles in the target lane must execute intense braking maneuvers, thereby leading to a sharp degradation in the macroscopic average speed of that lane. Conversely, vehicles departing from the current lane to an off-ramp or changing to an adjacent lane both manifest as the extraction of vehicles from the current lane, which also exerts a certain impact on the average speed of that lane; naturally, this impact is far less severe than that caused by vehicle insertions, necessitating separate considerations. Consequently, we define a generalized inflow $\phi_{m,n}^{i,in}$ and a generalized outflow $\phi_{m,n}^{i,out}$, which respectively represent the number of entering and departing vehicles accepted by lane $n$ of segment $m$ within time step $i$, thereby better quantifying the speed impacts triggered by lane-changing behaviors.

As mentioned, lateral lane-changing behaviors are categorized into free lane-changing and forced lane-changing: (1) Free lane-changing represents the lateral migrations across lanes executed by vehicles to achieve higher operational efficiency. Let $\phi_{m,n+1\to n}^{i,fr}$ denote the traffic flow transferring from lane $n-1$ to lane $n$ within segment $m$ due to free lane-changing behaviors at time step $i$; (2) Forced lane-changing refers to the compulsory merging and diverging maneuvers dictated by on-ramp and off-ramp operations, which comprise two components. The first component is the on-ramp inflow $r_{m,n}^{i}$, representing the traffic flow entering lane $n$ of segment $m$ from the adjacent on-ramp under a one-to-one mapping. The second component is the traffic flow forcing its way lane-by-lane from various inner lanes to depart the mainline via the off-ramp. Specifically, when $s_{m,n}^{i} < q_{m,n}^{i}$, $\phi_{m,n+1\to n}^{i,fo} = 0$; when $s_{m,n}^{i} > q_{m,n}^{i}$, vehicles in the more inner lanes must first laterally change into lane $n$ before diverging from lane $n$ to exit the mainline, where $\phi_{m,n+1\to n}^{i,fo}$ represents this specific flow component. It should be noted that under extreme conditions, $s_{m,n}^{i}$ can become exceptionally large, requiring vehicles across multiple inner lanes to execute successive lateral lane changes to satisfy off-ramp diverging demands. Grounded in these analyses, the calculation methods for $\phi_{m,n}^{i,in}$ and $\phi_{m,n}^{i,out}$ are formulated in Eqs. (12) to (13). Within these two equations, the first two terms represent the traffic flow induced by free lane-changing, whereas the latter two terms represent the flow driven by forced lane-changing behaviors.

$$\phi_{m,n}^{i,in} = \phi_{m,n-1\to n}^{i,fr} + \phi_{m,n+1\to n}^{i,fr} + r_{m,n}^{i} + \phi_{m,n+1\to n}^{i,fo} + \phi_{m,n-1\to n}^{i,fo} \tag{12}$$

$$\phi_{m,n}^{i,out} = \phi_{m,n\to n-1}^{i,fr} + \phi_{m,n\to n+1}^{i,fr} + s_{m,n}^{i} + \phi_{m,n\to n-1}^{i,fo} + \phi_{m,n\to n+1}^{i,fo} \tag{13}$$

$$\begin{aligned} v_{m,n}^{i+1} = {} & v_{m,n}^{i} + \frac{\Delta t}{\tau}\left[V\left(k_{m,n}^{i}\right) - v_{m,n}^{i}\right] + \frac{\Delta t}{\Delta x_m}\cdot v_{m,n}^{i}\cdot\left(v_{m-1,\tilde{n}}^{i} - v_{m,n}^{i}\right) - \frac{\eta\cdot\Delta t\cdot\left(k_{m+1,\tilde{n}}^{i} - k_{m,n}^{i}\right)}{\tau\cdot\Delta x_m\cdot\left(k_{m,n}^{i}+\kappa\right)} \\ & - \frac{\delta_{in}\cdot\Delta t\cdot v_{m,n}^{i}\cdot\phi_{m,n}^{i,in}}{\Delta x_m\cdot\left(k_{m,n}^{i}+\kappa\right)} - \frac{\delta_{out}\cdot\Delta t\cdot v_{m,n}^{i}\cdot\phi_{m,n}^{i,out}}{\Delta x_m\cdot\left(k_{m,n}^{i}+\kappa\right)} \end{aligned} \tag{14}$$

Finally, we construct the dynamic speed equation that accounts for the asymmetric speed impacts of lateral merging and diverging behaviors. The modified lane-level dynamic speed equation is formulated in Eq. (14). Because the impacts of vehicle insertion and extraction on speed differ significantly, with insertion exerting a substantially greater impact than extraction, the parameter calibration must strictly satisfy $\delta_{in} > \delta_{out}$. This mathematically guarantees the asymmetric traffic flow characteristic wherein deceleration occurs more readily than speed recovery within weaving segments.

**Improvement 3:** The impacts of VSL and HDV driver compliance rates.

In a connected and automated environment implementing lane-level ATM, the LVSL commands dynamically displayed on the roadside (or directly received by RSU) artificially truncate the desired speed of vehicles in a free-flow state. However, the actual traffic stream on the roadway is a mixed traffic flow composed of CAVs and HDVs. There is a fundamental difference in the response mechanisms of these two classes of vehicles to management and control commands: CAVs can obtain speed limit commands in real time via V2X communications and strictly execute them with 100% compliance, whereas HDVs exhibit only partial compliance with speed limit commands due to individual driver heterogeneity, perception delays, and subjective driving intentions. To scientifically quantify the comprehensive impact of these heterogeneous compliance characteristics on macroscopic traffic flow, this model introduces the CAV penetration rate parameter $\alpha \in [0,1]$ and the average speed limit compliance rate pa-rameter of HDVs $\beta \in (0,1)$. Accordingly, a comprehensive speed limit compliance index for mixed traffic flow $\beta_{mix}$ is designed, as shown in Eq. (15). Based on this comprehensive index that governs the variations in desired speed caused by VSL, the steady-state speed fundamental diagram equation is modified to take the minimum value between the speed limit and the spontaneous desired speed, thereby accounting for the heterogeneous characteristics of mixed traffic flow, as shown in Eq. (16).

$$\beta_{mix} = \alpha \cdot 1 + (1-\alpha) \cdot \beta \tag{15}$$

$$V(k_{m,n}^{i}) = \min\left\{ \beta_{mix} \cdot v_{m,n}^{VSL,i} + (1-\beta_{mix}) \cdot v_{f,m,n} \cdot \exp\left[ -\frac{1}{a_{m,n}} \left( \frac{k_{m,n}^{i}}{k_{cr,m,n}} \right)^{a_{m,n}} \right], v_{f,m,n} \cdot \exp\left[ -\frac{1}{a_{m,n}} \left( \frac{k_{m,n}^{i}}{k_{cr,m,n}} \right)^{a_{m,n}} \right] \right\} \tag{16}$$

It should be noted that the application of VSL fundamentally alters the shape of the traffic fundamental diagram, causing the capacity of the road segment to evolve dynamically. According to traffic flow theory, the steady-state flow-density relationship function $q_{m,n}^{VSL,i}(k)$ of a lane unit equals the product of density and the desired speed. Substituting Eq. (16) fully into this relationship yields a piecewise flow function, as shown in Eq. (17). The dynamic capacity of the roadway corresponds to the global maximum of this flow-density function, and the density at this extreme point defines the dynamic critical density $k_{cr,m,n}^{VSL,i}$. Since this requires calculating the derivative of the min function, solving for this extreme point necessitates differentiating the flow with respect to density in Eq. (17) and setting the derivative to zero (i.e., $\mathrm{d}q_{m,n}^{VSL,i}(k) / \mathrm{d}k = 0$), thereby obtaining the stationary points of the min function. According to the mathematical theory of optimization for piecewise continuous functions, this global maximum point must originate from the following set of three candidate critical points $\mho = \{k_A, k_B, k_C\}$:

(1) The stationary point of the speed-limit branch $k_A$: Differentiating $q_{m,n,1}^{VSL,i}(k)$ and setting $\mathrm{d}q_{m,n,1}^{VSL,i}(k) / \mathrm{d}k = 0$ yields a transcendental equation. Although this transcendental equation lacks an analytical solution, $k_A$ can be obtained by solving the equation via numerical algorithms.

(2) The stationary point of the natural branch $k_B$: Differentiating $q_{m,n,2}^{VSL,i}(k)$ and setting $\mathrm{d}q_{m,n,2}^{VSL,i}(k) / \mathrm{d}k = 0$. Because $q_{m,n,2}^{VSL,i}(k)$ represents the uncontrolled fundamental diagram function, its stationary point solution remains identically equal to the original natural critical density, i.e., $k_B = k_{cr,m,n}$.

(3) The intersection point of the two branches $k_C$: Equating the two flow curves, i.e., $q_{m,n,1}^{VSL,i}(k) = q_{m,n,2}^{VSL,i}(k)$, yields the non-differentiable turning point $k_C$.

To determine the final global maximum point, a physical validity check must be performed on the points within the candidate set $\mho$ to eliminate spurious peaks that are truncated by the envelope and cannot be achieved in real-world traffic flow. Valid density points must strictly comply with the definition of the lower envelope: if $k_A$ satisfies $q_{m,n,1}^{VSL,i}(k_A) \le q_{m,n,2}^{VSL,i}(k_A)$, then $k_A$ is valid; if $k_B$ satisfies $q_{m,n,2}^{VSL,i}(k_B) \le q_{m,n,1}^{VSL,i}(k_B)$, then $k_B$ is valid; the intersection point $k_C$ is inherently valid. Let $\mho_{valid}$ denote the set of valid candidate points. Among all valid candidate points that pass the verification, the flow values corresponding to all candidates are extracted; the maximum flow value represents the traffic capacity under the influence of LVSL and mixed fleet characteristics $Q_{m,n}^{cap,i}$, as shown in Eq. (19). Its corresponding valid candidate point defines the dynamic critical density $k_{cr,m,n}^{VSL,i}$.

$$
\begin{aligned}
&q_{m,n}^{VSL,i}(k) = V(k_{m,n}^{i}) \cdot k \\
&= k \cdot \min\left\{ \beta_{mix} \cdot v_{m,n}^{VSL,i} + \left(1-\beta_{mix}\right) \cdot v_{f,m,n} \cdot \exp\left[ -\frac{1}{a_{m,n}} \left( \frac{k}{k_{cr,m,n}} \right)^{a_{m,n}} \right],\ v_{f,m,n} \cdot \exp\left[ -\frac{1}{a_{m,n}} \left( \frac{k}{k_{cr,m,n}} \right)^{a_{m,n}} \right] \right\} \\
&= \begin{cases} q_{m,n,1}^{VSL,i} = (k)k \cdot \beta_{mix} \cdot v_{m,n}^{VSL,i} + \left(1-\beta_{mix}\right) \cdot v_{f,m,n} \cdot \exp\left[ -\frac{1}{a_{m,n}} \left( \frac{k}{k_{cr,m,n}} \right)^{a_{m,n}} \right], & q_{m,n,1}^{VSL,i}(k) \le q_{m,n,2}^{VSL,i}(k) \\ q_{m,n,2}^{VSL,i} = (k)k \cdot v_{f,m,n} \cdot \exp\left[ -\frac{1}{a_{m,n}} \left( \frac{k}{k_{cr,m,n}} \right)^{a_{m,n}} \right], & q_{m,n,1}^{VSL,i}(k) > q_{m,n,2}^{VSL,i}(k) \end{cases}
\end{aligned} \tag{17}
$$

$$
\beta_{mix} \cdot v_{m,n}^{VSL,i} + (1-\beta_{mix}) \cdot V_{f,m,n} \cdot \exp\left[ -\frac{1}{a_{m,n}} \left( \frac{k}{k_{cr,m,n}} \right)^{a_{m,n}} \right] \cdot \left[ 1 - \left( \frac{k}{k_{cr,m,n}} \right)^{a_{m,n}} \right] = 0 \tag{18}
$$

$$
Q_{m,n}^{cap,i} = \max_{k \in \mho_{valid}} \left\{ \min\left[ q_1(k), q_2(k) \right] \right\} \tag{19}
$$

**Improvement 4:** Upstream and downstream state interdependencies under cross-sectional geometric variations.

In METANET, when the speed of the upstream segment $m-1$ is greater than that of the current segment $m$, the convection term is positive. Without considering other influencing factors, the speed of the current segment in the next sampling period will increase. However, this is an idealized condition and is not fully applicable to realistic, complex weaving segment scenarios: if the upstream segment $m-1$ is in a free-flow state, the downstream segment $m+1$ is in a congested state, and the current segment $m$ is in an unstable car-following state. If the upstream speed is greater than the current segment speed at this moment, the original model would conclude that the speed of segment $m$ will increase in the next period; however, empirical traffic flow experience indicates that vehicles in the current segment will instead decelerate due to the backward propagation effect of downstream congestion waves. This contradiction demonstrates that the convection term in the original equation is primarily applicable to situations where both upstream and downstream segments are in free-flow states. Grounded in this, a geometric mean is adopted to improve the convection term. Combining the dynamic constraints imposed by VSL on the desired speed introduced in the previous improvement (see Eq. (16)) and the asymmetric friction effect exerted on the subject lane by lateral lane-changing behaviors (see Eq. (14)), the comprehensive lane-level dynamic speed evolution equation ultimately constructed in this study is formulated as Eq. (20).

$$
\begin{aligned}
v_{m,n}^{i+1} = v_{m,n}^{i} + \frac{\Delta t}{\tau}\left[ V\left(k_{m,n}^{i}\right) - v_{m,n}^{i} \right] + \frac{\Delta t}{\Delta x_m} \cdot v_{m,n}^{i} \cdot \left( \sqrt{v_{m-1,\tilde{n}}^{i} \cdot v_{m,n}^{i}} - v_{m,n}^{i} \right) - \frac{\eta \cdot \Delta t \cdot \left( k_{m+1,\tilde{n}}^{i} - k_{m,n}^{i} \right)}{\tau \cdot \Delta x_m \cdot \left( k_{m,n}^{i} + \kappa \right)} \\
- \frac{\delta_{in} \cdot \Delta t \cdot v_{m,n}^{i} \cdot \phi_{m,n}^{i,in}}{\Delta x_m \cdot \left( k_{m,n}^{i} + \kappa \right)} - \frac{\delta_{out} \cdot \Delta t \cdot v_{m,n}^{i} \cdot \phi_{m,n}^{i,out}}{\Delta x_m \cdot \left( k_{m,n}^{i} + \kappa \right)}
\end{aligned} \tag{20}
$$

**Improvement 5:** Dynamic flow estimation for lateral lane-changing behaviors involving free and forced lane-changing, which can be viewed in **Appendix A.**

*4.1.2 Calibration method of L-METANET*

Following the formulation of the theoretical L-METANET model, its parameters must be precisely calibrated using road network data to ensure the model accurately captures the spatiotemporal evolution characteristics of traffic flow within the expressway network. Drawing upon systematic frameworks for macroscopic traffic flow modeling, this section establishes a parameter calibration system across three dimensions: parameter set definition, evaluation variable selection, and an efficient solution algorithm.

(1) Parameter decoupling and calibration set formulation: the variables of L-METANET are decoupled into two distinct subsets to reduce optimization complexity: First, the global dynamic parameters encompass the driver reaction time $\tau$, the speed-density relationship coefficient $\eta$, the anticipation elasticity coefficient parameter $\kappa$, and the weight coefficients $\delta^{in}$ and $\delta^{out}$ that characterize the lateral asymmetric friction effects across lanes. Second, the local fundamental diagram parameters are utilized to characterize the geometric heterogeneity across road space (between different segments or different lanes). This subset includes the free-flow speed $v_{f,m,n}$, natural critical density $k_{cr,m,n}$, and shape parameter $a_{m,n}$ for a specific segment-lane unit $(m,n)$. All parameter vectors to be identified $\theta$ constitute the decision variable set for the calibration optimization problem, as shown in Eq. (21). To avoid calibration intractability caused by high-dimensional over-parameterization, this study assumes uniform fundamental diagram parameters across the entire network.

(2) Evaluation variable and objective function construction: Parameter calibration is essentially a parameter estimation problem for a nonlinear dynamic system. Referencing relevant studies (Wang et al., 2022), this study selects the average speed $v_{m,n}^{i}$ as the core evaluation metric. The objective function $J(\mathbf{\theta})$ is formulated to minimize the root-mean-square error (RMSE), as shown in Eq. (22).

(3) Hybrid global optimization algorithm based on GA-VNS: Owing to the complex extremum calculations and lane-level coupling embedded in the L-METANET model, its objective function space exhibits highly nonlinear and non-convex properties. Conventional gradient-based optimization algorithms are highly susceptible to trapping in local optima when solving such problems. Synthesizing the advantages of multiple meta-heuristic algorithms, this study adopts a hybrid global optimization strategy combining a genetic algorithm (GA) and variable neighborhood search (VNS) (Ma et al., 2025), which was specifically designed in our previous work for parameter calibration. This algorithm framework has been implemented multiple times in our previous studies, demonstrating robust population evolution capabilities and a strong capacity to escape local optima.

$$\theta = \left\{ \tau, \eta, \kappa, \delta^{in}, \delta^{out}, \mu_{CAV}, \mu_{HDV}, v_{f,m,n}, k_{cr,m,n}, a_{m,n} \right\}_{m=1,\ldots,M}^{n=1,\ldots,N} \tag{21}$$

$$J(\mathbf{\theta}) = \sqrt{\frac{1}{I \cdot \sum_{m=1}^{M} N_m} \sum_{i=1}^{K} \sum_{m=1}^{M} \sum_{n=1}^{N_m} \left( \hat{v}_{m,n}^{i} - v_{m,n}^{i} \right)^2} \tag{22}$$

where $v_{m,n}^{i}$ and $\hat{v}_{m,n}^{i}$ denote the calculated speed output by the L-METANET model and the empirically observed speed from detectors, respectively, for lane $n$ of segment $m$ at time step $i$; $I$ represents the total number of time steps within the calibration period, and $M$ represents the total number of calibrated segments.

4.2. Accuracy- and controllability-driven risk prediction and assessment model for weaving segments

Coordinated VSL and RM control requires a rigorous risk-based objective function balancing risk-characterization accuracy with operational controllability. Traditional macroscopic risk models neglect microscopic conflict mechanisms, while high-fidelity microscopic models are disconnected from macroscopic control parameters. To bridge this gap, this study identifies macroscopic precursors of merging and diverging risks to synthesize microscopic fidelity with macroscopic controllability. Furthermore, to resolve the trade-off between black-box machine learning models lacking closed-form expressions and traditional logit models suffering from multicollinearity, we develop a hybrid XGBoost-SHAP-RPBL framework. This framework achieves high-dimensional feature reduction, accounts for unobserved heterogeneity, and derives an analytical, closed-form risk probability output suitable for optimization.

*4.2.1 The XGBoost-SHAP-RPBL framework and implementation steps*

The implementation procedure of XGBoost-SHAP-RPBL consists of five core steps, as illustrated in the flowchart in Fig. 4. Note that Steps 1 through 3 were completed in our previous study (Ma et al., 2026a) and are thus only briefly outlined here, whereas Steps 4 and 5 represent the novel contributions of this study, with their detailed mathematical procedures elaborated in Section 4.2.2.

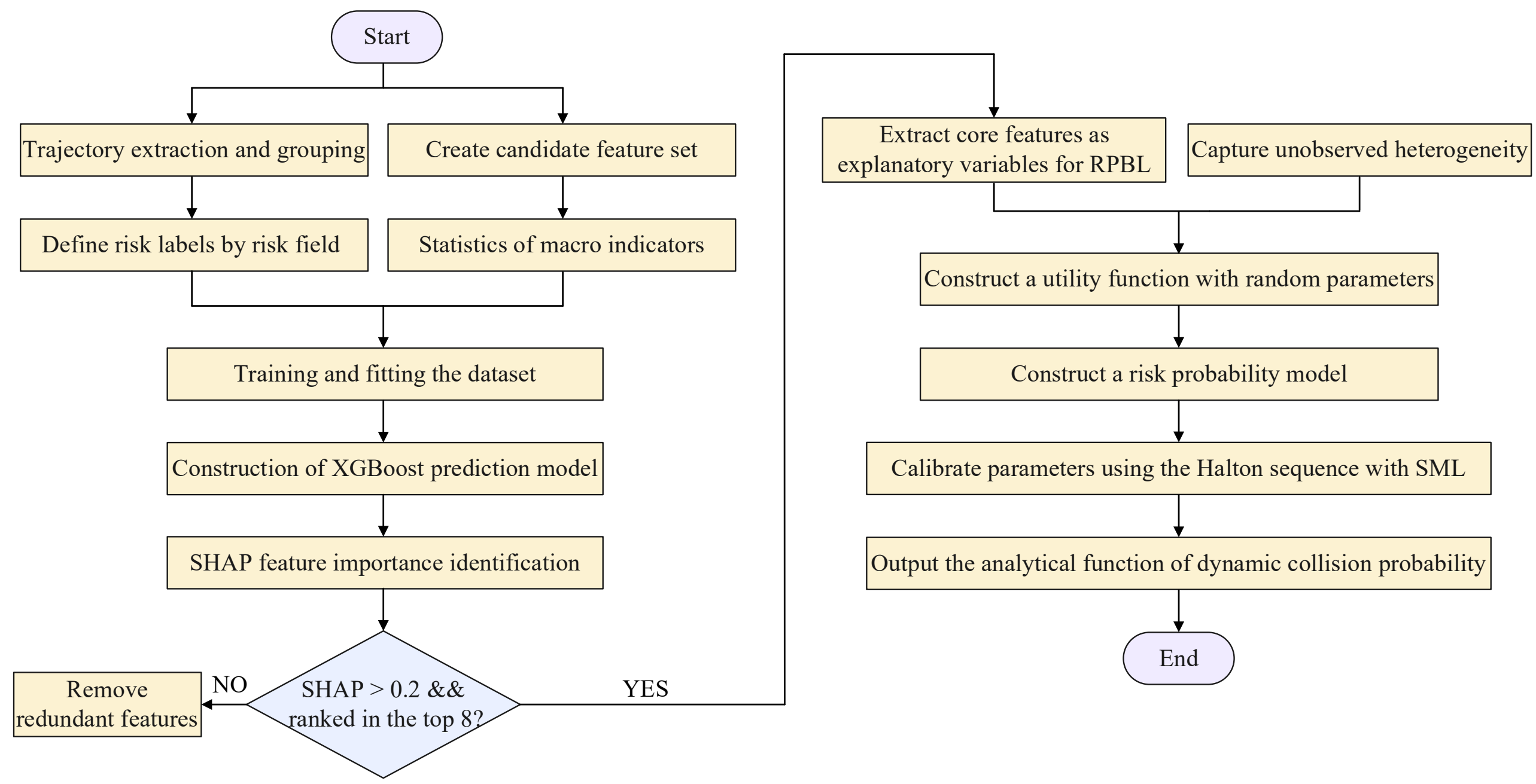

Fig. 4 Fitting process of XGBoost-SHAP-RPBL

**Step 1: Candidate feature set construction:** A three-dimensional feature space (lane-segment-factor) is designed to capture the complex spatiotemporal heterogeneity within weaving segments. Spatially, the weaving segment is partitioned longitudinally into upstream (Segment 1), midstream (Segment 2), and downstream (Segment 3) zones, and laterally discretized into the inner fast lanes and outer weaving lanes (Lane 1, Lane 2, etc.). Factorially, core macroscopic traffic parameters: cell density, average speed, flow rate, maximum speed, and speeding vehicle counts are extracted for each spatiotemporal cell.

**Step 2: dataset construction for training and fitting:** Empirical data are processed into a standardized dataset, aligning microscopic vehicle interactions with macroscopic traffic states via three substeps: Step 2.1: microscopic trajectory extraction: Vehicle trajectories are extracted from 3.6 hours of UAV aerial video data capturing seven short weaving segments on Changchun expressways under various congestion levels. Microscopic trajectory data for all vehicles within the weaving segments were extracted, and all merging and diverging vehicles, along with their surrounding vehicles, were identified. Each merging or diverging vehicle and its corresponding surrounding vehicles form a group, which constitutes a single sample; Step 2.2: macroscopic parameter aggregation: The weaving segment is discretized into continuous macroscopic time slices. Microscopic trajectories within each lane-segment grid cell are aggregated per time slice to extract cell-level traffic flow indicators, forming the independent variable candidate set; Step 2.3: STRF-based weaving risk classification: Based on the STRF model developed in our previous study (Ma et al., 2026b), the STRF is utilized to risk classification. A sample with a risk exceeding the STRF threshold is defined as high-risk, whereas the remaining samples are classified as low-risk.

**Step 3: construction of the XGBoost-based risk prediction model:** The dataset constructed in Step 2 is input into XGBoost for preliminary classification modeling. Utilizing a powerful ensemble of base learners and regularization mechanisms, XGBoost automatically suppresses the weights of redundant features while accurately capturing the complex nonlinear threshold responses and synergetic interactions between macroscopic traffic flow characteristics and merging/diverging risks.

**Step 4: SHAP-Driven Feature Selection for RPBL Modeling:** To resolve the black-box limitation of XGBoost and establish a parsimonious input set for subsequent statistical modeling, the SHAP framework calculates the global marginal contribution of each feature. Highly multicollinear or low-contribution features are eliminated based on two criteria: a SHAP importance threshold of 0.2 and a maximum allocation of eight variables. This dimensional reduction retains only the core physical variables with decisive impacts on crash risk (listed in **Appendix C**), which serve as the explanatory variables for the final Logit-based RPBL model.

**Step 5: RPBL model estimation and analytical risk equation derivation:** The core feature subset $\mathbf{X}^i$ selected via SHAP is utilized as the explanatory variables, and all samples from Step 2 are input into the RPBL model for joint estimation. By introducing random parameters governed by specific probability distributions to capture the unobserved heterogeneity in traffic flow, the model yields a smooth, differentiable analytical equation for the overall crash risk over the prediction horizon, thereby establishing a safety cost function for the control system.

*4.2.2 Risk prediction and assessment model based on random parameters binary Logit*

The final stage of the hybrid framework involves converting the reduced-dimensional core feature set $\mathbf{X}^i$ into a smooth, differentiable continuous probability equation. Since traffic flow within weaving segments is continuously exposed to highly uncertain open environments, random variations in microscopic driver aggressiveness, localized transient weather, and vehicle mechanical performance can cause the crash probability induced by identical macroscopic traffic states to fluctuate drastically. Conventional fixed-parameter Logit models fail to capture such unobserved heterogeneity within the data. To address this limitation, a random parameters binary logit (RPBL) model is developed as the final analytically expressed risk assessment and prediction model. The framework utilizes the feature state $\mathbf{X}^i$ at the current time step $i$ to evaluate the utility function $U^i$ for the entire weaving segment evolving into a high-risk conflict state at time step $i$, as formulated in Eq. (23).

$$U^i = \beta_0 + \sum_{k=1}^{K} \beta_k \cdot \mathbf{X}_k^i + \varepsilon^i \tag{23}$$

where $\beta_0$ is the constant intercept, and $\varepsilon^i$ is the random error term following an extreme value distribution. The $\mathbf{X}_k^i$ denotes the $k^{th}$ core traffic flow variable standardized via Z-score transformation to eliminate scale effects and ensure the convergence of the simulated maximum likelihood estimation algorithm, as expressed in Eq. (24).

$$\mathbf{X}_k^i = \frac{x_k^i - \mu_k}{\sigma_k} \tag{24}$$

where $x_k^i$ represents the raw value of the $k^{th}$ core traffic flow variable at time step $i$, while $\mu_k$ and $\sigma_k$ denote the mean and standard deviation of the $k^{th}$ core traffic flow variable, respectively.

The RPBL model relaxes the restrictive assumption of fixed regression coefficients by endowing the feature parameters $\beta_k$ with stochastic properties to capture spatiotemporal fluctuations. The random parameter $\beta_k$ is defined as Eq. (25).

$$\beta_k = \begin{cases} \bar{\beta}_k + \varphi_k, k \in K \\ \bar{\beta}_k, k \notin K \end{cases} \tag{25}$$

where $\bar{\beta}_k$ represents the fixed mean of the parameter across the population, characterizing the average positive impact of the feature on the weaving segment crash risk. Meanwhile, $\varphi_k$ is defined as a normally distributed random disturbance term, such that $\varphi_k \sim (0, \hat{\sigma}_k^2)$. Crucially, rather than treating all candidate variables as random parameters, only the core variables selected via SHAP constitute the explanatory variable set for the RPBL model. A fixed-parameters binary logit model is first estimated as the baseline. Subsequently, the significance of the random standard deviation for each core variable parameter is tested sequentially. The variable is identified as a random parameter, denoted as $k \in K$, only if its standard deviation is statistically significant at a given significance level ($p = 0.05$). Otherwise, it remains in the final model as a fixed parameter, denoted as $k \notin K$. This parameter-level formulation enables the model to effectively absorb unobserved stochastic noise within the physical traffic system.

Based on the utility formulation above, the dynamic probability of a merging/diverging crash occurring across the entire weaving segment at time step $i$, denoted as $CR^i$, is obtained by integrating the Logit cumulative distribution function, as shown in Eq. (26).

$$CR^i = \int \frac{\exp(U^i)}{1+\exp(U^i)} \cdot f(\beta) d\beta \approx \frac{1}{N_C} \sum_{d=1}^{N_C} \left[ \frac{1}{1+\exp(-U_d^i)} \right] \tag{26}$$

where $f(\beta)$ represents the joint probability density function of the random parameters. Because this integration lacks a closed-form analytical solution, this study employs SML estimation, utilizing Halton sequences to execute $N_C$ (set to 500) Monte Carlo to achieve high-precision parameter calibration.

The hybrid XGBoost-SHAP-RPBL model successfully distills massive, discrete traffic flow data characterized by heterogeneous fluctuations into a continuous dynamic analytical equation for crash probability bounded between $(0,1)$. This function directly serves as the core objective for the coordinated VSL-RM optimization, driving the VSL and RM actuators to proactively manage the macroscopic states of the expressway, thereby mitigating merging and diverging crash risks within the weaving segment.

## 5. Hierarchical coordinated control framework for lane-level variable speed limits and ramp metering

Despite the high accuracy of the developed L-METANET model, macroscopic formulations inevitably suffer from environmental disturbances and minor model mismatches. Because MPC relies heavily on model fidelity, these perturbations challenge its control stability. Conversely, while standalone MARL provides robust high-frequency adaptability, its lack of physical boundaries and domain knowledge renders it susceptible to triggering traffic breakdowns in unseen scenarios. To address these challenges, we develop a hierarchical ATM framework for joint LVSL and RM control, establishing a closed-loop architecture where the upper-level MPC provides a rigid baseline safeguard and the lower-level MARL injects residual compensation.

### 5.1. Overall architecture

The hierarchical MPC-MARL coordinated control framework is illustrated in Fig. 5. The hierarchical architecture operates across three distinct temporal scales: (1) Macroscopic state evolution step $\triangle T$: The fundamental time step at which the underlying L-METANET model updates its macroscopic state variables; (2) MPC control horizon $T_c$: The lower-frequency operational step at which the upper-level MPC baseline controller executes its rolling horizon optimization; (3) MARL execution cycle $T_s$: The higher-frequency sampling and actuation interval at which the lower-level MARL residual compensator performs adaptive fine-tuning. These three time scales satisfy the multi-frequency temporal relationship formulated in Eq. (27).

$$T_c = \lambda_1 \cdot T_s = \lambda_1 \cdot \lambda_2 \cdot \Delta T \tag{27}$$

where $\lambda_1$ and $\lambda_2$ are both positive integers, satisfying: $T_c > T_s > \Delta T$.

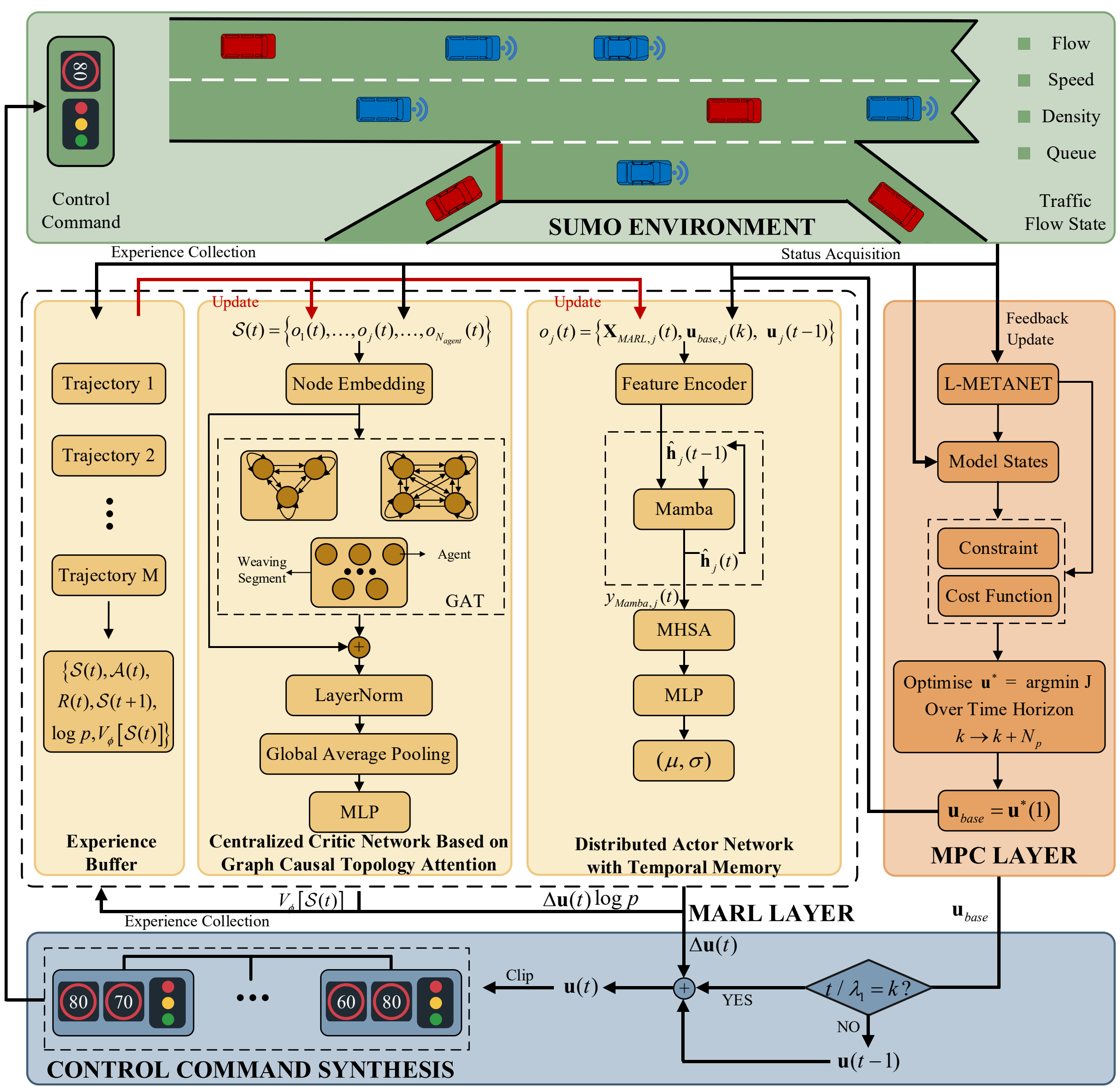

Fig. 5 The network architecture of MPC-STMAPPO

In the physical road network, control devices are deployed in a spatially discrete manner. Let $N_{VSL}$ be the total number of lane-level VSL gantries, defining the VSL actuator set as $c \in \{1,2,\ldots,N_{VSL}\}$, where each actuator $c$ maps to a specific segment-lane combination. Let $N_{RM}$ be the total number of ramp meters, defining the RM actuator set as $g \in \{1,2,\ldots,N_{RM}\}$, where each actuator $g$ maps to a specific on-ramp. Within the MPC framework, the state space and control inputs must be defined at the MPC control step $k$. The macroscopic state vector at MPC step $k$ is defined as $\mathbf{X}_{MPC}(k)$, which is physically equivalent to the output of the L-METANET model at the physical step $k \cdot \lambda_1 \cdot \lambda_2$. This mapping relation is defined in Eq. (28). Similarly, for the MARL framework, the state space $\mathbf{X}_{MARL}(k)$ at the high-frequency control step $t$ is physically equivalent to the output of the L-METANET model at the physical step $t \cdot \lambda_2$, as defined in Eq. (29).

$$\mathbf{X}_{MPC}(k) = \mathbf{X}_{L-METANET}(i)\big|_{i=k\cdot\lambda_1\cdot\lambda_2} \tag{28}$$

$$\mathbf{X}_{MARL}(t) = \mathbf{X}_{L-METANET}(i)\big|_{i=t\cdot\lambda_2} \tag{29}$$

Regarding the overall operational logic of the hierarchical coordinated architecture, at the beginning of each MPC step, the upper-level MPC solver predicts the system evolution within the prediction horizon via rolling optimization based on the current macroscopic traffic flow state combined with the L-METANET model, and issues a low-frequency fixed baseline control command $\mathbf{u}_{base} = [\mathbf{v}_{base}, \mathbf{r}_{base}]^T$ for the current step. Within each MARL step, the lower-level MARL agents deployed at each lane and ramp output high-frequency residual actions $\Delta\mathbf{u} = [\Delta\mathbf{v}, \Delta\mathbf{r}]^T$ online, according to the high-frequency sampled macroscopic state and the baseline control command determined by the MPC. The final issued control command is related to $\mathbf{u}_{base}$ and $\Delta\mathbf{u}$, but instead of a static superposition of the two, it is reconstructed as a piecewise chronological rolling cumulative generation law. Its action synthesis and truncation protection mechanism are shown in Eqs. (30) to (32). This control superposition mechanism mathematically ensures that when the MARL step exactly satisfies $t/\lambda_1 = k$, it indicates that the control loop has transitioned to the MPC cycle boundary, at which time $\mathbf{u}_{ref}(t)$ is forcibly anchored to the latest global optimal solution $\mathbf{u}_{\text{base}}(k)$ issued by the upper-level MPC; whereas when $t/\lambda_1 \neq k$, it indicates that the system is in the intra-cycle rolling fine-tuning stage, and $\mathbf{u}_{ref}(t)$ seamlessly transitions to the actual physical execution command of the previous MARL step $\mathbf{u}(t-1)$. This superposition mechanism ensures that, on one hand, the MPC can correct and guide the MARL actions in every cycle, keeping the final issued values permanently within an acceptable range; on the other hand, the system eliminates discrete abrupt mutations in speed limits and metering rates within the dynamic execution structure, securing a smooth physical transition of the overall control boundary.

$$\mathbf{u}(t) = \text{Clip}\left(\mathbf{u}_{ref}(t) + \Delta\mathbf{u}(t)\right) \tag{30}$$

$$Clip(\mathbf{u}) = \begin{cases} \mathbf{u}_{\max}, & u > u_{\max} \\ \mathbf{u}_{\min}, & u < u_{\min} \\ \mathbf{u}, & \text{otherwise} \end{cases} \tag{31}$$

$$\mathbf{u}_{\text{ref}}(t) = \begin{cases} \mathbf{u}_{\text{base}}(k), & t/\lambda_1 = k \\ \mathbf{u}(t-1), & t/\lambda_1 \neq k \end{cases} \tag{32}$$

where $\Delta\mathbf{u}(t)$ is the physical increment vector transformed through spatial mapping from the continuous residual actions output by the MARL layer; $\mathbf{u}_{\min}, \mathbf{u}_{\max}$ are the legal boundary extreme vectors of the discrete physical actuators; and the upper-level low-frequency step index $k$ and the lower-level high-frequency step index $t$ satisfy the chronological rounding mapping relationship: $k = \left[t/\lambda_1\right]$.

5.2. Upper-level baseline controller based on model predictive control

The core task of the upper-level MPC is to deduce the spatiotemporal evolution of traffic flow using L-METANET macroscopic dynamic equations. Combined with the RPBL dynamic crash probability equations constructed in Section 4.2, it establishes a baseline control boundary balancing safety and efficiency for the weaving segments from a global perspective over a long prediction horizon.

*5.2.1 State-space and prediction model formulation*

Defining Eq. (33) as the state space of the MPC, the upper-level baseline control vector $\mathbf{u}_{base}(k)$ represents the command set issued to all independent physical actuators, as shown in Eq. (34). In the MPC rolling optimization, since the issued baseline command $\mathbf{u}_{base}(k)$ remains constant throughout the entire low-frequency control cycle $T_c$, the underlying L-METANET dynamic model must absorb this constant command to continuously execute state updates at step size $\triangle T$. By inputting the baseline traffic flow operational data of the expressway system and the environmental boundary demands $\mathbf{D}(i)$, the improved L-METANET model predicts future traffic flow states across sampling periods. The state prediction transition equation is formulated as Eq. (35).

$$\mathbf{X}_{MPC}(k) = \left[ k_{m,n}^{k\cdot\lambda_1\cdot\lambda_2}, v_{m,n}^{k\cdot\lambda_1\cdot\lambda_2}, w_g^{k\cdot\lambda_1\cdot\lambda_2} \right]^T \tag{33}$$

$$\mathbf{u}_{base}(k) = \left[ v_1^{VSL}(k),\ldots,v_{N_{VSL}}^{VSL}(k), \quad r_1^{RM}(k),\ldots,r_{N_{RM}}^{RM}(k) \right]^T \tag{34}$$

$$\mathbf{X}_{L-METANET}(i+1) = F_{L-METANET}\left[\mathbf{X}_{L-METANET}(i), \mathbf{u}_{base}(k), \mathbf{D}(i)\right],\ k\cdot\lambda_1\cdot\lambda_2 \le i \le (k+1)\cdot\lambda_1\cdot\lambda_2 - 1 \tag{35}$$

where $\mathbf{X}_{L-METANET}(i+1)$ is the predicted basic traffic flow state vector at sampling step $i+1$; $\mathbf{X}_{L-METANET}(i)$ is the input basic traffic flow state vector at sampling step $i$; $\mathbf{D}(i)$ is the traffic boundary demand vector (arrival flows of the upstream mainline and ramps) at sampling step $i$; and $F_{L-METANET}(\cdot)$ represents the mapping function of L-METANET.

In the MPC rolling optimization, the controller predicts forward from the current control step $k$ to $k+p+1$ (where the prediction step is $p \in \left[0, N_p - 1\right]$ and $N_p$ is the MPC prediction horizon). Based on the zero-order hold principle, the upper-level baseline command remains constant throughout the entire low-frequency control cycle, while the underlying L-METANET model executes $\lambda_1 \cdot \lambda_2$ physical iterative steps with a step size of $\triangle T$ within this cycle.

*5.2.2 Multi-objective comprehensive cost function formulation*

The core pain point of weaving segment management lies in the inherent trade-off between traffic efficiency and operational safety. To address this, the upper-level MPC constructs a comprehensive cost function $J_{MPC}$ across the prediction horizon, as shown in Eq. (36).

$$J_{MPC} = \sum_{i=k\cdot\lambda_1\cdot\lambda_2}^{(k+N_p)\cdot\lambda_1\cdot\lambda_2-1} \gamma_i \left[\omega_1 \cdot CR(i) + \omega_2 \cdot J_{eff}(i) + \omega_3 \cdot J_{queue}(i)\right] + J_{smooth}(k) \tag{36}$$

$$\gamma_i = \frac{\exp\left(-\dfrac{i-k\cdot\lambda_1\cdot\lambda_2}{\eta_d}\right)}{\sum\limits_{j=0}^{P=N_p\cdot\lambda_1\cdot\lambda_2-1} \exp\left(-\dfrac{j}{\eta_d}\right)} \tag{37}$$

where $\omega_1$, $\omega_2$, and $\omega_3$ are priority weights, and $J_{\text{penalty}}(k)$ is the smoothing penalty term for temporal and spatial variations. Within the total prediction horizon, the normalized temporal discount coefficient $\gamma_i$ at any prediction step $i$ is expressed as in Eq. (37), where $\eta_d$ represents the temporal discount adjustment parameter reflecting the decay of predictive efficacy over the macroscopic prediction horizon. In Eq. (36), the first term represents the safety component, the second is the mainline efficiency component, the third denotes the on-ramp vehicle efficiency component, and the fourth constitutes the additional smoothing penalty for control execution. These four components are detailed below:

(1) XGBoost-SHAP-RPBL based merge and diverge risk penalty $CR(i)$: The lane-level density and speed feature subsets deduced by the prediction model at time $i$ are directly substituted into the RPBL dynamic crash probability integral equation derived in Eq. (26). Through this cost term, the MPC optimizer gains risk-foreseeing capabilities, forcing it to proactively search for VSL and RM combinations that minimize merge and diverge conflicts.

(2) Mainline traffic efficiency delay cost $J_{eff}(i)$: The efficiency evaluation is restricted to the target macroscopic segments composed of the weaving area and its upstream sections (the set of all evaluation zones is denoted as $\mathcal{M}_{\text{target}}$). The average speed within the target segments is utilized as the efficiency metric, with the specific structure formulated in Eq. (38).

(3) On-ramp queue delay cost $J_{queue}(i)$: The total queuing time on the on-ramps within the weaving segment is adopted as the ramp queue penalty, as shown in Eq. (39).

(4) Actuator discrete smoothing penalty $J_{smooth}(k)$: To prevent severe command jumps when spatially discrete actuators issue control signals, strict smoothing constraints must be imposed. Specifically, for VSL control, the system must restrict not only the temporal command discrepancy of a single speed limit gantry between adjacent time steps but also the spatial command discrepancy between adjacent gantries on the same lane at the same timestamp, thereby preventing rear-end collision risks induced by sudden spatial speed limit drops. To this end, let $\varepsilon_{VSL}$ be the set of spatially adjacent VSL actuator pairs on the same lane. If actuator $c'$ is physically located immediately downstream of actuator $c$ on the same lane, the pair $(c,c')\in\varepsilon_{VSL}$ exists. A quadratic penalty for control increments is introduced across all independent discrete actuators within the control horizon, as shown in Eq. (40).

$$J_{eff}(i) = -\frac{\sum_{m\in\mathcal{M}_{\text{target}}}\sum_{n=1}^{N_m} q_{m,n}(i)}{\sum_{m\in\mathcal{M}_{\text{target}}}\sum_{n=1}^{N_m} k_{m,n}(i)\cdot v_{\max} + \epsilon} \tag{38}$$

$$J_{queue}(i) = \sum_{g=1}^{N_{RM}} w_{og}(i)\cdot \Delta T / 3600 \tag{39}$$

$$\begin{aligned} J_{smooth}(k) = \sum_{p=0}^{N_c-1}\Bigg\{ & \sum_{c=1}^{N_{VSL}}\left[\frac{v_c^{VSL}(k+p)-v_c^{VSL}(k+p-1)}{v_{max}}\right]^2 + \\ & \sum_{(c,c')\in\mathcal{E}_{VSL}}\left[\frac{v_c^{VSL}(k+p)-v_{c'}^{VSL}(k+p)}{v_{max}}\right]^2 + \sum_{g=1}^{N_{RM}}\left[r_g^{RM}(k+p)-r_g^{RM}(k+p-1)\right]^2\Bigg\} \end{aligned} \tag{40}$$

*5.2.3 Traffic physical boundaries and control execution constraints*

To ensure the physical feasibility and absolute safety of the control sequence solved throughout the prediction horizon of the upper-level MPC rolling optimization, a strict system of equality and inequality constraints must be imposed, encompassing the following aspects:

(1) Macroscopic traffic flow dynamics equality constraints. The evolution of all state vectors within the prediction horizon must strictly obey the conservation laws of L-METANET.

(2) VSL spatial execution boundary and control rate constraints. For any speed limit gantry $c\in\{1,\ldots,N_{VSL}\}$, its issued command must fall between the maximum speed limit $v_{\max}$ allowed by the physical road network and the minimum safe driving speed $v_{\min}$, as shown in constraint (41). Meanwhile, to prevent rear-end collisions caused by sharp drops in speed limits, the speed limit variation within a single low-frequency step is restricted from exceeding the physical threshold $\Delta v_{step}$ (constraint (42)), and the command difference between adjacent upstream and downstream speed limit gantries on the same lane at the same prediction timestamp is restricted from exceeding the spatial safety threshold $\Delta v_{spacc}$ (constraint (43)). Generally, the magnitudes of $\Delta v_{step}$ and $\Delta v_{spacc}$ remain consistent.

(3) Queue constraints. At any underlying physical time step within the prediction horizon, the number of queuing vehicles not exceed the maximum physical storage capacity of the ramp, eliminating deadlocks and spillovers at the hard constraint level, as shown in constraint (44).

(4) RM rate fluctuation constraints. Sudden and large variations in ramp inflow will severely degrade the traffic efficiency of the weaving segment and elevate accident risks. Therefore, high-frequency fluctuations must be suppressed, as shown in constraint (45).

$$v_{min} \le v_c^{VSL}(k+p) \le v_{max}, \quad \forall p \in [0, N_c - 1] \tag{41}$$

$$\left| v_c^{VSL}(k+p) - v_c^{VSL}(k+p-1) \right| \le \Delta v_{step}, \quad \forall c \in \{1, \ldots, N_{VSL}\}, \forall p \in [0, N_c - 1] \tag{42}$$

$$\left| v_c^{VSL}(k+p) - v_{c'}^{VSL}(k+p) \right| \le \Delta v_{spacc}, \quad \forall (c, c') \in \mathcal{E}_{VSL}, \forall p \in [0, N_c - 1] \tag{43}$$

$$w_{o_g}^i \le w_{max}, \quad \forall i \in \left[ k \cdot \lambda_1 \cdot \lambda_2, (k+N_p) \cdot \lambda_1 \cdot \lambda_2 \right] \tag{44}$$

$$\left| r_g^{RM}(k+p) - r_g^{RM}(k+p-1) \right| \le \Delta r_{step}, \quad \forall p \in [0, N_c - 1] \tag{45}$$

*5.2.4 Model optimization and solution*

At each low-frequency time step $k$, the upper-level baseline controller aims to solve for an optimal control sequence $U^*$ that minimizes the multi-objective comprehensive cost function $J_{MPC}$ (Eq. (36)) while strictly adhering to constraints (41) to (45). Because this optimization task constitutes a characteristically high-dimensional, non-linear, and non-convex optimization problem, conventional gradient-based deterministic solvers, such as sequential quadratic programming, are highly susceptible to becoming trapped in local optima and frequently fail to converge. To circumvent these limitations, this study employs a hybrid GA-VNS algorithm for optimization. Upon completing the hybrid optimization for the current step, the system extracts only the first element of the optimal control sequence $\mathbf{u}_{base}(k) = \mathbf{U}^*(1)$ as the deterministic baseline command to be dispatched to the physical road network. At the next low-frequency control step $k+1$, the closed-loop rolling horizon optimization process is re-initialized by introducing the fresh lane-level macroscopic state feedback collected via roadside sensors.

5.3. Lower-level residual compensator based on multi-agent reinforcement learning

This study models the lower level as a fully cooperative multi-agent system. Within each low-frequency baseline cycle, agents distributed across each lane and ramp of the weaving segment continuously sample local macroscopic states at a higher frequency (step size $T_s$) and output continuous residual compensation actions online, thereby agilely smoothing model prediction mismatches and achieving adaptive closed-loop control of the macroscopic traffic flow.

In the complex spatiotemporal coupled road network of the weaving segment, distributed speed limit gantries and ramp signals cannot acquire the global state of the entire network and must make coordinated decisions based on local information. Therefore, this study rigorously models the lower-level residual compensation-based multi-agent coordinated control problem as a decentralized partially observable Markov decision process (Dec-POMDP). This system can be mathematically defined by a 7-tuple $(\mathcal{J}, \mathcal{S}, \mathcal{A}, \mathcal{P}, R, \Omega, \gamma_d)$, with the specific mapping of each element under this framework as follows:

(1) Agent set $\mathcal{J}$: Consists of $N_{VSL}$ variable speed limiters and $N_{RM}$ ramp meters in the road network, with a total number of $N_{agent}=N_{VSL}+N_{RM}$ and the agent index as $j \in \mathcal{J}$;
(2) Global state space $\mathcal{S}$: Characterizes the true global macroscopic traffic flow state of the road network at any high-frequency step $t$.
(3) Set of partially observable state spaces for each agent $\Omega$: Agents cannot acquire the global $\mathcal{S}$ and can only acquire a local observation vector $o_j(t) \in \Omega$ containing high-dimensional features.
(4) Joint action space $\mathcal{A}$: The set of continuous residual actions $\Delta \mathbf{u}(t)$ output by all agents at high-frequency step $t$.
(5) State transition probability $\mathcal{P}$: Implicitly driven by the underlying simulation environment, which here refers to SUMO.
(6) Reward function $R$: A fully cooperative feedback signal that guides the agents to jointly compromise toward the global macroscopic optimum.
(7) Discount factor $\gamma_d$: $\gamma_d \in [0,1)$, which measures the current importance of future spatiotemporal traffic flow rewards.

*5.3.1 State-space formulation incorporating upper-level MPC outputs*

The inputs of MARL require extensive feature engineering to enhance the Markov property of the environment. This study deeply extends the local observation space $o_j(t)$ of agent $j$, constructing it into a high-dimensional feature vector containing three core information matrices, as shown in Eq. (46). The global state space at step $t$ is then the concatenation of the local observation spaces of all agents, as shown in Eq. (47).

$$o_j(t) = \left[ \mathbf{X}_{MARL,j}(t), \mathbf{u}_{base,j}(k),\ \mathbf{u}_j(t-1) \right] \tag{46}$$

$$\mathcal{S}(t) = \left[ o_1(t), \ldots, o_j(t), \ldots, o_{N_{agent}}(t) \right] \tag{47}$$

The specific definitions of the aforementioned observation features are as follows:

**Current macroscopic physical state $\mathbf{X}_{MARL,j}(t)$:** Since weaving segment congestion exhibits a significant downstream-to-upstream propagation characteristic, to enable the deep neural network to precisely capture the causal topological relationships of traffic shock waves and maximize feature noise elimination, this study classifies the agent set into three categories based on the physical attributes of the actuators: inner-lane VSL agents ($\mathcal{J}_{VSL}^{in}$), outer-lane VSL agents ($\mathcal{J}_{VSL}^{out}$), and on-ramp RM agents ($\mathcal{J}_{RM}$), constructing their states respectively, as shown in Eq. (48). For inner-lane VSL agents, the density, speed, and flow at three critical nodes—the VSL zone, acceleration zone, and weaving segment—are extracted. The outer lane is not only a severe conflict zone for longitudinal weaving but also directly faces the pressure of the merging flow from the on-ramp. Therefore, in addition to having a longitudinal three-segment perception capability completely symmetrical to that of the inner lane, the outer VSL agents must incorporate the observation of the queue storage on their connected on-ramp $oj$. The core logic of RM is to control the merging flow rate based on the residual capacity of the mainline outer lane to prevent mainline breakdown or ramp gridlock. Consequently, the RM agent focuses on the queue storage and high-frequency dynamic arrival demand of its own ramp, as well as the loading state of the connected mainline target merging lane. The representations of these macroscopic physical states can be observed as Eq. (48).

**Upper-level baseline control output $\mathbf{u}_{base,j}(k)$:** The baseline command $\mathbf{u}_{base,j}(k)$ received by agent $j$ within its current $k^{th}$ low-frequency control step, which is issued by the upper-level MPC rolling optimization. Consequently, the neural network can clearly identify the current macroscopic safety anchor, thereby avoiding ineffective exploration into invalid action spaces that violate physical bottom lines when generating high-frequency residuals action $\Delta\mathbf{u}_j(t)$ under the guidance of the MPC layer.

**Previous total control input $\mathbf{u}_j(t-1)$:** This feature extracts the final synthesized command actually applied to the road network by the agent at time $t-1$. In a Markov decision process, the current state of physical actuators is a crucial component of the environmental state. Feeding the action of the previous step as the observation input of the current step forces the Actor network to possess temporal memory. Coupled with the residual smoothing penalty in the reward function, this feature effectively guides the neural network to self-constrain its action span, preventing the output of fluctuating commands and ensuring a safe, smooth transition of vehicle operations from the algorithm input end.

Notably, all observation states are normalized, and padding zeros are applied to missing elements to maintain consistent state dimensions.

$$\mathbf{X}_{MARL,j}(t) = \begin{cases} \left[k_{AZ,nj}(t), v_{AZ,nj}(t), q_{AZ,nj}(t), k_{mj,nj}(t), v_{mj,nj}(t), q_{mj,nj}(t), k_{WS,nj}(t), v_{WS,nj}(t), q_{WS,nj}(t)\right]^T, \ j \in \mathcal{J}_{VSL}^{in} \\ \left[k_{AZ,nj}(t), v_{AZ,nj}(t), q_{AZ,nj}(t), k_{mj,nj}(t), v_{mj,nj}(t), q_{mj,nj}(t), k_{WS,nj}(t), v_{WS,nj}(t), q_{WS,nj}(t), w_{oj}(t)\right]^T, \ j \in \mathcal{J}_{VSL}^{out} \\ \left[w_{oj}(t), d_{oj}^{i}(t), k_{\tilde{m}j,\tilde{n}j}(t), v_{\tilde{m}j,\tilde{n}j}(t), q_{\tilde{m}j,\tilde{n}j}(t)\right]^T, \ j \in \mathcal{J}_{RM} \end{cases} \tag{48}$$

where the subscript $AZ, nj$ denotes the traffic flow parameters of the acceleration zone in the same lane $n$ as agent $j$; $mj, nj$ denotes the traffic flow parameters of the same lane $n$ within the same segment $m$ as agent $j$; $WS, nj$ denotes the traffic flow parameters of the weaving segment in the same lane $n$ as agent $j$; $\tilde{m}j, \tilde{n}j$ denotes the traffic flow parameters of the segment $m$ and lane $n$ connected to agent $j$; and $oj$ denotes the traffic flow parameters of the ramp where the RM agent is located.

*5.3.2 Joint action space design*

During the decentralized execution phase, all VSL and RM agents within the weaving segment jointly output decisions at each high-frequency time step $t$, forming the joint action space $\mathcal{A}(t) = \left[\Delta\mathbf{u}_1(t), \Delta\mathbf{u}_2(t), \ldots, \Delta\mathbf{u}_{N_{agent}}(t)\right]$. To guarantee the absolute safety bottom line of road network control, the MARL agents under the proposed framework do not directly output the final absolute control speed limits or absolute metering rates. Instead, agent $j$ outputs a continuous action vector bounded between $[-1, 1]$, which is mapped into a physical correction value after denormalization. The core task of the agents is to perform high-frequency, small-scale fine-tuning on the basis of $\mathbf{u}_{ref}(t)$, thereby agilely absorbing and smoothing out traffic shock wave disturbances triggered by model mismatches or sudden demand surges. Ultimately, the residual actions generated by the agents are superposed with the corresponding baseline commands to synthesize the theoretically optimal control law. Before being dispatched to the physical actuators, these synthesized commands uniformly pass through the underlying physical boundary truncation mechanism of the system, securing the operational safety bottom line of the weaving segment from the physical layer while granting MARL full exploration freedom.

*5.3.3 Fully cooperative shared global reward function formulation*

To prevent the multi-agent system from falling into destructive, self-interested local games, this system adopts a fully shared reward mechanism. At any high-frequency control step $t$, all VSL and RM agents are assigned a completely identical global comprehensive instantaneous reward. This reward references the upper-level MPC, as shown in Eq. (49), and is jointly composed of four dimensions: weaving segment safety, mainline travel efficiency, weaving segment on-ramp queuing, and control smoothness (Eqs. (50) to (51)). Among these, the first three components are completely identical to the corresponding terms in the cost function of the MPC layer.

$$
\begin{aligned}
&R_{global}(t) = R_j(t) = \\
&\frac{1}{\lambda_2}\sum_{i=t\cdot\lambda_2}^{(t+1)\cdot\lambda_2-1}\left[-\omega_1\cdot CR(i)-\omega_2\cdot J_{eff}(i)-\omega_3\cdot J_{queue}(i)\right]+\omega_{smooth}\cdot R_{smooth}(t),\quad \forall j\in\mathcal{J}
\end{aligned}
\tag{49}
$$

$$
\begin{aligned}
R_{smooth}(t) = -\Bigg[&\sum_{c\in\mathcal{J}_{VSL}}\left(\frac{\Delta v_c^{VSL}(t)-\Delta v_c^{VSL}(t-1)}{v_{max}}\right)^2+\sum_{(c,c')\in\mathcal{E}_{VSL}}\left[\frac{v_c^{VSL}(t)-v_{c'}^{VSL}(t)}{v_{max}}\right]^2 \\
&+\sum_{g\in\mathcal{J}_{RM}}\left(\frac{\Delta r_g^{RM}(t)-\Delta r_g^{RM}(t-1)}{r_{max}}\right)^2\Bigg]
\end{aligned}
\tag{50}
$$

$$
v_c^{VSL}(t) = v_{ref,c}^{VSL}(t)+\Delta v_c^{VSL}(t) \tag{51}
$$

### *5.3.4 ST-MAPPO: MAPPO algorithm architecture with spatiotemporal attention mechanisms*

Under the joint multi-lane multi-ramp control of weaving segments, the evolution of traffic congestion exhibits significant time lags and spatial causal traceability. Conventional DRL algorithms often struggle to effectively capture strongly coupled spatiotemporal features when processing such high-dimensional, non-stationary traffic states. To address this, this study proposes a multi-agent proximal policy optimization algorithm with spatiotemporal attention mechanisms (ST-MAPPO) integrating the Mamba selective state-space model and graph causal topological attention (GCTA), applying it within a centralized training, decentralized execution (CTDE) architecture.

**Distributed actor network with temporal memory.** Distributed agents in the weaving segment can only obtain local observations $o_j(t)$, positioning the system within a POMDP. To eliminate environmental non-Markovian properties, this framework reconstructs the Actor network structure by introducing Mamba blocks based on selective state spaces. At each high-frequency decision step $t$, the agent inputs not only the current normalized $o_j(t)$ but also incorporates a temporally maintained hidden state memory vector $\hat{\mathbf{h}}_j(t-1)$. The forward propagation and action output distribution parameter computation of the Actor network are formulated as Eq. (52).

$$
\left[\mu_j(t),\sigma_j(t),\hat{\mathbf{h}}_j(t)\right] = MambaCell\left[o_j(t),\hat{\mathbf{h}}_j(t-1)\right] \tag{52}
$$

where $MambaCell(\cdot)$ represents the computational operator for parameter selective decay and state updates through time-varying control matrices; $\mu_j(t)$ and $\sigma_j(t)$ are the mean and standard deviation of the Gaussian distribution for the continuous residual actions output by agent $j$, respectively; and $\hat{\mathbf{h}}_j(t)$ is the updated temporal hidden state vector. The agent ultimately samples the high-frequency residual action $\Delta\mathbf{u}_j(t)\sim\mathcal{N}(\mu_j(t),\sigma_j(t))$ from this Gaussian distribution.

**Centralized critic network based on graph causal topological attention (GCTA).** On the training side, the centralized Critic network evaluates the value of the global state matrix $\mathcal{S}(t)$. To prevent fully connected structures from introducing spatially uncorrelated traffic flow noise, the Critic network incorporates a spatial causal topological edge index architecture. The underlying static graph topology $\mathcal{G}=(\mathcal{V},\mathcal{E})$ is constructed using the physical correlations of the weaving bottleneck segments. Where the node set $\mathcal{V}$ represents the hardware actuators in the entire road network; the edge set $\mathcal{E}$ binds the lane-level VSL gantry nodes and their adjacent on-ramp RM signal nodes within the same weaving bottleneck region into a locally fully connected spatial clique, with self-loop edges configured for all nodes. For any control node pair $(j,l)\in\mathcal{E}$ with an edge connection relationship, their dynamic self-attention weight distribution coefficient $\alpha_{j,l}^{h}(t)$ under an independent attention head $h$ is solved as shown in Eq. (53). In the multi-head attention feature aggregation, a residual connection mechanism is introduced to ensure the training stability of the deep network, and the expression of the updated node feature $\mathbf{z}_j^{'}(t)$ is shown in Eq. (54).

$$\alpha_{j,l}^{h}(t)=\frac{\exp\left(\mathbf{a}_h^T\left(\text{LeakyReLU}\left(\mathbf{W}_{\text{src}}^{h}\cdot\mathbf{z}_j(t)+\mathbf{W}_{\text{dst}}^{h}\cdot\mathbf{z}_l(t)\right)\right)\right)}{\sum_{m\in\mathcal{N}_j\cup\{j\}}\exp\left(\mathbf{a}_h^T\left(\text{LeakyReLU}\left(\mathbf{W}_{\text{src}}^{h}\cdot\mathbf{z}_j(t)+\mathbf{W}_{\text{dst}}^{h}\cdot\mathbf{z}_m(t)\right)\right)\right)} \tag{53}$$

$$\mathbf{z}_j^{'}(t)=\text{LayerNorm}\left(\prod_{h=1}^{H}\sum_{l\in\mathcal{N}_j\cup\{j\}}\alpha_{j,l}^{h}(t)\cdot\mathbf{W}_{\text{src}}^{h}\mathbf{z}_l(t)\right)+\mathbf{W}_{\text{res}}\cdot\mathbf{z}_j(t) \tag{54}$$

where $\mathbf{z}_j(t)$ is the node feature vector of the current node transformed by the pre-encoding layer; $\mathbf{W}_{\text{src}}^{h}$ and $\mathbf{W}_{\text{dst}}^{h}$ are the attention transformation projection matrices; $\mathbf{a}_h^T$ is the learnable attention vector; and $\mathcal{N}_j$ is the neighborhood set of node $j$. $\text{LayerNorm}(\cdot)$ is the layer normalization operator; and $\mathbf{W}_{\text{res}}$ is the linear projection matrix for the residual connection.

The Critic network ultimately outputs the global state value estimate $V_\phi[\mathcal{S}(t)]$ through a global average pooling layer and a multi-layer perceptron (MLP), as shown in Eq. (55).

$$V_\phi[\mathcal{S}(t)]=\text{MLP}\left(\frac{1}{N_{\text{agent}}}\sum_{j\in\mathcal{V}}\mathbf{z}_j^{'}(t)\right) \tag{55}$$

In addition to the network structure, the following operations are implemented:

(1) Generalized advantage estimation (GAE) based on shared rewards. The attention-enhanced Critic network utilizes the fully shared global reward $R_{global}(t)$ to calculate the temporal difference error $\delta(t)$. Combined with the future reward discount factor $\gamma_d$ and the GAE smoothing parameter $\lambda_{GAE}$, the global advantage function $\hat{A}(t)$ is computed, as shown in Eqs. (56) to (57). Since all agents share the same global advantage signal $\hat{A}(t)$, the individual Actor networks are forced to align during updates, iteratively updating weights uniformly toward maximizing the macro-coordinated benefits of the entire weaving segment road network.

(2) Clipped update mechanism and loss function design. To ensure smooth convergence of policy updates and avoid severe oscillations in the road network caused by excessive residual exploration, the ST-MAPPO algorithm introduces an importance sampling ratio $\rho_j(t)$ and a proximal policy clipping mechanism into the Actor network. The loss function of the Actor is given in Eqs. (58) to (59). $S[\cdot]$ is the policy entropy regularization term (with coefficient $\beta$), encouraging agents to fully explore the boundaries of high-frequency residual actions during early training. The Critic loss function employs a mean squared error mechanism to approximate the true return value, as shown in Eq. (58).

(3) Implementation of other training tricks. First, sequential time-slice experience buffer: To preserve the continuous temporal dependencies required by the Mamba selective state-space model, trajectory replay avoids random shuffling. Instead, mini-batches are strictly partitioned by chronological time slices to facilitate sequential deductive learning. Second, dynamic feature and advantage normalization: Input state features are adaptively scaled using a moving average operator to accelerate optimization. Simultaneously, advantage functions undergo batch normalization within each mini-batch to maintain a standard distribution, stabilizing gradients against abrupt traffic flow fluctuations. Third, orthogonal initialization and linear decay: Network hidden weights are orthogonally initialized. Furthermore, the learning rates of the Actor and Critic networks, alongside the policy entropy coefficient, decay linearly across training epochs to balance aggressive early-stage exploration with stable late-stage convergence.

$$\delta(t) = R_{global}(t) + \gamma_d \cdot V_\phi\left[\mathcal{S}(t+1)\right] - V_\phi\left[\mathcal{S}(t)\right] \tag{56}$$

$$\hat{A}(t) = \sum_{l=0}^{T_{\max}-t-1} (\gamma_d \cdot \lambda_{GAE})^l \cdot \delta(t+l) \tag{57}$$

$$\rho_j(t) = \frac{\pi_{\theta j}\left[\Delta \mathbf{u}_j(t) \mid o_j(t), \mathbf{h}_j(t-1)\right]}{\pi_{\theta j,old}\left[\Delta \mathbf{u}_j(t) \mid o_j(t), \mathbf{h}_j(t-1)\right]} \tag{58}$$

$$L(\theta_j) = -\mathbb{E}\left\{\min\left[\rho_j(t)\cdot \hat{A}(t), \mathrm{Clip}\left(\rho_j(t), 1-\epsilon, 1+\epsilon\right)\cdot \hat{A}(t)\right]\right\} - \zeta \cdot \mathcal{H}(\pi_{\theta_j}) \tag{59}$$

$$L(\phi) = \mathbb{E}_t\left\{\max\left[\left(V_\phi(\mathcal{S}(t)) - R_t^{\text{target}}\right)^2, \left(\mathrm{Clip}(V_\phi(\mathcal{S}(t))) - R_t^{\text{target}}\right)^2\right]\right\} \tag{60}$$

where $\epsilon$ is the policy clipping threshold; $\mathcal{H}(\pi_{\theta_j})$ is the policy entropy regularization term, whose coefficient $\zeta$ linearly decays with the training process to maintain control activity during the early phase of exploration. The Critic loss function approximates the true return target value $R_t^{\text{target}}$ via a dual-clipped mean squared error to maintain the robustness of value function estimation over a long horizon.

## 6. Simulation experiments and performance verification

### 6.1. Experimental setup

To systematically validate the fidelity and control efficacy of the proposed macroscopic traffic flow model (L-METANET), crash risk assessment model (XGBoost-SHAP-RPBL), and hierarchical coordinated control framework (MPC-STMAPPO) under realistic and complex conditions, a high-fidelity macro-micro fused traffic simulation platform was constructed. The platform replicates an 18 km network topology using multi-source detector data from the Eastern Expressway in Changchun. This section details the experimental environment across three dimensions: simulation platform development, multi-source heterogeneous data composition, and hierarchical control parameter configurations.

**Simulation platform development and network topology replication:** The simulation environment is built on the open-source microscopic platform SUMO. To circumvent the TCP latency of standard TraCI and support high-frequency MARL interactions, the Python control algorithms interface directly with the SUMO engine via the libsumo C++ API. As illustrated in Fig. 6, the network replicates the topology of the 18 km Changchun Eastern Expressway at a 1:1 scale, encompassing mainline segments, on/off ramps, and weaving bottlenecks. The testbed features 9 VSL gantries and 3 RMs. Virtual loop (E1) and area (E2) detectors are positioned across all critical nodes (mainline origins, ramps) and within individual segment-lane cells to ensure full-state network perception.

**Multi-source heterogeneous data and scenario configurations:** empirical multi-source heterogeneous traffic data drive the calibration and evaluation tasks across different modules through the following data pipelines: 1) L-METANET Calibration (Section 6.2): High-frequency corridor volumes collected from 06:00 to 19:00 on November 10, 2022, covering all mainline cross-sections and ramps, establish the dynamic boundary conditions in SUMO (illustrated in **Appendix B**). This demand is directly input into SUMO as boundary conditions. Under conditions without management intervention, the natural driving data output by the SUMO are extracted as the ground truth to calibrate the L-METANET parameters; 2) Weaving segment risk assessment model fitting (Section 6.3): High-resolution vehicle trajectories extracted from UAV aerial videos at weaving segments WS3, WS4, and WS6 are paired with the STRF model to generate crash conflict labels, providing the microscopic data baseline to train the macroscopic XGBoost-SHAP-RPBL risk prediction model; 3) Hierarchical coordinated control framework performance evaluation (Section 6.4): The identical 13-hour empirical flow acts as the network boundary demand. Operational metrics under various control algorithms are evaluated to quantify the control efficacy of the proposed framework; 4) Transferability analysis (Section 0): The 13-hour baseline flow is adapted into five distinct synthetic demand scenarios to evaluate the framework's generalization resilience under unobserved traffic variations.

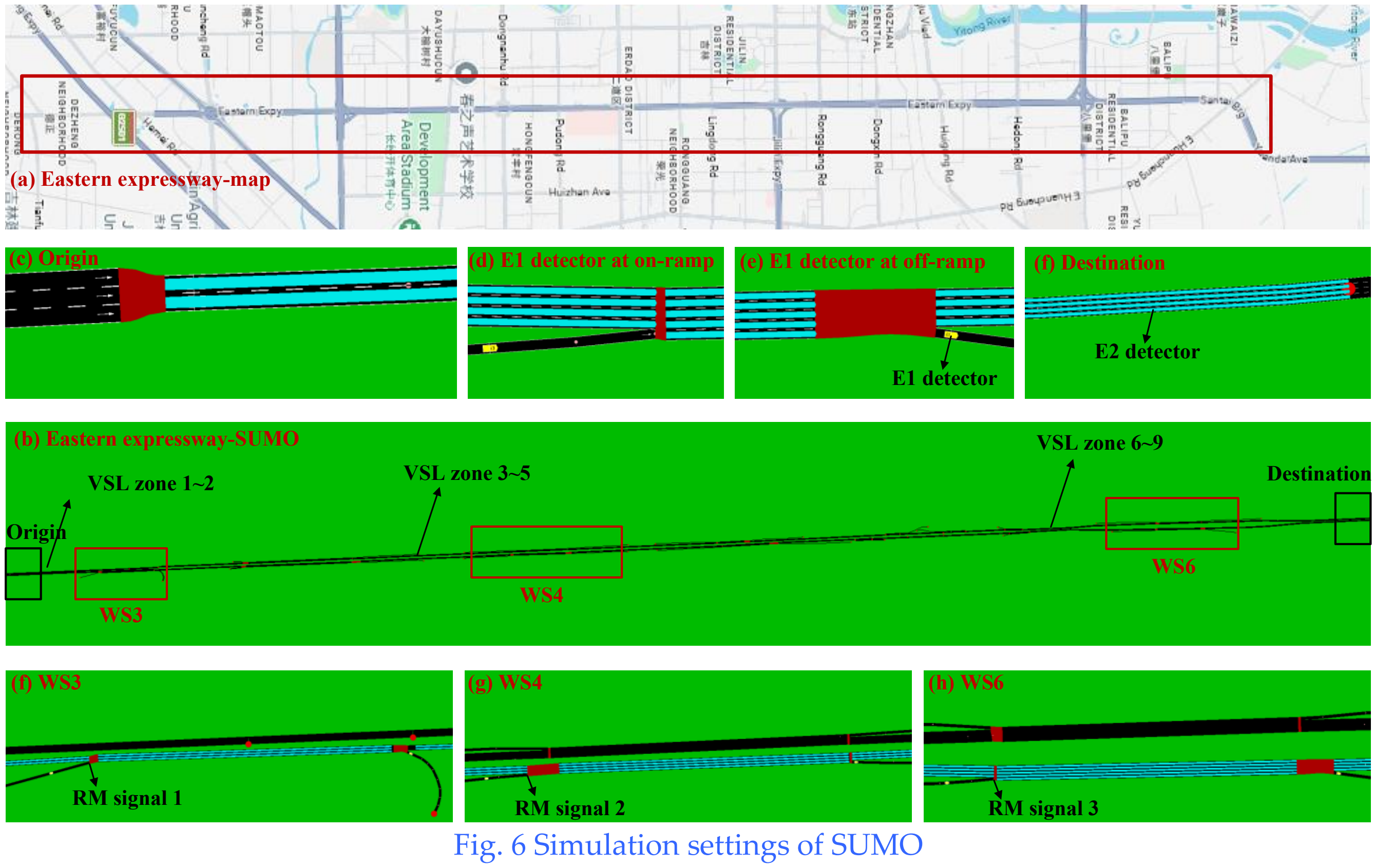

Fig. 6 Simulation settings of SUMO

**Hierarchical time scales and hyperparameter configurations:** (1) To simulate a mixed traffic environment during the transition toward CAVs, the CAV penetration rate is set to $\alpha = 0.5$, and the spontaneous compliance rate of HDVs is set to $\beta = 0.8$. (2) The simulation step of the underlying SUMO engine is $0.5s$, while the prediction step of L-METANET is set to $\Delta T = 5s$. Within the hierarchical control framework, the high-frequency sampling and action execution interval of the lower-level ST-MAPPO agents is $T_s = 120s$, and the update cycle of the upper-level MPC is $T_c = 600s$ with a rolling prediction horizon of $N_p = 3$. (3) To guarantee system safety and prevent the MARL agents from issuing catastrophic commands during the exploration phase, multi-dimensional physical constraints and penalty mechanisms are imposed on the action space: The absolute speed limit range of the VSL is restricted to [40, 80] km/h, with the maximum speed jump between adjacent control steps limited to $\Delta v_{step} = 20$ km/h. The metering rate of the RM is bounded within $[0.1, 1.0]$, with a maximum step change of $\Delta r_{step} = 0.3$. Determined through extensive preliminary experiments, the weight allocation for the cost and reward functions serves a dual purpose. First, it normalizes the disparate physical units and mathematical scales, preventing any single indicator from dominating the optimization. Second, it balances the competing objectives of efficiency and safety to ensure control decisions converge toward a Pareto optimal point. Ultimately, $\omega_1, \omega_2, \omega_3$ are set to a 1:5:0.2 ratio. Meanwhile, the global feedback reward is uniformly scaled to the order of $[-1,1]$ before being fed into the network, preventing gradient explosion in the Critic network when encountering severe penalties. To address the variance collapse problem common to the PPO algorithm in continuous control tasks, the log of the initial standard deviation output by the Actor network of ST-MAPPO is restricted to -1.0 to suppress abrupt action fluctuations during early exploration. The initial learning rates for the Actor and Critic networks are set to $6\times10^{-4}$ and $1\times10^{-4}$, respectively, both of which employ a linear decay mechanism toward a final value of $1\times10^{-5}$.

6.2. Calibration and validation of L-METANET

To validate the fidelity of the proposed L-METANET model in real-world weaving segments, an empirical analysis is conducted using high-frequency detector data from the Changchun expressway network. Real-time flow from the Eastern Expressway serves as dynamic boundary conditions in the SUMO simulation environment. Data generated by virtual detectors within SUMO are subsequently utilized to calibrate both the L-METANET and METANET models through the following steps:

**Step 1:** Based on a road network exported from OpenStreetMap, the Eastern Expressway network is fully replicated in SUMO at a 1:1 scale. This process strictly maps the physical topology, including lane configurations, weaving segment lengths, and the geometric characteristics of entrances and exits.

**Step 2:** The Eastern Expressway is partitioned into multiple macroscopic physical segments and specific lanes. The partitioning criteria dictate: (1) ensuring absolute geometric homogeneity within each sub-segment, meaning no lane additions/drops or ramp variations occur internally; (2) establishing independent physical boundary cut-off points at critical nodes featuring on-ramp merging, off-ramp diverging, or weaving bottlenecks; and (3) strictly constraining and uniforming the lower bound of physical segment lengths to satisfy the Courant–Friedrichs–Lewy numerical stability condition.

**Step 3:** Virtual detectors are systematically deployed. E2 detectors are deployed within each discretized segment-lane cell to gather lane-level density, speed, and flow data. Concurrently, E1 detectors are installed at all mainline origin/destination points and ramp terminals to record the vehicle arrival and departure flow rates for boundary conditions.

**Step 4:** The obtained real-time traffic volumes are input as boundary demands into the SUMO simulator to model the spatiotemporal evolution of the mixed traffic flow within the weaving network.

**Step 5:** Data collection and cleaning are performed on the E1 and E2 detector outputs. Vehicle trajectories are aggregated into macroscopic statistical indicators within discrete time steps, serving as the boundary input drivers and calibration ground truth for both L-METANET and conventional METANET.

**Step 6:** Prediction models for L-METANET and METANET are executed independently. Using the cleaned E1 data as physical boundary conditions and the E2 data as fitting ground truths, the parameter sets for both models are strictly calibrated utilizing the GA-VNS algorithm detailed in Section 4.2.1. The calibrated parameter values for the L-METANET and METANET are summarized in Table 1.

**Step 7:** Post-calibration, 5-minute aggregated flow and speed data are randomly sampled across different segment-lanes. The outputs of both prediction models are compared against the SUMO ground truth to systematically verify the validity and high fidelity of the constructed models.

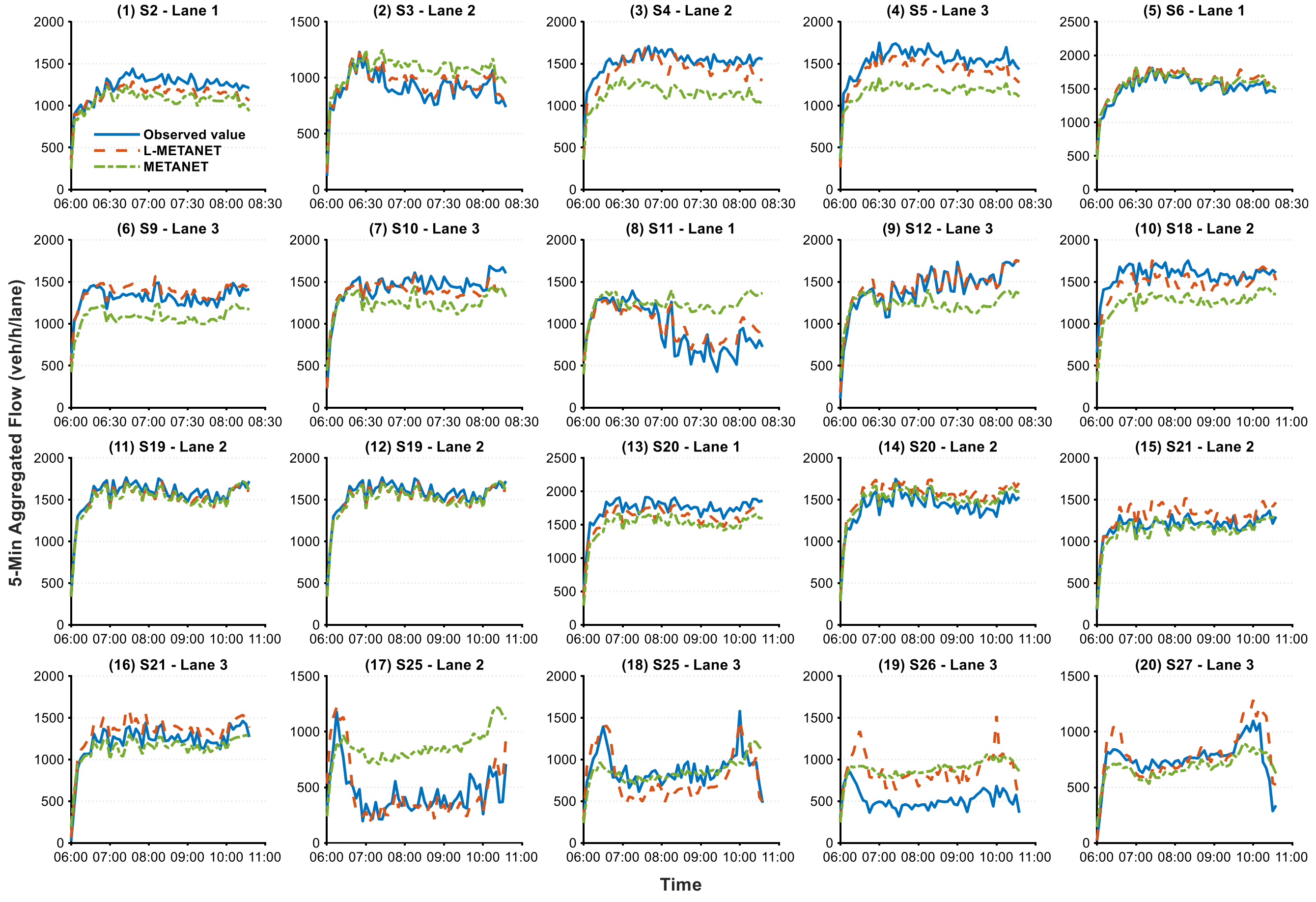

Fig. 7 Comparison of simulated flow for different macroscopic traffic flow models

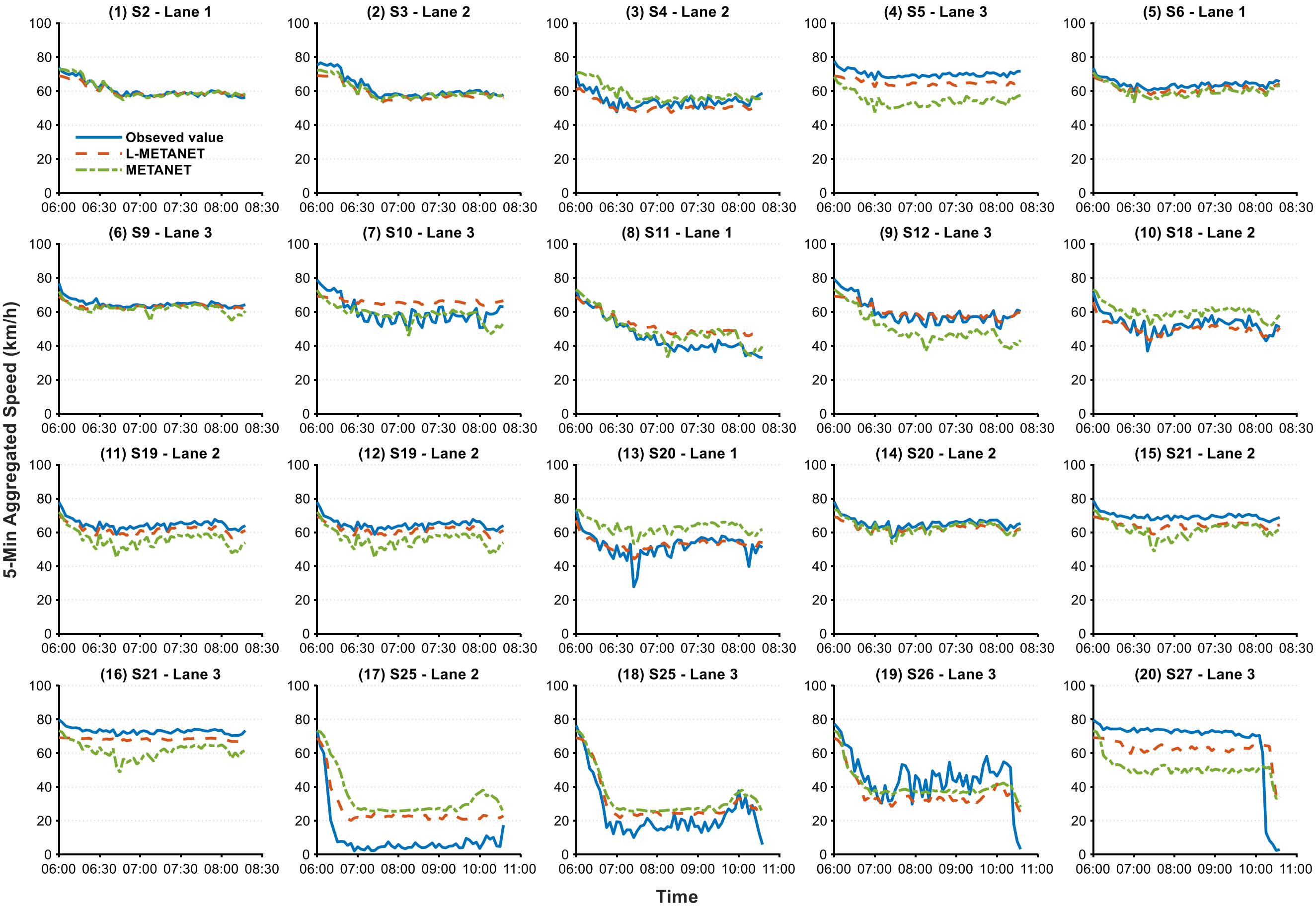

Fig. 8 Comparison of simulated speed for different macroscopic traffic flow models

Table 1 Calibrated parameter values of L-METANET and METANET

| Model | $\tau$ | $\eta$ | $\kappa$ | $\delta^{in}$ | $\delta^{out}$ | $\mu_{CAV}$ | $\mu_{HDV}$ | $v_f$ | $k_{cr}$ | $a$ |
|---|---|---|---|---|---|---|---|---|---|---|
| L-METANET | 0.0277 | 12.78 | 13.8 | 0.59 | 0.02 | 0.52 | 0.32 | 74.17 | 41.59 | 1.49 |
| METANET | 0.0200 | 28.54 | 21.30 | - | - | - | - | 72.89 | 22.3 | 3.37 |

Table 2 Average error of L-METANET and METANET

| Model | Indicator | Speed | | Flow | |
|---|---|---|---|---|---|
| | | Lane 1 | Other lanes | Lane 1 | Other lanes |
| L-METANET | RMSE | 4.25 | 8.98 | 158 | 133 |
| | MAPE | 6.7% | 19.4% | 12.0% | 12.1% |
| METANET | RMSE | 7.04 | 11.51 | 253 | 256 |
| | MAPE | 10.8% | 27.5% | 19.5% | 24.9% |

Fig. 7 and Fig. 8 show the comparison results, and Table 2 summarizes the average errors of the two models, from which several key conclusions can be drawn:

(1) Flow evolution (Fig. 7): While both models generally track the macroscopic fluctuations of the total network throughput, L-METANET exhibits a significantly more precise dynamic response during peak-to-valley transition periods. This advantage is particularly evident in outer weaving lanes heavily impacted by ramp merging and diverging. Due to its cross-sectional homogeneity assumption, conventional METANET uniformly distributes traffic across lanes, thereby passively smoothing localized lane-level lateral flow variations. As shown in Fig. 7(17) and Fig. 7(18), METANET generates identical flow profiles for Lane 2 and Lane 3 of Segment 25, inevitably leading to distortion. In contrast, L-METANET accurately reproduces the lateral flow redistribution and sudden capacity drops induced by forced lane-changing.

(2) Speed evolution (Fig. 8): The structural theoretical advantages of L-METANET are further confirmed. During morning and evening peak hours within the weaving segments, frequent merging and diverging maneuvers induce severe lateral asymmetric friction. This phenomenon triggers sharp speed drops in the outer and adjacent lanes, accompanied by high-frequency stop-and-go shock waves. As illustrated in Fig. 8, because conventional METANET omits lane heterogeneity and lateral vehicle physics, it misestimates traffic platoon speeds under varying traffic states, which is highly visible in regions characterized by volatile temporal fluctuations, such as Fig. 8(13), (17), (18), and (20). Conversely, the L-METANET model relying on space allocation mechanisms for free and forced lane-changing and asymmetric speed penalties, accurately replicates the nonlinear dynamic transition from free-flow states to congestion shock waves. L-METANET achieves high alignment with the ground truth across both the disturbance-resistant, high-speed operations of inner lanes and the complex turbulent flows of outer lanes.

(3) Quantitative evaluation of average errors (Table 2): Across all evaluation metrics, the prediction accuracy of L-METANET consistently outperforms the METANET. Specifically, regarding speed prediction, L-METANET yields Mean Absolute Percentage Errors (MAPEs) as low as 6.7% for Lane 1 and 19.4% for other lanes, whereas the conventional METANET exhibits significantly higher errors at 10.8% and 27.5%, respectively. For flow prediction, this accuracy advantage is particularly prominent in other lanes, where the flow MAPE drops sharply from 24.9% under METANET to 12.1% under L-METANET, while the corresponding RMSE decreases from 256 to 133, nearly halving the modeling error. These statistical results strongly demonstrate the critical role of explicitly incorporating lane heterogeneity and lane-changing space allocation mechanisms, which effectively eliminates conventional cross-sectional homogeneity distortions and substantially enhances lane-level traffic flow modeling precision.

These results demonstrate that L-METANET possesses a high-fidelity capability to capture complex, nonlinear physical phenomena of multi-stream weaving alongside sudden capacity drops. This establishes a reliable predictive foundation for subsequent closed-loop coordinated VSL-RM control.

### 6.3. Fitting and performance verification of XGBoost-SHAP-RPBL

Following the steps described in Section 4.2, the analytical equations for risk assessment and prediction across all tasks were derived, as detailed in **Appendix C**. The underlying data originate from our previous study involving 14 tasks across 7 weaving segments for both merging and diverging risk types (Ma et al., 2026a). This comprehensive dataset encompasses the six specific tasks evaluated in this study: WS3-M, WS3-D, WS4-M, WS4-D, WS6-M, and WS6-D.

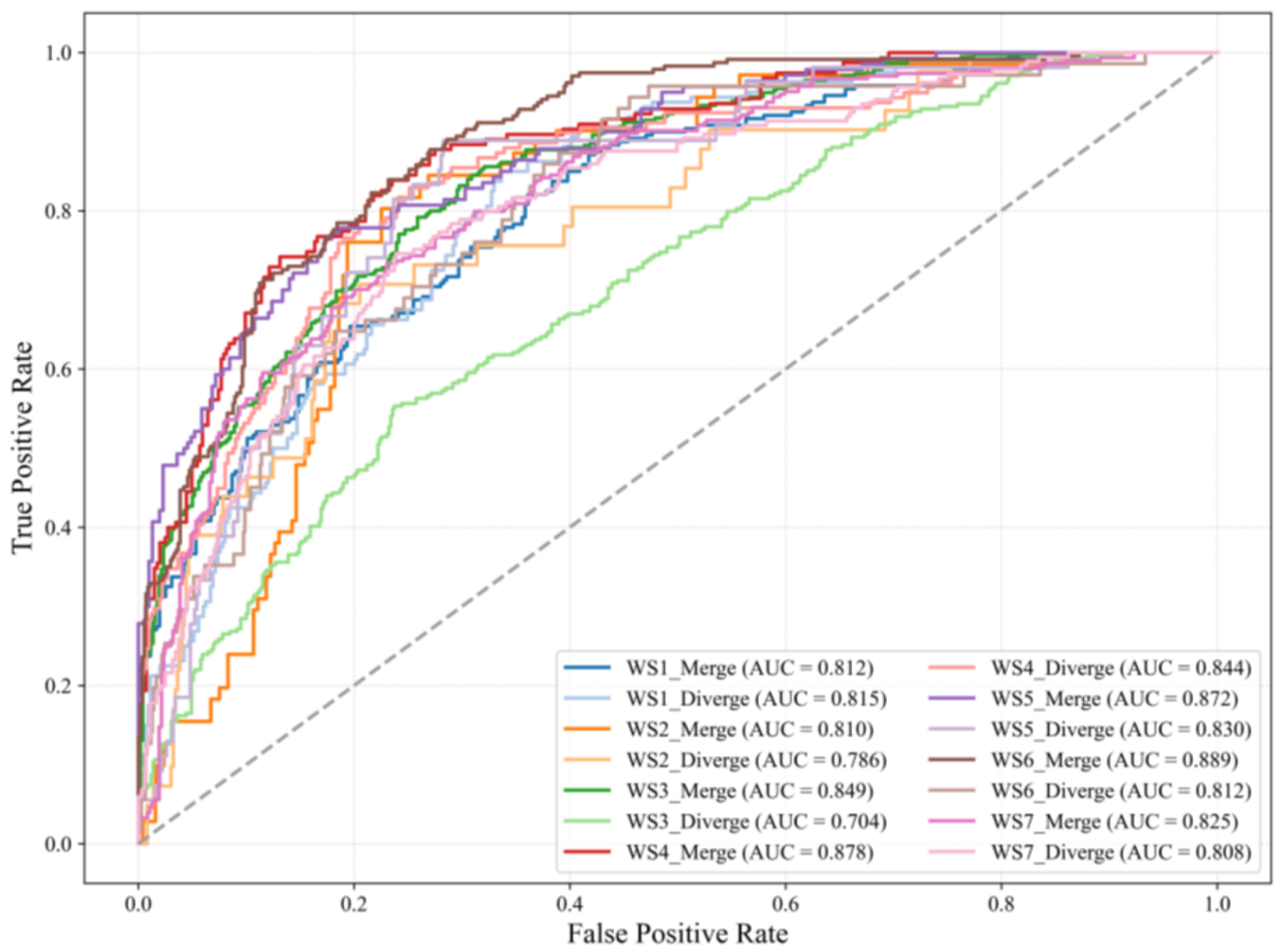

Fig. 9 ROC curves of the XGBoost-SHAP-RPBL

Fig. 9 illustrates the ROC curves of the proposed model across the 14 risk identification tasks. All tasks achieve an AUC above 0.70, demonstrating robust classification capabilities. Except for a few tasks (WS1_D and WS3_D), the AUC values exceed 0.80, with the merging risk task at WS6 reaching a peak AUC of 0.889. This performance confirms that the core feature set selected via XGBoost-SHAP accurately captures the deep non-linear mapping mechanisms between macroscopic traffic flow parameters and microscopic conflict risks, thereby ensuring high predictive accuracy and generalization robustness.

To further evaluate the performance of the XGBoost-SHAP-RPBL framework, comparative experiments were conducted against multiple baselines, as illustrated in Fig. 10. In addition to the proposed model, five representative baselines were introduced for comprehensive comparison: (1) A hybrid model considering only fixed effects (XGBoost-SHAP-FPBL); (2) A random parameters model with features selected based on Pearson correlation (Pearson-RPBL); (3) A random parameters model with features selected based on Spearman correlation (Spearman-RPBL); (4) A random parameters model retaining the full set of candidate features (All Params-RPBL); (5) A baseline model with randomly selected feature inputs (Random-RPBL). Four representative tasks (WS1-M, WS3-M, WS6-D, and WS7-M) were randomly selected for analysis. The comparison of their ROC curves demonstrates the following insights:

(1) Superiority of the XGBoost-SHAP interpretability framework. The proposed XGBoost-SHAP-RPBL model significantly outperforms the Pearson-RPBL, Spearman-RPBL, and Random-RPBL models. Unlike conventional feature selection methods that rely on linear correlation assumptions, the SHAP-based interpretability framework captures high-order synergetic interactions among variables, effectively eliminating the interference of redundant features and severe multicollinearity. This feature extraction mechanism enables the model to achieve high predictive accuracy within a highly parsimonious feature space that retains only key physical variables.

(2) Capability of random parameters to capture unobserved heterogeneity. The predictive accuracy of the proposed model slightly surpasses that of the XGBoost-SHAP-FPBL baseline. Although the performance gain is incremental, it underscores the necessity of accounting for unobserved heterogeneity in traffic safety modeling. By introducing normally distributed random parameters, the RPBL model effectively absorbs utility fluctuations stemming from individual driver differences, microscopic behavioral stochasticity, and environmental noise, thereby yielding superior probabilistic fitting performance from a statistical perspective.

(3) Dynamic trade-off between accuracy and controllability: Although the proposed model's absolute accuracy is slightly lower than the All Params-RPBL baseline in certain scenarios, its predictive performance remains highly acceptable. Crucially, the proposed framework exhibits superior practical engineering value. The baseline incorporates dozens of candidate variables, causing it to suffer from the curse of dimensionality, an impractical data collection threshold, and slow online optimization convergence. Conversely, the proposed framework leverages SHAP selection to compress input dimensionality to single digits, enabling real-time risk estimation and closed-loop control feedback.

In summary, the XGBoost-SHAP-RPBL model optimally balances risk accuracy with operational controllability. This establishes a rigorous foundation for embedding the analytical risk function into the hierarchical, closed-loop MPC-MARL framework.

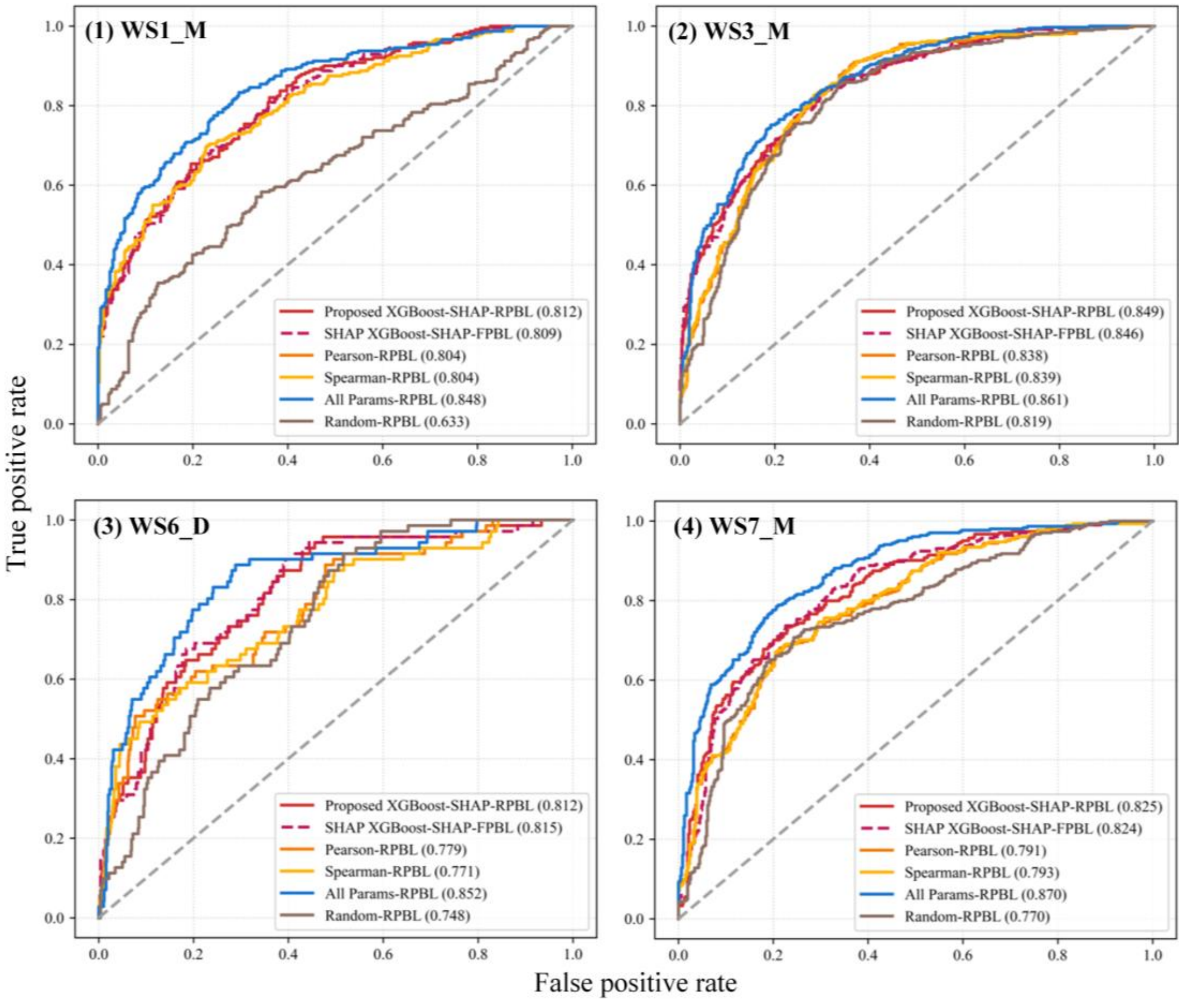

Fig. 10 Comparison of ROC curves for different algorithms

## 6.4. Performance verification of the MPC-STMAPPO hierarchical coordinated control framework

### *6.4.1 Baseline control strategies*

To implement a rigorous controlled variable analysis, this study designs 11 baseline comparison schemes. These strategies comprehensively cover various architectural paradigms, including pure model-driven, pure data-driven, hybrid model-data-driven, and other reference frameworks:

**Type 1: Pure model-driven frameworks.** This category comprises three control schemes: (1) MPC + L-METANET + Lane-level (MPC-LL): Utilizes the L-METANET model and an MPC controller to execute LVSL and RM. This represents the proposed upper-level baseline MPC controller but operates without the lower-level MARL fine-tuning; (2) MPC + L-METANET + Segment-level (MPC-LS): Utilizes the L-METANET model and an MPC controller to execute segment-level VSL and RM. Lane-level predicted states from L-METANET are aggregated into segment-level metrics before executing segment-level VSL control; (3) MPC + METANET + Segment-level (MPC-SS): Utilizes the conventional METANET model and an MPC controller to execute segment-level VSL and RM.

**Type 2: Pure data-driven frameworks.** This category comprises three MARL schemes: (1) Pure ST-MAPPO: Implements MARL using the proposed lower-level MARL residual compensator, operating independently without upper-level MPC guidance; (2) Pure MAPPO: Employs the standard MAPPO algorithm, where the proposed Mamba-embedded policy network and GCTA-embedded Critic network are replaced by conventional multilayer perceptron (MLP) architectures; (3) Pure MADDPG: Implements the multi-agent deep deterministic policy gradient (MADDPG) algorithm to analyze the influence of alternative MARL paradigms on training efficacy and control results.

**Type 3: Hybrid model-data-driven frameworks.** This category evaluates different combinations within the hierarchical structure: (1) Proposed MPC + ST-MAPPO Framework (MPC-STMAPPO): The integrated framework features the upper-level MPC baseline controller and the lower-level ST-MAPPO residual compensator; (2) MPC + MAPPO (MPC-MAPPO): Combines the proposed MPC baseline controller with a standard MAPPO-based residual compensator; (3) MPC + MADDPG (MPC-MADDPG): Combines the proposed MPC baseline controller with a MADDPG-based residual compensator.

**Type 4: Baseline benchmarks.** This category includes two strategies: (1) No-control: Vehicles navigate naturally within SUMO based on empirical boundary inputs. VSL and RM values are maintained at their maximum allowable limits; (2) Random control: VSL and RM actuation commands are generated randomly within their respective physical upper and lower bounds.

### *6.4.2 Training process and convergence mechanism analysis*

To evaluate the internal evolution, training efficiency, and robustness of the proposed hierarchical coordinated ATM framework, this section quantitatively compares training convergence across different control architectures (Fig. 11). Based on preliminary trials, the maximum training duration is set to 500 episodes for the pure data-driven group to accommodate its slower convergence, and 200 episodes for the hybrid model-data-driven group due to MPC-guided acceleration. Although pure model-driven frameworks require no training, they are evaluated using the identical reward system by embedding an inactive MARL layer (yielding zero residual outputs) to ensure a consistent baseline. The resulting reward curves yield the following insights:

(1) Pure data-driven group: As shown in Fig. 11(a), without upper-level MPC guidance, the proposed ST-MAPPO algorithm significantly outperforms baseline MARL methods over the 500-episode horizon. Despite the high-dimensional spatiotemporal state space and non-stationary disturbances across multiple weaving segments, ST-MAPPO's total reward climbs rapidly from approximately -46 during the initial phase (0–100 episodes) and converges around episode 280, stabilizing near -25. Conversely, standard MAPPO traps in a local optimum near -33 due to the lack of spatiotemporal attention, while MADDPG exhibits severe global oscillations (fluctuating between -25 and -48) and fails to converge. This confirms the efficacy of the network architecture designed in Section 5.3.4: the embedded Mamba selective state-space model filters high-frequency non-Markovian noise to preserve long-sequence temporal context, while the GCTA-enhanced Critic network injects a strong spatial inductive bias that mitigates competitive agent behaviors during decentralized execution, steering the policy toward the global Pareto front.

(2) Pure model-driven group: Serving as static baselines derived from deterministic rolling optimization (Fig. 11(b)), the proposed MPC-LL yields the highest average reward (-29.81), outperforming the MPC-LS (-32.29) and MPC-SS (-36.58). Conventional METANET's cross-sectional homogeneity assumption fails under highly heterogeneous traffic—such as outer-lane turbulence from forced cut-ins paired with free-flowing inner lanes—where its uniform segment-level commands waste inner-lane capacity and jeopardize outer-lane safety. In contrast, the L-METANET explicitly models free/forced lane-changing mechanisms and lateral asymmetric friction penalties within the momentum relaxation equation. This high-fidelity forecasting of lane heterogeneity enables fine-grained, differentiated LVSL, boosting system performance.

(3) Hierarchical coordinated control group: Fig. 11(c) evaluates the training curves over a 200-episode horizon. Incorporating upper-level MPC baseline commands as guiding anchors for rolling fine-tuning substantially lifts initial rewards and accelerates learning. The proposed MPC-STMAPPO framework converges within just 90 episodes, stabilizing securely at approximately -26. Conversely, MPC-MAPPO stagnates in a local optimum, and MPC-MADDPG exhibits severe global oscillations. This indicates that while physical boundary constraints contract the invalid exploration space during blind exploration, a well-matched MARL architecture remains crucial for effective model-data synergy.

(4) Ablation analysis: Fig. 11(d) presents the ablation study for the proposed framework. MPC-STMAPPO exhibits a higher initial reward, faster convergence, and reduced volatility compared to ST-MAPPO, confirming that MPC guidance ensures stable policy behavior and superior training efficiency. However, pure ST-MAPPO slightly outperforms MPC-STMAPPO in final steady-state reward by a margin of approximately 1. This numerical discrepancy reveals a slight detriment to MARL optimality caused by predictive model mismatch: since actual control commands must anchor to the upper-level baseline, the limited step-by-step residual fine-tuning margin struggles to completely overcome the systemic bias introduced by the MPC due to the spatial constraints of the truncation protection mechanism, causing a minor compromise in the ultimate reward of the hybrid framework. Nonetheless, pure ST-MAPPO's slight edge comes at the expense of severe reward oscillations during the first 150 episodes, complicating stable training. More importantly, the hierarchical framework trades this slight system performance for structural generalization resilience. Pure data-driven methods risk overfitting to specific scenarios, whereas the physical conservation laws embedded in the hierarchical architecture guarantee a rigid performance lower bound. As detailed in the Section 0, under varying loads or sudden demand mutations, pure MARL suffers from variance collapse or action overshooting due to a lack of safety guardrails. In contrast, the hierarchical framework exhibits robust transfer performance, proving its viability for real-world engineering deployment.

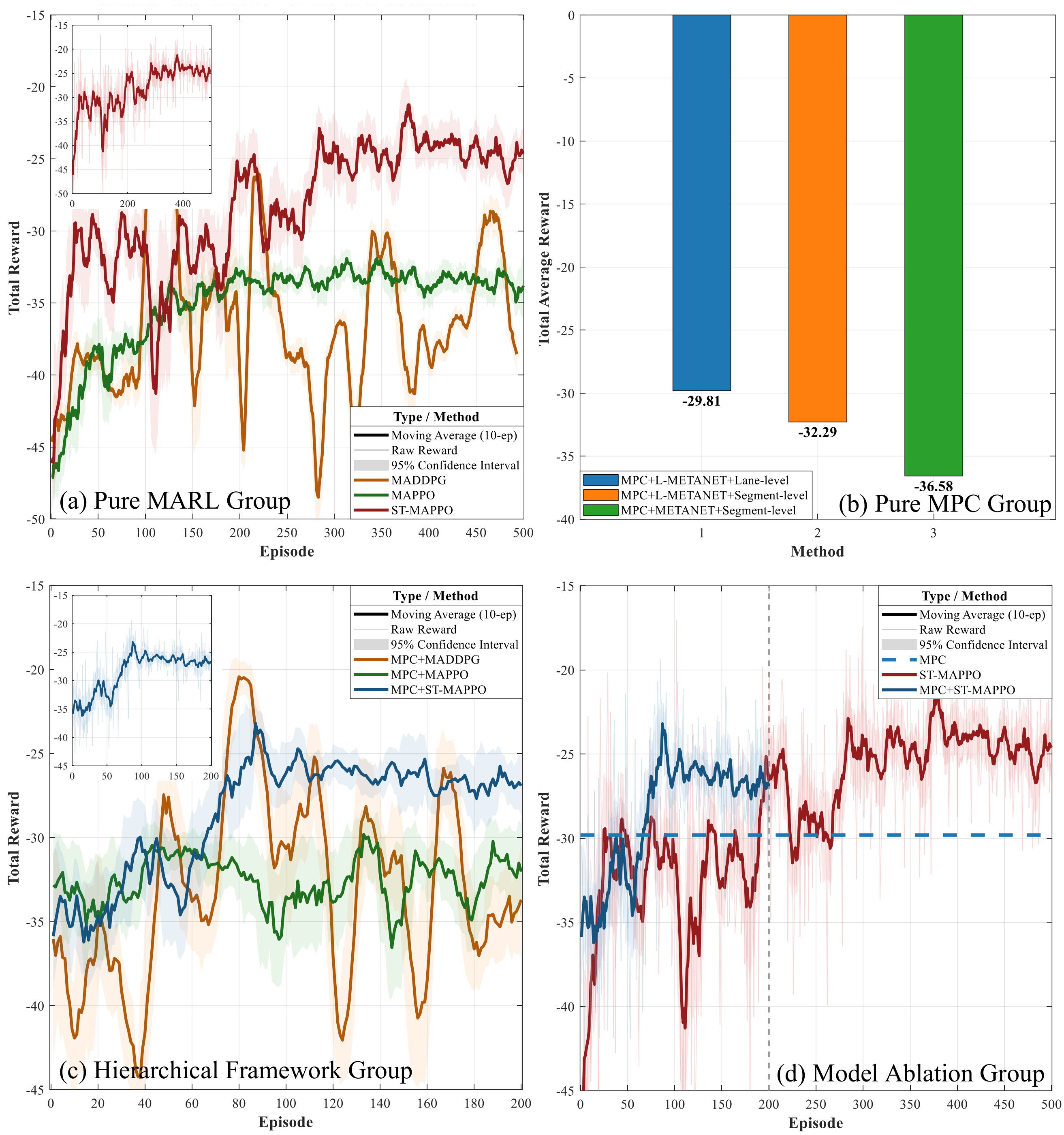

Fig. 11 Reward curves of different methods

*6.4.3 Quantitative evaluation of multi-objective control efficacy*

Table 3 provides a comprehensive quantitative evaluation of the operational efficiency and macroscopic safety multi-objective metrics across 11 control schemes under empirical dynamic traffic demands for three core weaving segments (WS3, WS4, and WS6). These metrics include mainline average speed (Speed), on-ramp queue length (Queue), macroscopic weaving segment risk (Risk), total time loss (TL), and the total number of simulated collisions (Collisions). The following key observations can be made:

(1) Across the entire experimental dataset, all ATM strategies demonstrate varying degrees of system optimization compared to the No-Control and Random Control benchmarks. Under the unmanaged natural driving state (No-Control), intense lateral disruptions among the three weaving streams at the bottleneck cause the mainline average speed at the downstream weaving segment WS6 to plummet to 14.4 km/h, escalating the total time loss (TL) to 2,947 veh·h and triggering 312 network-wide micro-collisions, which pushes the system into widespread breakdown and severe safety crises. In contrast, the proposed MPC-STMAPPO hierarchical coordinated control framework exhibits exceptional global coordination and proactive traffic state reconstruction capabilities. It slashes the total network-wide collisions to 25 (a 91.99% reduction relative to No-Control), lifts the mainline speed at the heavily congested recurrent bottleneck WS6 to 38.5 km/h, and suppresses the corresponding TL to 627 veh·h. These results validate that a coordinated VSL-RM control strategy at weaving bottlenecks effectively achieves bi-objective optimization for both traffic safety and operational efficiency;

(2) To analyze the impacts of control granularity and predictive model fidelity on the closed-loop system, the pure model-driven group (MPC-SS, MPC-LS, and MPC-LL) provides clear empirical evidence. Due to the cross-sectional lane homogeneity assumption, the conventional segment-level MPC-SS model fails to capture the lateral friction heterogeneity within the lane cross-sections induced by diverging and lane-changing. Consequently, its uniform control commands trigger severe secondary intervention disruptions. Although the L-METANET segment-level controller (MPC-LS) outperforms MPC-SS, its coarse-grained execution limits further performance gains. Conversely, when the control granularity is refined down to the lane level via the developed L-METANET model (MPC-LL), lane-level differentiated LVSL can be executed precisely to smooth traffic fluctuations. As a result, the queue length at WS4 is minimized to 0.22 veh, and the bottleneck speed at WS6 is doubled to 32.4 km/h. This strongly substantiates the necessity of reconstructing the macroscopic traffic flow model at the lane level and applying fine-grained lane-by-lane control, as established in Section 4.1;

(3) Regarding the performance of the pure data-driven groups (MADDPG, MAPPO, and ST-MAPPO), the ST-MAPPO leverages its embedded Mamba selective state-space memory and GCTA mechanism, demonstrating robust capabilities in spatiotemporal feature decoupling and exploratory policy optimization. At the recurrent core bottleneck WS6—where weaving friction is most intense and congestion is frequent—ST-MAPPO maintains a mainline operating speed of 37.1 km/h and restricts the total number of micro-collisions to 30 over the full simulation cycle. In contrast, the standard MAPPO algorithm without the spatiotemporal attention mechanism suffers from degraded traffic efficiency at WS6. Furthermore, the deterministic policy gradient-based MADDPG algorithm exhibits severe competitive, self-interested behaviors and policy fitting limitations, yielding an unacceptable expected macroscopic crash risk of 1.751 at WS6, indicating an unsustainably high crash vulnerability within this weaving segment;

(4) Within the hybrid model-data-driven group (MPC-MADDPG, MPC-MAPPO, and MPC-STMAPPO), the proposed MPC-STMAPPO achieves the highest WS6 mainline operating speed (38.5 km/h) and the shortest on-ramp queue length (21.82 veh) among all 11 control schemes, while substantially reducing network-wide simulated collisions to 25. Other hybrid combinations show consistently inferior performance compared to MPC-STMAPPO.

The quantitative results in Table 3 demonstrate the superior performance of the proposed framework from a global multi-objective Pareto optimization perspective. Among the 13 network-wide core evaluation metrics, the proposed MPC-STMAPPO strategy ranks in the top three for six metrics and in the bottom three for only one. Analysis of the ranking distribution reveals that the framework's optimal metrics are concentrated in heavily congested weaving segments under severe bottleneck pressure, notably at WS6. This indicates that coupling a physical safeguard with high-frequency residual compensation yields highly robust performance under severe spatiotemporal shock waves. At the less congested bottlenecks (WS3 and WS4), MPC-STMAPPO remains highly competitive despite not being the best-performing scheme. This minor compromise reflects the network’s multi-segment coordination mechanism, which strategically introduces localized upstream sacrifices to efficiently mitigate severe downstream congestion at WS6. Additionally, the standalone upper-level MPC-LL and lower-level

Table 3 Quantitative multi-objective performance indicators of the system under different methods

| Indicator / Method | WS3 | | | | WS4 | | | | WS6 | | | | Collisions (Number) | Top 3 (Number) | Bottom 3(Number) |
|---|---|---|---|---|---|---|---|---|---|---|---|---|---|---|---|
| | Speed (km/h) | Queue (veh) | Risk | TL (veh·h) | Speed (km/h) | Queue (veh) | Risk | TL (veh·h) | Speed (km/h) | Queue (veh) | Risk | TL (veh·h) | | | |
| No Control | **64.7 ±0.0** | **0.24 ±0.00** | 0.42 ±0.00 | **563 ±0** | **42.0 ±0.0** | 2.88 ±0.00 | 0.58 ±0.00 | **776 ±0** | **14.4 ±0.0** | 31.98 ±0.00 | **1.879 ±0.00** | **2947 ±0** | **312 ±0** | 1 | 8 |
| Random Control | **63.9 ±0.4** | **6.48 ±0.38** | 0.35 ±0.05 | 305 ±43 | **64.7 ±0.8** | 14.07 ±0.85 | 0.59 ±0.15 | 512 ±132 | **14.7 ±0.2** | **32.90 ±0.25** | **1.632 ±0.14** | **1684 ±134** | **172 ±22** | 1 | 7 |
| MPC-SS | **72.3 ±0.0** | **14.91 ±0.00** | **0.62 ±0.00** | **380 ±4** | **71.8 ±0.0** | **36.11 ±0.00** | 0.92 ±0.00 | 564 ±6 | 15.4 ±0.2 | **32.84 ±0.75** | 1.556 ±0.05 | 1677 ±56 | **62 ±8** | 2 | 6 |
| MPC-LS | **72.1 ±0.2** | **14.95 ±0.03** | **0.30 ±0.02** | **256 ±20** | **68.3 ±0.3** | 14.77 ±0.64 | **0.56 ±0.02** | **486 ±20** | **14.9 ±0.2** | 31.50 ±0.54 | 1.547 ±0.04 | 1565 ±54 | 60 ±10 | 6 | 2 |
| MPC-LL | 68.2 ±0.5 | 13.24 ±0.10 | **0.27 ±0.02** | **249 ±20** | 65.3 ±0.2 | **0.22 ±0.02** | **0.50 ±0.04** | 461 ±33 | 32.4 ±0.6 | **33.23 ±0.10** | **0.791 ±0.01** | **570 ±12** | 37 ±4 | 6 | 1 |
| MADDPG | 69.7 ±0.0 | 11.65 ±0.00 | 0.46 ±0.00 | **378 ±0** | 68.1 ±0.0 | 0.24 ±0.00 | **0.32 ±0.00** | **269 ±0** | 17.9 ±0.0 | **25.38 ±0.00** | **1.751 ±0.00** | **1753 ±0** | **24 ±0** | 4 | 3 |
| MAPPO | **65.6 ±0.6** | 11.07 ±0.65 | **0.55 ±0.12** | 340 ±80 | 67.3 ±0.3 | **25.56 ±1.32** | 0.88 ±0.19 | 543 ±123 | **32.8 ±1.2** | 31.76 ±0.44 | **0.791 ±0.03** | 634 ±63 | 36 ±12 | 2 | 3 |
| ST-MAPPO | **72.3 ±0.1** | **14.91 ±0.16** | 0.52 ±0.03 | 355 ±38 | **69.6 ±0.1** | **0.21 ±0.00** | **0.97 ±0.04** | **663 ±71** | **37.1 ±1.3** | **23.26 ±1.16** | 0.893 ±0.23 | **599 ±146** | **30 ±7** | 7 | 3 |
| MPC+ MADDPG | 67.8 ±0.3 | **8.66 ±0.08** | **0.20 ±0.03** | **174 ±40** | **62.6 ±0.2** | 0.26 ±0.03 | **1.08 ±0.04** | **877 ±91** | 22.3 ±1.1 | 31.62 ±0.49 | 0.899 ±0.03 | 969 ±77 | 37 ±7 | 3 | 3 |
| MPC+ MAPPO | 67.2 ±0.5 | 11.21 ±0.42 | **0.53 ±0.10** | 310 ±47 | 68.2 ±0.4 | **25.44 ±2.00** | 0.80 ±0.20 | **477 ±132** | 31.3 ±1.2 | 30.82 ±0.72 | 1.03 ±0.07 | 686 ±41 | 32 ±9 | 1 | 2 |
| MPC+ ST-MAPPO | 71.8 ±0.2 | 14.13 ±0.25 | 0.47 ±0.03 | 316 ±40 | 66.5 ±0.7 | **0.22 ±0.09** | **0.93 ±0.05** | 626 ±77 | **38.5 ±3.3** | **21.82 ±3.08** | **0.833 ±0.27** | **627 ±168** | **25 ±6** | 6 | 1 |

Note: "SS" stands for segment model (METANET) with segment-level VSL; "LS" stands for lane model (L-METANET) with segment-level VSL; "LL" stands for lane model (L-METANET) with lane-level VSL. "TL" stands for Time Loss. Bold green text indicates the top 3 rankings; bold red text indicates the bottom 3 rankings. "X ± Y" stands for Mean ± Standard Deviation.

ST-MAPPO models also achieve high rankings across most metrics, confirming that the hierarchical coordinated architecture effectively mitigates the inherent trade-offs between operational efficiency and traffic safety across interconnected weaving segments.

6.5. Transferability and generalization performance analysis of MPC-STMAPPO

As demonstrated by the experimental results in Section 6.4.2 and 6.4.3, under the baseline historical traffic demand, the pure data-driven ST-MAPPO algorithm yields a final converged expected reward that slightly outperforms the proposed MPC-STMAPPO hierarchical coordinated framework, while displaying highly competitive optimization performance across the various microscopic quantitative metrics listed in Table 3. Consequently, superior training efficiency alone is insufficient to fully establish the advantages of the proposed hierarchical framework. This indicates that under stationary traffic demand and a fixed network topology, ST-MAPPO can thoroughly fit the empirical trajectory profiles of that specific condition. However, this single-scenario optimality relies heavily on deterministic time-series sample interactions, making the policy highly susceptible to overfitting. In real-world engineering deployments, traffic demands are inherently time-varying. To validate policy robustness against unknown perturbations or extreme conditions—and to evaluate the generalization resilience and transferability of the hierarchical framework—this section constructs diverse transfer testing scenarios. Pre-trained policy networks optimized under baseline conditions are directly deployed across four distinct dynamic demand scenarios without further fine-tuning. This setup examines the worst-case performance boundaries of the control loop under severe shifts in the solution space. These transferability experiments focus on the two top-performing control strategies: ST-MAPPO and MPC-STMAPPO.

To evaluate their multi-dimensional transfer capabilities, the dynamic traffic demand scenarios are designed as follows: (1) Baseline traffic demand (Base_Demand): Identical to the previous simulated traffic demand, serving as the benchmark reference scenario; (2) Low-load traffic demand (0.5×Base_Demand): Scales down the original time-series boundary flows uniformly by a factor of 0.5 to simulate sparse traffic states typical of late-night or early-morning periods; (3) Over-saturated traffic demand (1.5×Base_Demand): Scales up the original boundary flows uniformly by a factor of 1.5 to simulate extreme scenarios such as holiday peaks or sudden surge propagation, strictly testing the baseline capability of the controller to mitigate large-scale, network-wide traffic breakdown; (4) Random demand 1: Applies a uniform random scaling factor $U(0.3, 1.7)$ independently to the vehicle volume of each OD flow in the baseline route file. The random demand for the $k^{th}$ flow is formulated as $d_k(rand1) = \max\left[1, \left(\alpha_k \cdot d_k(base)\right)\right]$, where $\alpha_k \sim U(0.3, 1.7)$. This transformation disrupts the volume distribution and statistical characteristics of historical demand across spatiotemporal cross-sections, causing the traffic flow of each OD pair to fluctuate randomly between 30% and 170% of its original value; (5) Random demand 2: Completely removes the topological structural constraints of historical demand, independently applying a uniform random sampling scheme $U(1, 200)$ to the vehicle volume of each OD flow in the route file, denoted as the random demand $d_k(rand2) \sim U(1, 200)$ for the $k^{th}$ flow. This transformation generates the traffic volume for each flow within entirely independent and uncorrelated intervals, thereby exposing the generalization limits of the controller under completely unknown, unstructured environmental disturbances.

The results in Table 4 utilize identical metrics to Table 3. The experiments demonstrate that while ST-MAPPO achieves excellent operational performance under specific known conditions, it is susceptible to generalization failure induced by sample overfitting and control action overshooting when subjected to unseen traffic demands in open environments. Specifically, under the over-saturated traffic demand (1.5×Base_Demand), the pure MARL strategy causes the mainline average speed at the weaving bottleneck WS6 to plummet to 14.3 km/h, driving the TL up to 3,112 veh·h and triggering 109 network-wide collisions over the simulation horizon. Under Random demand 2, where the topological structure

Table 4 Transfer application performance between ST-MAPPO and MPC-STMAPPO under different traffic demands

| Demand | Method | WS3 Speed (km/h) | WS3 Queue (veh) | WS3 Risk | WS3 TL (veh·h) | WS4 Speed (km/h) | WS4 Queue (veh) | WS4 Risk | WS4 TL (veh·h) | WS6 Speed (km/h) | WS6 Queue (veh) | WS6 Risk | WS6 TL (veh·h) | Collisions (Number) |
|---|---|---|---|---|---|---|---|---|---|---|---|---|---|---|
| 0.5×Base_ Demand | STMAPPO | **74.0 ±0.1** | 14.90 ±0.12 | **0.60 ±0.03** | **210 ±10** | **73.6 ±0.1** | **0.08 ±0.00** | 0.75 ±0.03 | 261 ±10 | **70.6 ±0.1** | 27.65 ±0.30 | 0.64 ±0.02 | 262 ±8 | 4 ±1 |
| | MPC+ STMAPPO | 73.8 ±0.2 | **14.67 ±0.13** | 0.65 ±0.02 | 218 ±14 | 72.1 ±0.4 | 0.13 ±0.01 | **0.74 ±0.02** | **247 ±5** | 70.0 ±0.3 | **20.48 ±3.34** | **0.61 ±0.04** | **240 ±9** | **3 ±1** |
| Base_ Demand | STMAPPO | **72.3 ±0.1** | 14.91 ±0.16 | 0.52 ±0.03 | 355 ±38 | **69.6 ±0.1** | **0.21 ±0.00** | 0.97 ±0.04 | 663 ±71 | 37.1 ±1.3 | 23.26 ±1.16 | 0.893 ±0.23 | **599 ±146** | 30 ±7 |
| | MPC+ STMAPPO | 71.8 ±0.2 | **14.13 ±0.25** | **0.47 ±0.03** | **316 ±40** | 66.5 ±0.7 | 0.22 ±0.09 | **0.93 ±0.05** | **626 ±77** | **38.5 ±3.3** | **21.82 ±3.08** | **0.833 ±0.27** | 627 ±168 | **25 ±6** |
| 1.5×Base_ Demand | STMAPPO | **70.8 ±0.15** | 14.67 ±0.18 | **0.28 ±0.02** | **461 ±25** | 62.9 ±0.32 | 0.47 ±0.16 | 0.60 ±0.08 | 975 ±116 | 14.3 ±0.04 | 36.07 ±0.23 | 1.61 0.09 | 3112 ±247 | 109 ±18 |
| | MPC+ STMAPPO | 64.2 ±0.50 | **11.48 0.56** | 0.40 ±0.12 | 546 ±132 | **65.5 ±0.93** | **0.37 ±0.23** | **0.56 ±0.08** | **886 ±124** | **16.4 ±0.26** | **34.60 ±0.26** | **1.59 ±0.17** | **2621 ±238** | **88 ±6** |
| Random Demand 1 | STMAPPO | **72.6 ±0.15** | 14.81 ±0.06 | **0.48 ±0.06** | **301 ±61** | 66.9 ±0.41 | 0.27 ±0.04 | **0.92 ±0.05** | 602 ±78 | 30.8 ±6.32 | 32.12 ±1.23 | **0.88 ±0.13** | 651 ±119 | 33 ±8 |
| | MPC+ STMAPPO | 72.1 ±0.09 | **13.89 ±0.16** | 0.49 ±0.05 | 325 ±57 | **69.5 ±0.59** | **0.22 ±0.01** | **0.92 ±0.09** | **582 ±97** | **36.8 ±3.04** | **24.41 ±1.25** | 0.94 ±0.17 | **627 ±91** | **28 ±5** |
| Random Demand 2 | STMAPPO | 40.3 ±8.98 | 5.65 ±0.89 | 0.63 ±0.11 | 2092 ±385 | 28.5 ±5.47 | 26.9 ±2.57 | **1.13 ±0.36** | **3797 ±1326** | 12.1 ±1.35 | 24.93 ±1.46 | 1.55 ±0.30 | 5180 ±1068 | 154 ±26 |
| | MPC+ STMAPPO | **60.6 ±8.4** | **14.92 ±0.17** | **0.47 ±0.26** | **1644 ±908** | **30.4 ±4.19** | **8.40 ±0.74** | 1.17 ±0.18 | 4079 ±663 | **12.8 ±1.96** | **22.36 ±0.71** | **1.42 ±0.12** | **4288 ±456** | **134 ±14** |

Note: "TL" stands for Time Loss; Bold text indicates the best performance.

is completely decoupled and lacks historical statistical attributes, the control brittleness of the pure MARL policy becomes more pronounced. Even at the WS3 bottleneck, where upstream weaving pressure is relatively mild, the mainline speed deteriorates to 40.3 km/h. Furthermore, the TL at the critical bottleneck WS6 reaches 5,180 veh·h, causing total network collisions to climb to 154. In contrast, although MPC-STMAPPO underperforms the pure MARL strategy in a few metrics, it delivers superior overall performance, particularly when handling unseen traffic demands. Specifically, under the 1.5×Base_Demand scenario, MPC-STMAPPO successfully maintains the WS6 mainline speed at 16.4 km/h. Although this operating speed is still low, it represents a improvement over ST-MAPPO. Furthermore, network-wide collisions are restricted to 88, marking a significant reduction compared to the 109 collisions under ST-MAPPO. Under the extreme conditions of Random demand 2, the hierarchical framework reshapes the spatiotemporal safety manifold of the traffic system through proactive control, maintaining the mainline speed at WS3 at 60.6 km/h, which represents a 50.37% efficiency increase over pure MARL, and reducing delay at the critical WS6 bottleneck by 17.22% relative to the pure MARL baseline.

In summary, MPC-STMAPPO breaks the trade-off between system performance and transfer resilience under unknown conditions. Leveraging the rigid protection of the upper-level MPC, the framework establishes a physical safety lower bound and ensures high generalization capability when the traffic system encounters unseen, abrupt, and over-saturated turbulence. Consequently, the proposed framework offers superior viability for real-world engineering deployment compared to ST-MAPPO.

## 7. Conclusion

Under a hybrid model-data-driven framework, this study addresses the core scientific challenge of managing the spatiotemporal trade-off between operational efficiency and traffic safety, as well as the difficulties of online real-time adaptive coordinated control within mixed traffic flows across consecutive multiple weaving segments on urban expressways. The main research conclusions and scientific contributions are summarized as follows:

First, a high-fidelity lane-level macroscopic traffic flow prediction model, L-METANET, is established. To overcome the cross-sectional lane homogeneity assumption inherent in the conventional METANET model, L-METANET explicitly incorporates space allocation mechanisms for both free and forced lane-changing behaviors. Field validation using empirical detector data from Changchun demonstrates that, compared with the conventional segment-level METANET model, L-METANET precisely replicates the dynamic cross-lane flow redistribution and sudden capacity drops induced by forced lane-changing, showing high alignment with ground-truth flow and speed evolution. Second, a weaving segment risk assessment and prediction model is developed and validated. By integrating data-driven machine learning with rigorous statistical causal inference, an analytical assessment and prediction equation for dynamic merging and diverging crash risks is constructed using XGBoost-SHAP feature selection and a RPBL model. This hybrid risk assessment framework exhibits excellent high-risk identification capabilities, achieving an AUC greater than 0.70 across all merging and diverging risk classification tasks, with most tasks falling into a high-performance interval above 0.80, significantly outperforming various conventional Logit models. Third, the efficacy of the hierarchical coordinated control strategy in balancing system performance and resilience is comprehensively verified. The convergence speed and steady-state reward values of the proposed MPC-STMAPPO framework surpass baseline methods, rapidly stabilizing at a reward value of approximately -26 within only 90 episodes. Furthermore, the proposed strategy ranks within the top three across 6 out of 13 core evaluation metrics. It sharply reduces the total number of simulated network-wide crashes from 312 to 25, efficiently mitigates turbulent weaving flows, and proactively suppresses localized shock waves at the core weaving bottleneck WS6, effectively breaking the inherent conflict between efficiency and safety in conventional ATM. In cross-load zero-shot transfer performance tests characterizing system resilience, where pure data-driven strategies suffer from performance degradation in unknown scenarios, the proposed framework significantly outperforms pure MARL strategies in transfer robustness and worst-case performance lower-bound protection due to the rigid protection provided by the embedded physical conservation laws within the upper-level MPC, demonstrating substantial potential for real-world industrial deployment.

Despite these breakthroughs, several directions remain for future exploration and refinement. The control horizon of this study primarily focuses on macroscopic LVSL and RM. Future work will extend this framework across scales to integrate the proposed macroscopic ATM strategies with our previously developed microscopic trajectory planning for heterogeneous multi-vehicle merging in weaving segments (Ma et al., 2026c), which is expected to further enhance multi-objective system performance. Consequently, future efforts will utilize collaborative aerial photography via multi-UAV swarms to capture microscopic trajectories over longer continuous stretches, combined with finer macroscopically aggregated data, to further demonstrate the synergistic benefits of macroscopic ATM strategies and microscopic trajectory control in real-world engineering deployments.

### CRediT authorship contribution statement

**Guodong Ma:** Writing – original draft, Writing – review & editing, Software, Methodology, Conceptualization, Formal analysis. **Baofeng Sun:** Writing – review & editing, Supervision, Funding acquisition. **Wenyu Yang:** Writing – review & editing. **Zhihong Yao:** Writing – review & editing.

### Declaration of competing interest

The authors declare that they have no known competing financial interests or personal relationships that could have appeared to influence the work reported in this paper.

### Acknowledgments

This research was supported by the National Natural Science Foundation of China (grant number 52472313).

### Data availability

Data will be made available on request.

**Appendix A.** Flow estimation of free and forced lane-changing in L-METANET

**Part 1: estimation of free lane-changing flow.** In the lane-level extended dynamic equations, the lateral net free lane-changing flow $\phi^{i}_{m,n\to n'}$ (where $n'$ denotes the adjacent $n-1$ or $n+1$ lane) between adjacent lanes within the same segment serves as the core dynamic that breaks the independent operation of each lane and induces lateral friction effects. Following the multi-lane traffic flow modeling logic proposed by Papageorgiou et al., free lane-changing behaviors are primarily driven by the density differences between adjacent lanes and are strictly constrained by the remaining available space of the target lane (Roncoli et al., 2015). Incorporating the heterogeneity of mixed traffic flow into the baseline Papageorgiou macroscopic lane-changing model yields the formulation L-METANET in this study.

First, the lane-changing attractiveness of target lane $n'$ to current lane $n$ is defined as $A^{i}_{m,n\to n'}$. In a mixed traffic environment comprising CAVs and HDVs, lane-changing aggressiveness varies significantly by vehicle type. Let $\mu_{CAV}$ and $\mu_{HDV}$ denote the lane-changing aggressiveness coefficients for CAVs and HDVs, respectively. Combined with the CAV penetration rate $\alpha$, the comprehensive lane-changing aggressiveness parameter $\mu_{mix}$ for the current segment is expressed as Eq. (61). Based on $\mu_{mix}$, the mixed-flow lane-changing attractiveness $A^{i}_{m,n\to n'}$ is calculated as shown in Eq. (62). Note that the parameter $P^{i}_{m,n\to n'}$ is a spatial weighting coefficient, which typically equals 1, but can be calibrated in specific zones such as on-ramps and off-ramps to capture localized driving preferences. Based on $A^i_{m,n\to n'}$, the lateral free lane-changing demand flow $D^i_{m,n\to n'}$ migrating from cell $(m,n)$ to cell $(m,n')$ is the product of the total number of vehicles in the current lane and the attractiveness index, as formulated in Eq. (63).

However, due to the finite physical space of target lane $n'$, not all lane-changing demands can be satisfied. The maximum remaining acceptable space (i.e., supply flow rate) $S^i_{m,n'}$ that target lane $n'$ can provide at the current time step depends on its jam density $k_{jam,m,n'}$, as defined in Eq. (64).

When adjacent lanes on both sides of target lane $n'$ (i.e., $n'-1$ and $n'+1$) simultaneously generate lane-changing demands that exceed the remaining space of the target lane, the actual lane-changing flow rates must be scaled down and allocated proportionally based on demand. Therefore, the final actual net lane-changing flow $\phi^i_{m,n\to n'}$ from the current lane $n$ to target lane $n'$ is calculated via Eq. (65). This allocation mechanism ensures that under extreme congestion, the target lane does not experience density violation beyond its physical limits due to excessive vehicle admissions.

$$\mu_{mix} = \alpha \cdot \mu_{CAV} + (1-\alpha)\cdot \mu_{HDV} \tag{61}$$

$$A^i_{m,n\to n'} = \mu_{mix} \cdot \max\left[0, \frac{P^i_{m,n\to n'}\cdot k^i_{m,n} - k^i_{m,n'}}{P^i_{m,n\to n'}\cdot k^i_{m,n} + k^i_{m,n'}}\right] \tag{62}$$

$$D^i_{m,n\to n'} = A^i_{m,n\to n'} \cdot \frac{\Delta x_m}{\Delta t}\cdot k^i_{m,n} \tag{63}$$

$$S^i_{m,n'} = (k_{jam,m,n'} - k^i_{m,n'})\frac{\Delta x_m}{\Delta t} \tag{64}$$

$$\phi^{i,fr}_{m,n\to n'} = \min\left[1, \frac{S^i_{m,n'}}{D^i_{m,n'-1\to n'} + D^i_{m,n'+1\to n'}}\right]\cdot D^i_{m,n\to n'} \tag{65}$$

**Part 2: Estimation of forced lane-changing flow.** Unlike free lane-changing driven by localized density differentials, forced lane-changing is primarily governed by network topology constraints, such as lane drops and ramps, alongside the route choices of drivers. Traditional macroscopic models typically calculate forced lane-changing flows continuously within a road segment. However, under complex weaving segment topologies characterized by sudden lane drops or high-volume off-ramps, lane-level models often suffer from numerical errors such as mass leakage or duplicate counting due to the coupling between convection and friction terms. To resolve this issue, this study introduces a cascading vehicle-borrowing mechanism to precisely quantify forced lane-changing flows and their associated lateral friction penalties. Let $Q^i_m$ denote the actual arrival flow at the downstream cross-section of segment $m$ after the redistribution by internal free lane-changing, satisfying $Q^i_{m,n} = q^i_{m,n} + \phi^{i,fr}_{m,n}$. The difference $\Delta q_{m,n}$ dictates the physical direction and cascading intensity of forced lane-changing maneuvers, as formulated in Eq. (66).

$$\Delta q_{m,n} = Q^i_{m,n} - s^i_{m,n} \tag{66}$$

**Scenario A:** Flow surplus ($\Delta q_{m,n} \geq 0$) without downstream lane continuity. When the arrival flow exceeds the diverging demand, vehicles originally traveling in the diverging lane intend to continue straight. Due to downstream geometric cut-offs or narrowing, this surplus flow $\Delta q_{m,n}$ must be forced into the adjacent inner through lane (lane $n+1$). This forced cut-in maneuver imposes a lateral asymmetric friction penalty on both lanes: the inner through lane absorbs the forced vehicle insertions, whereas the outer lane experiences forced vehicle extractions, as expressed in Eq. (67).

$$\phi^{i,fo}_{m,n\to n+1} = \Delta q_{m,n} \tag{67}$$

**Scenario B:** Cascading vehicle borrowing induced by flow deficit ( $\Delta q < 0$ ). When the arrival flow fails to meet the off-ramp demand, a flow deficit $-\Delta q_{m,n}$ occurs. To satisfy the diverging demand, traffic must be forceded from the inner through lanes ( $n-1, n-2\ldots$ ) toward the outer side. The model evaluates lanes sequentially from the outermost to the innermost. Taking lane $n+1$ as an example, its actual borrowable flow $b_{m,n+1}$ is bounded by the absolute flow capacity threshold of the lane, as shown in Eq. (68). For the $(n+k)^{th}$ inner through lane, when selected for vehicle borrowing, the remaining unfulfilled deficit $E_{m,n+k}$ equals the total deficit minus the cumulative flow already extracted from all lanes to its outer side, as shown in Eq. (69). The actual borrowed flow $b_{m,n+k}$ from this lane is then determined as the minimum of the remaining deficit and the arrival flow of the lane, as formulated in Eq. (70). Crucially, vehicles cannot instantaneously teleport across multiple lanes in physical space. If an inner lane (e.g., $n+k, k \ge 1$ ) is subject to vehicle borrowing, the forced lane-changing traffic stream must sequentially traverse all intermediate lanes to reach the outer diverging lane $n$. Consequently, the actual forced lane-changing flow occurring between any two adjacent lanes $n+k$ and $n+k-1$ equals the flow directly borrowed from lane $n+k$ plus the cumulative borrowed flow originating from all lanes further inner to lane $n+k$ that must transit through lane $n+k$, as formulated in Eq. (71).

$$b_{m,n+1} = \min(-\Delta q_{m,n}, Q^{i}_{m,n+1}) \tag{68}$$

$$E_{m,n+k} = -\Delta q_{m,n} - \sum_{j=1}^{k-1} b_{m,n+j} \tag{69}$$

$$b_{m,n+k} = \min\left(E_{m,n+k}, Q^{i}_{m,n+k}\right) \tag{70}$$

$$\phi^{i,fo}_{m,n+k\to n+k-1} = \sum_{j=n+k}^{N_m} b_{m,j} \tag{71}$$

**Scenario C:** Other conditions. For situations outside Scenarios A and B, no forced lane-changing behavior is triggered, as defined in Eq. (72).

$$\phi^{i,fo}_{m,n\to n+1} = 0, \forall n \tag{72}$$

### Appendix B. Traffic demand of the Eastern Expressway

The boundary traffic demands for the SUMO simulation are derived from empirical traffic volumes collected between 06:00 and 19:00 on November 10, 2022, across all mainline cross-sections and ramps along the Eastern Expressway. Inflow demands are illustrated in Fig. 12, where Fig. 12(1) denotes the mainline origin volume and the remaining subfigures show the on-ramp demands. Similarly, Fig. 13 depicts outflow profiles, with Fig. 13(10) representing the mainline destination demand and the remaining subfigures displaying the off-ramp demands.

### Appendix C. Identification of key variables and RPBL fitting results

Table 5 lists the key feature variables extracted via the XGBoost-SHAP framework, and Table 6 presents the estimated parameter values calibrated using the RPBL model.

## Appendix D. Data used in this study

The micro-trajectory data used has been fully made publicly available at the following address: https://huggingface.co/datasets/InterestingITS/U-EASWS/tree/main. The macro-traffic demand data used can be obtained by contacting the first author at: magd22@mails.jlu.edu.cn.

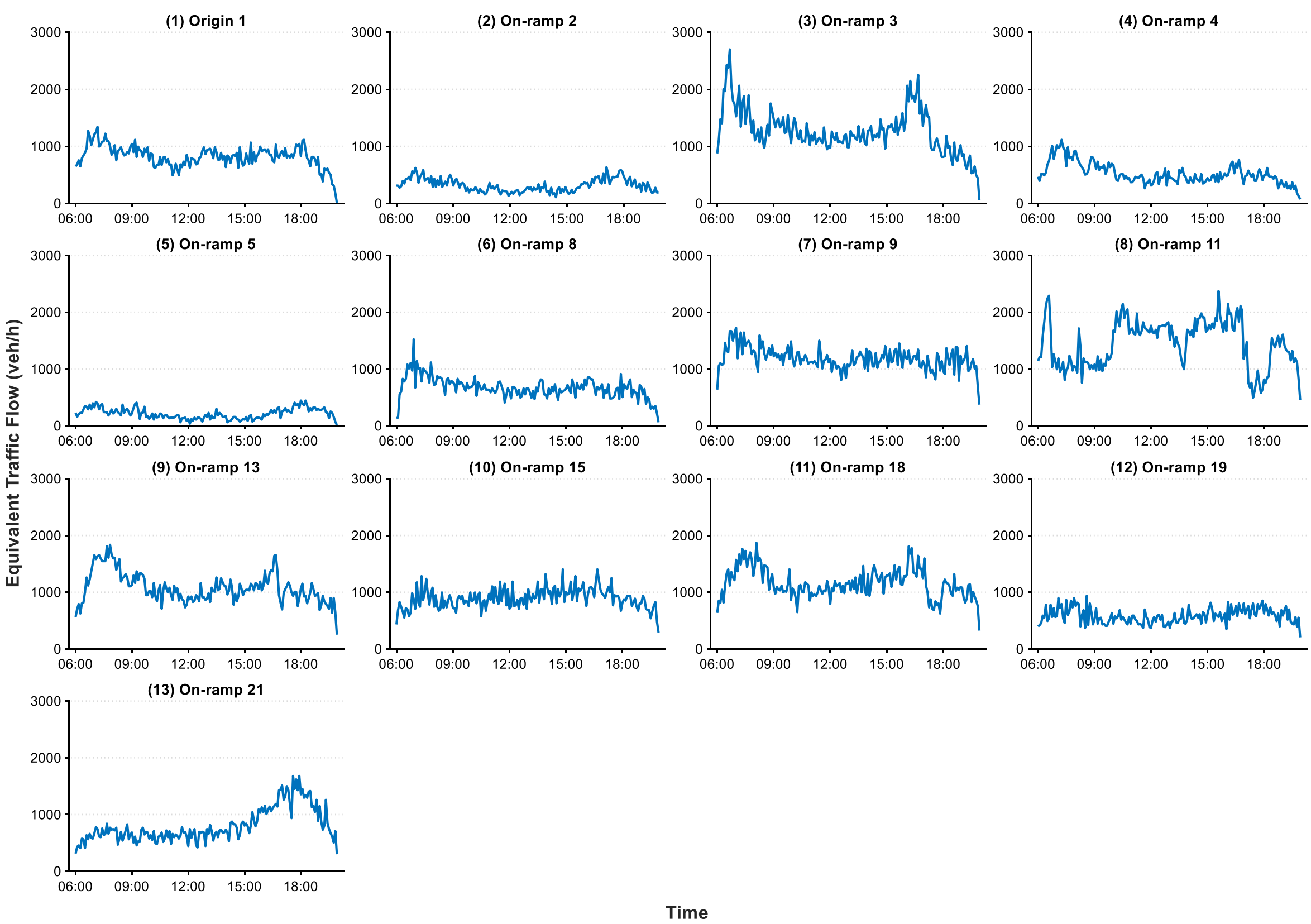

Fig. 12 Distribution of inflow at each entrance over time

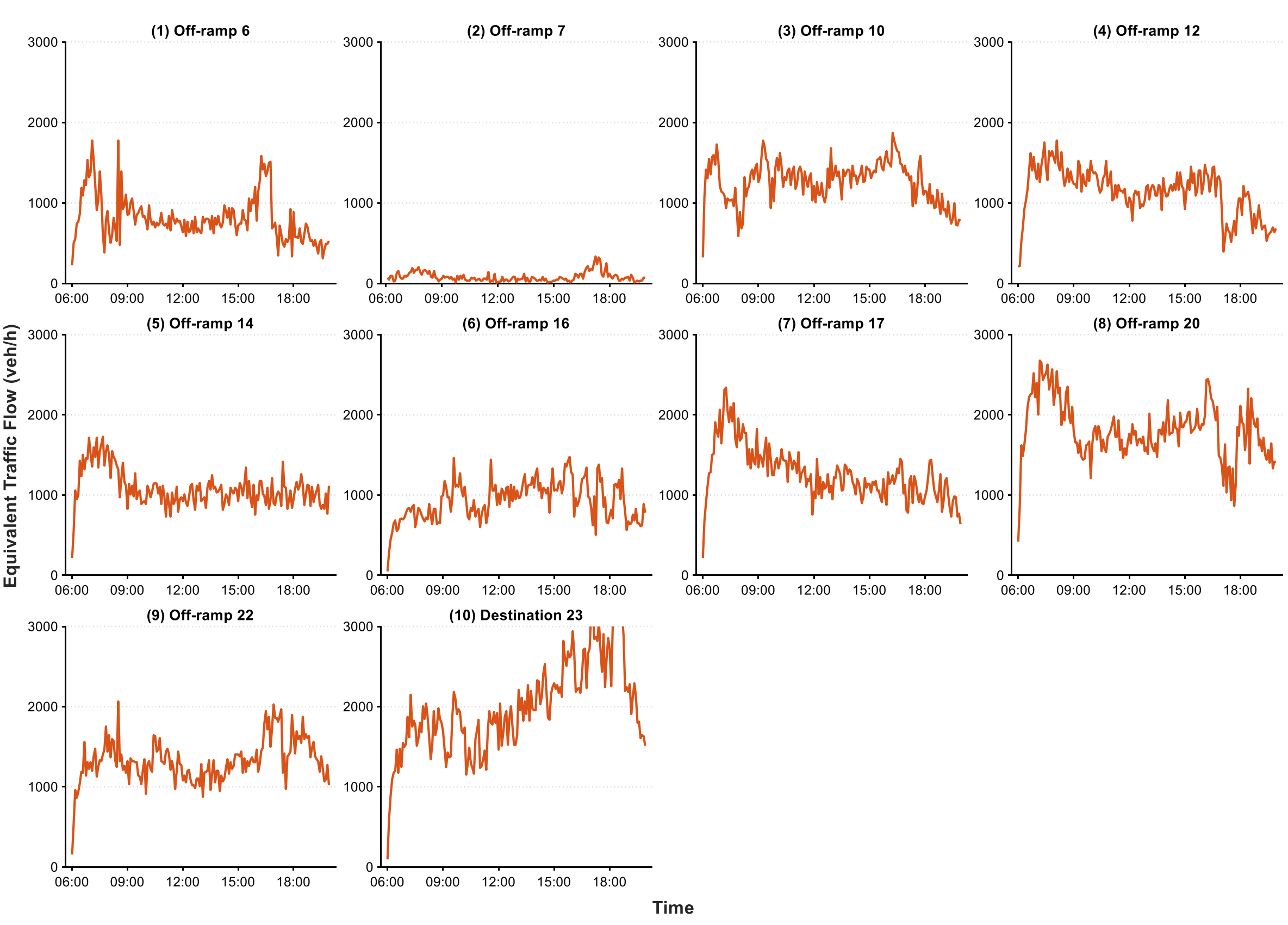

Fig. 13 Distribution of flow at each exit over time

Table 5 Extraction of key feature variables

| | WS 1-M | WS 1-D | WS 2-M | WS 2-D | WS 3-M | WS 3-D | WS 4-M | WS 4-D | WS 5-M | WS 5-D | WS 6-M | WS 6-D | WS 7-M | WS 7-D |
|---|---|---|---|---|---|---|---|---|---|---|---|---|---|---|
| V1 | D11 | D21 | S32 | D12 | D21 | D21 | D21 | D21 | D21 | D31 | D22 | D21 | D31 | D31 |
| V2 | D21 | D22 | D11 | S21 | D11 | D32 | D11 | S43 | F21 | D23 | D11 | D11 | M32 | D21 |
| V3 | D32 | D12 | D32 | S23 | D31 | D31 | S11 | S11 | D22 | D41 | D41 | D41 | S31 | D22 |
| V4 | D22 | S32 | F22 | O32 | D32 | D11 | S31 | D43 | D43 | D42 | D12 | S41 | D21 | S51 |
| V5 | D12 | F22 | F11 | S22 | D13 | D12 | D22 | D22 | D33 | S42 | D23 | S42 | D12 | D11 |
| V6 | D31 | D11 | F43 | S13 | F31 | S12 | D43 | O33 | D31 | D11 | D42 | O12 | S51 | S41 |
| V7 | S41 | D31 | S23 | S41 | D12 | O13 | D32 | D11 | D23 | F23 | S11 | S11 | D32 | D12 |
| V8 | S12 | D13 | O32 | S12 | D22 | D23 | F11 | D41 | S53 | D51 | S23 | D12 | D23 | D33 |

Table 6 Summary of risk assessment and prediction functions for all tasks

| Task | CT | Fixed parameter | | | | | | | | | Random parameter | | |
|---|---|---|---|---|---|---|---|---|---|---|---|---|---|
| | | FP | $(\mu, \sigma)$ | $\beta\ (\bar{\beta})$ | FP | $(\mu, \sigma)$ | $\beta\ (\bar{\beta})$ | FP | $(\mu, \sigma)$ | $\beta\ (\bar{\beta})$ | RP | $(\mu, \sigma)$ | $\beta\ (\bar{\beta}, \hat{\sigma})$ |
| WS 1-M | -0.448 | D11 | (25.5, 8.9) | 1.091 | D22 | (18.8, 4.6) | -0.068 | S12 | (13.8, 2.0) | 0.236 | D12 | (25.2, 9.9) | (0.202, 1.002) |
| | | D21 | (22.5, 7.5) | 0.619 | D31 | (16.3, 3.9) | 0.352 | | | | | | |
| | | D32 | (16.3, 3.3) | 0.192 | S41 | (18.9, 1.6) | 0.107 | | | | | | |
| WS 1-D | -1.367 | D31 | (16.2, 3.9) | 0.579 | F22 | (816.7, 292.1) | 0.192 | D22 | (18.7, 4.6) | -0.301 | D21 | (22.3, 7.8) | (1.123, 0.770) |
| | | D11 | (25.2, 8.8) | 0.300 | S32 | (17.2, 1.9) | -0.043 | | | | | | |
| | | D13 | (25.3, 8.0) | 0.238 | D12 | (24.9, 9.7) | -0.102 | | | | | | |
| WS2_M | -1.673 | D11 | (14.0, 3.4) | 0.979 | F43 | (1144.9, 399.3) | 0.178 | O32 | (0.3, 0.6) | -0.1775 | F22 | (720.4, 318.1) | (0.162, 0.789) |
| | | D32 | (20.9, 6.1) | 0.524 | S32 | (17.4, 1.6) | 0.042 | | | | | | |
| | | S23 | (17.7, 1.7) | 0.229 | F11 | (478.6, 315.3) | -0.163 | | | | | | |
| WS2_D | -2.824 | S13 | (17.3, 1.9) | 0.6574 | S41 | (18.6, 1.2) | 0.175 | S12 | (16.5, 1.9) | -0.1939 | - | | |
| | | S23 | (17.8, 1.9) | 0.3781 | D12 | (18.2, 5.9) | -0.159 | S21 | (15.7, 2.0) | -1.0772 | | | |
| | | S22 | (17.0, 1.8) | 0.2201 | O32 | (0.3, 0.6) | -0.168 | | | | | | |
| WS3_M | -0.577 | D11 | (28.7, 10.0) | 1.470 | | | | | | | D22 | (23.0, 6.5) | (0.335, 0.550) |
| | | D21 | (22.4, 6.9) | 1.311 | | | | | | | D32 | (21.9, 6.6) | (0.244, 0.342) |
| | | F31 | (1342.7, 352.1) | 0.720 | | | | | | | D13 | (27.9, 10.9) | (-0.356, 1.514) |
| | | D31 | (20.7, 7.5) | -0.348 | | | | | | | D12 | (26.2, 8.2) | (-0.819, 0.008) |
| WS3_D | -1.054 | D21 | (22.2, 6.9) | 0.879 | | | | | | | D12 | (26.3, 8.1) | (0.659, 0.576) |
| | | S12 | (16.5, 2.4) | 0.366 | | | | | | | D31 | (20.5, 7.4) | (0.552, 0.597) |
| | | D11 | (28.4, 10.1) | 0.311 | | | | | | | D23 | (19.7, 6.9) | (0.052, 0.511) |
| | | O13 | (0.5, 1.0) | -0.351 | | | | | | | D32 | (21.7, 6.6) | (-1.061, 0.888) |
| WS4_M | -2.895 | F11 | (972.8, 346.4) | 1.267 | S11 | (15.0, 1.6) | -0.904 | | | | D21 | (21.5, 5.8) | (0.8501, 0.474) |
| | | D32 | (19.0, 4.8) | 0.917 | S31 | (17.3, 1.6) | -1.443 | | | | D22 | (23.1, 6.0) | (0.6799, 1.260) |
| | | D11 | (15.5, 4.8) | 0.805 | | | | | | | D43 | (17.5, 4.1) | (-0.3750, 2.826) |
| WS4_D | -1.260 | D21 | (21.6, 5.8) | 1.147 | S43 | (19.7, 1.4) | 0.203 | D22 | (23.0, 6.1) | -0.215 | S11 | (14.9, 1.6) | (0.143, 0.003) |
| | | D11 | (15.5, 4.7) | 0.679 | O33 | (0.9, 1.1) | -0.034 | | | | | | |
| | | D41 | (17.2, 4.5) | 0.345 | D43 | (17.5, 4.1) | -0.061 | | | | | | |
| WS5_M | -1.489 | D22 | (25.2, 9.2) | 1.145 | D43 | (16.3, 4.5) | 0.519 | D33 | (16.3, 5.3) | -0.363 | S53 | (18.5, 1.3) | (-0.351, 1.218) |
| | | F21 | (1319.0, 401.4) | 0.893 | D21 | (24.6, 8.6) | 0.105 | | | | | | |
| | | D31 | (12.7, 3.6) | 0.537 | D23 | (14.9, 5.2) | -0.124 | | | | | | |
| WS5_D | -2.940 | D23 | (14.9, 5.4) | 0.468 | D41 | (16.0, 4.7) | 0.284 | F23 | (536.4, 282.7) | -0.618 | - | | |
| | | D42 | (17.8, 4.8) | 0.422 | D11 | (10.1, 3.2) | -0.075 | S42 | (18.9, 1.6) | -0.658 | | | |
| | | D31 | (12.5, 3.6) | 0.306 | D51 | (15.6, 5.0) | -0.188 | | | | | | |
| WS6_M | -1.377 | D22 | (18.3, 6.3) | 1.130 | D42 | (22.7, 6.6) | 0.361 | S11 | (13.0, 2.2) | 0.075 | D41 | (19.3, 5.5) | (0.448, 1.059) |
| | | D11 | (15.4, 7.1) | 1.023 | D12 | (17.1, 6.7) | 0.146 | | | | | | |
| | | S23 | (16.1, 2.0) | 0.562 | D23 | (12.4, 4.1) | 0.082 | | | | | | |
| WS6_D | -3.264 | D21 | (15.7, 5.7) | 0.664 | S42 | (15.9, 1.8) | -0.038 | D12 | (17.1, 6.5) | -1.395 | D11 | (15.1, 7.0) | (1.351, 0.931) |
| | | D41 | (19.3, 5.5) | 0.088 | O12 | (0.1, 0.2) | -0.177 | | | | | | |
| | | S11 | (13.1, 2.3) | 0.049 | S41 | (15.4, 2.0) | -0.355 | | | | | | |
| WS7_M | -0.475 | D31 | (13.7, 5.4) | 0.727 | | | | | | | D12 | (23.1, 9.4) | (0.931, 0.851) |
| | | M32 | (19.5, 1.8) | 0.685 | | | | | | | D21 | (15.8, 6.7) | (0.849, 1.811) |
| | | D23 | (12.2, 3.4) | 0.389 | | | | | | | D32 | (18.3, 5.7) | (0.385, 0.837) |
| | | S51 | (17.9, 1.6) | -0.091 | | | | | | | S31 | (16.2, 1.9) | (-0.774, 0.347) |
| WS7_D | -1.769 | D21 | (15.5, 6.7) | 0.620 | D12 | (22.7, 9.1) | -0.080 | S51 | (18.0, 1.6) | -0.327 | - | | |
| | | D31 | (13.5, 5.3) | 0.485 | S41 | (17.2, 1.8) | -0.166 | D33 | (15.8, 4.6) | -0.682 | | | |
| | | D22 | (16.2, 4.9) | 0.403 | D11 | (14.5, 5.8) | -0.171 | | | | | | |

Note: M and D represent merging and diverging respectively. V $i$ refers to the $i^{th}$ critical variable. X $ij$ stands for factor X at the $j^{th}$ segment of the $i^{th}$ lane, where X can be D (average density), S (average speed), F (flow rate), or O (maximum speed). CT stands for Constant term, FP stands for Fixed parameter, and RP stands for Random parameter.